\documentclass[leqno,11pt]{article}

\usepackage{amsfonts,amssymb,amsgen,theorem,stmaryrd, amsopn}
\usepackage[leqno]{amsmath} 
\usepackage[T1]{fontenc}
\usepackage[latin1]{inputenc} 
\usepackage{graphicx}
\usepackage{graphpap}

 \def\pasdegrille{\let\grille =
\pasgrille}  \def\aat#1#2#3{ \divide
\dimen1 by 48 \dimen3=\dimen1 \multiply \dimen1 by #1 \advance \dimen1
by -\dimen3 \divide \dimen1 by 101 \multiply \dimen1 by 100 \divide
\dimen2 by \count11 \multiply \dimen2 by #2
\setbox0=\hbox{#3}\ht0=0pt\dp0=0pt \rlap{\kern\dimen1 \vbox
to0pt{\kern-\dimen2\box0\vss}}\dimen1= \wd1 \dimen2=\ht1}
\def\pasgrille{ \count12= \dimen1 \divide \count12 by 50 \divide
\dimen2 by \count12 \count11 =\dimen2 \ \divide \dimen1 by 48
\setlength{\unitlength}{\dimen1} \smash{\rlap{\ }} \dimen1= \wd1
\dimen2=\ht1 } \def\grille{ \count12= \dimen1 \divide \count12 by 50
\divide \dimen2 by \count12 \count11 =\dimen2 \ \divide \dimen1 by 48
\setlength{\unitlength}{\dimen1} \smash{\rlap{\graphpaper[1](0,0)(50,
\count11)}} \dimen1= \wd1 \dimen2=\ht1 }

\pasdegrille

\newtheorem{theoreme}{Theorem} \theorembodyfont{\sl}
\newtheorem{proposition}{Proposition} 
\newtheorem{lemme}[proposition]{Lemma}
\newtheorem{definition}[proposition]{Definition}
\newtheorem{rem}[proposition]{Remark} 
\newtheorem{cor}[proposition]{Corollary}
\theoremheaderfont{\sc}
\numberwithin{equation}{section}
\numberwithin{proposition}{section}

\newcommand\R{{\mathbb R}} \newcommand\N{{\mathbb N}}
\newcommand\Z{{\mathbb Z}}

\newcommand\supp{{\mathrm{supp}\hskip 0.05cm}}

\global
 \global
 \theoremstyle{break} \theorembodyfont{\it}

  \newcommand\s{\sigma}

\newcommand{\dem}{\mbox{\it Proof: }} \newcommand{\cqfd}{\mbox{ }
\hfill$\Box$}

\newcommand{\B}[3]{{B}^{{#1},{#3}}_{#2}}

\newcommand{\BL}[4]{{B}^{{#1},{#3}}_{#2}(\mathcal{L}^{#4}_t)}

\newcommand{\LB}[4]{\mathcal{L}^{#1}_t({B}^{{#2},{#4}}_{#3})}

  \newcommand{\la}{\lambda}

\def\ra{\rightarrow}

\def\e{\varepsilon} \def\cdotv{\raise 2pt\hbox{,}}

\begin{document}

 \bibliographystyle{plain} \title{\vskip -1cm On well-posedness for
   the Benjamin-Ono equation}

  \author{Nicolas Burq \footnote{Département de mathématiques, UMR
  8628 du CNRS, B\^at 425 Université Paris-Sud, F-91405 Orsay et Institut universitaire de France} \ and Fabrice Planchon \footnote{ Laboratoire Analyse, G\'eom\'etrie
  \& Applications, UMR 7539 du CNRS, Institut Galil\'ee, Universit\'e Paris
  13, 99 avenue J.B. Cl\'ement, F-93430 Villetaneuse}}
\date{}
 \maketitle
 \begin{abstract}
We prove existence of solutions for the Benjamin-Ono equation with data
in $H^s(\R)$, $s>0$. Thanks to conservation laws, this yields global solutions for $H^\frac 1
2(\R)$ data, which is the natural ``finite energy'' class. Moreover,
inconditional uniqueness is obtained in  $L^\infty_t(H^\frac 1 2(\R))$, which
includes weak solutions, while for $s>\frac 3 {20}$, uniqueness holds in
a natural space which includes the obtained solutions.
\end{abstract}
 \par \noindent
\section{Introduction}Let us consider
\begin{equation}
  \label{eq:bo}
  \partial_t u+H\partial^2_x u+u\partial_x u=0,\,\,\,\,
  u(x,t=0)=u_0(x),\,\,\, (t,x)\in \R^2.
\end{equation}
Here and hereafter, $H$ is the Hilbert transform, defined by
\begin{equation}
  \label{eq:Hilbert}
  H f(x)=\frac 1 \pi \int \frac{f(y)}{x-y}dy=\frac 1 \pi
  \mathrm{vp}\frac 1 x \star u=\mathcal{F}^{-1}(-i\mathrm{sgn}(\xi)
  \hat f(\xi)).
\end{equation}
We will restrict ourselves to real-valued $u_0$.
\begin{rem}
  One could deal with complex-valued $u_0$ at the expense of an
  additional condition: $u_0\in \dot B^{-\frac 1 2,1}_2$ would
  probably be sufficient in our setting, and this somehow encodes that
  $u_0$ has zero mean, ``taming'' the behavior of low frequencies.
\end{rem}
Equation \eqref{eq:bo} deals with wave propagation at the interface of
layers of fluids with different densities (see Benjamin~\cite{B} and Ono~\cite{O}), and it belongs
to a larger class of equation modeling this type of phenomena, some of
which are certainly more physically relevant. Mathematically, however,
\eqref{eq:bo} presents several interesting and challenging properties;
the exact balance between the degree of the nonlinearity and the
smoothing properties of the linear part preclude any hope to achieve
results through a direct fixed point procedure, be it in Kato
smoothing type of spaces or more elaborate conormal (``Bourgain'')
spaces. In fact, the flow associated to \eqref{eq:bo} fails
to be $C^\infty$ (Molinet-Saut-Tzvetkov~\cite{MST1}, and even uniformly continuous
(Koch-Tzvetkov~\cite{KT2}). By standard energy methods (ignoring therefore the
dispersive part), one may obtain local in time solutions for smooth
data, e.g. $u_0\in H^s$ with $s>\frac 3 2$ and reach $s=\frac 3 2$ by
taking in account some form of dispersion (Ponce~\cite{Ponce} and
references therein). On the other hand, \eqref{eq:bo} has global weak $L^2$
and $H^\frac 1 2$ solutions (Ginibre-Velo\cite{GVbo}) and this result relies
heavily on dispersive estimates for the nonlinear equation as well as
the following two conservation laws:
\begin{align}
  \int_\R u^2(x,t)\,dx = & \int_\R u^2_0(x)\, dx,\\
\int_\R |\sqrt{-\Delta}^{\frac 1 2} u(x,t)|^2\,dx+\frac 1 3 \int_\R
u^3(x,t)\,dx= & \int_\R |\sqrt{-\Delta}^{\frac 1 2} u_0(x)|^2\,dx+\frac 1 3 \int_\R u_0^3(x)\,dx.\label{ham}
\end{align}
Recently, progress has been achieved on the Cauchy problem
for data in Sobolev spaces, by
using more sophisticated methods: Koch and Tzvetkov~\cite{KT1} obtained $s>\frac 5 4$,
and subsequently Kenig and Koenig~\cite{KK} improved this result to $s>\frac 9 8$ (both
use Strichartz estimates which are tailored to the frequency, a procedure
directly inspired by work on quasilinear wave equations); Tao~\cite{TaoBO}
obtained $H^1$ solutions, using a (complex) variant of the
Hopf-Cole transform (which linearizes Burgers equation). These
solutions can be immediately extended to global ones thanks to another
conservation law controlling the $\dot H^1$ norm (equation 
\eqref{eq:bo}, being completely integrable, has an infinite hierarchy
of conservation laws, a fact which at the moment cannot be connected
with the Cauchy problem at low regularity). After completion of the
present work (which originally obtained $s>\frac 1 4$), we learned that Ionescu and Koenig \cite{IoKe} improved existence
all the way down to $s=0$, which yields global $L^2$ solutions. While
this obviously supersedes our result, uniqueness is meant in the class
of limits of smooth solutions; we obtain uniqueness in a natural
space (which includes the linear flow) provided $s>\frac 3 {20}$. Moreover, one
can then deduce an unconditional well-posedness result for
$L^\infty(H^\frac 1 2)$ solutions: this is the natural energy class
and it contains the aforementioned weak solutions. Our main result reads as follows:
\begin{theoreme}
For any $s>\frac 3 {20}$, there exists a unique strong solution of the
Benjamin Ono equation~\eqref{eq:bo}, which is $C^0_{\text{loc}}(
\mathbb{R}_t; H^s( \mathbb{R}_x))$. Furthermore, if $s\geq 1/2$, this
solution is global and unique in $L^\infty_{\text{loc}}( \mathbb{R}_t;
H^s( \mathbb{R}_x))$, while for $\frac 3 {20} <s<1/2$, uniqueness holds in suitable spaces (see Theorem~\ref{t1}). As a consequence, in the energy space, $H^{1/2}$, weak solutions are strong.
\end{theoreme}

In fact, in this paper we construct solutions for $s> 0$ and
prove along the way that they enjoy better estimates than the mere $L^\infty_t(H^s)$ bound mentioned above. Then using this knowledge, we are able
to either prove unconditional uniqueness if $s\geq \frac 1 2$, or
uniqueness in the construction class if $\frac 1 8 < s<\frac 1 2$. Let us
outline briefly the procedure, when $s>\frac 1 4$.
\begin{itemize}
\item We work with smooth solutions, and obtain a priori estimates in
  various spaces with low regularity. Classical procedures allow to
  pass to the limit later on, yielding solutions for low regularity
  data.
\item We perform a renormalization in the spirit of Tao \cite{TaoBO}. As
  far as we know, this trick goes back to Hayashi-Ozawa~\cite{HO} when dealing with nonlinear
  Schr\"odinger equations with derivatives: facing an operator
  $(\partial^2_x +a(x) \partial_x)\phi$, one may reduce it to $(\partial^2_x
  +(-\partial_x a(x)/2+a^2(x)/4))\psi$ through conjugation: $\psi=\exp(-\int^x
  a(y)\,dy/2)\phi$. As observed in \cite{TaoBO}, the Hilbert transform
  is nothing but multiplication by $-i$ on positive frequency; and we may
  reduce ourselves to positive frequency because the solution is
  real-valued. As such, the exponential factor will be purely
  imaginary (otherwise, one needs some form of decay at infinity to
  make it bounded after space integration) and is therefore irrelevant when dealing with Lebesgue
  norms. 
\item Rather than performing the above conjugation globally, we first
  paralinearize the equation and gauge away only the worst term, which
  is when a low-high frequencies interaction takes place with the
  derivative falling on the high frequencies. In effect, we are
  replacing the exponentiation procedure by a paraproduct with the
  exponential factor. While this creates a lot of error terms, it
  highlights clearly which are the terms one should focus on.
\item We use (a variant of) conormal spaces $X^{s,b}$; in fact, as
  remarked by Tao \cite{TaoBO},
  should one silently drop the low-high interaction in the original
  equation mentioned above, the resulting equation can be
  proved to be well-posed in $X^{0^+,\frac 1 2^+}$. However, we need
  to deal with the exponential factor coming from the gauge
  transformation, and
  this is where an $\frac 1 4 $ loss occurs in a natural way.
\item Inverting the gauge in conormal spaces will lose a $\frac 1 4$ factor
  in spatial regularity (as an interpolation between two ``crude''
  estimates, an $X^{s,0}$ one which does not lose anything and an
  $X^{s,1}$ which loses only $\frac 1 2$ thanks to a smoothing
  effect). Meanwhile, the gauge action only requires the exponential
  factor to be roughly in $X^{1,\frac 1 2}$ which matches
  exactly. On the other hand  we may still invert the gauge without loss in any
  Lebesgue-like space-time or time-space norm, which allows to retain
  our solution at $H^{s}$ level.
\item Obtaining an a priori estimate does not provide uniqueness, and
  one has to perform a separate argument. This requires taking
  differences of two solutions, and performing another gauge transform. Dissymetry leads to worst terms than before, but by
  using all the a priori knowledge one has on both solutions (and
  especially one being the limit of smooth solutions constructed
  before), we are able to close an a priori estimate in (a suitable
  version of) $X^{-\frac 1 2,\frac 1 2}$.
\item The unconditional uniqueness result in $H^s, s\geq 1/2$ is
  obtained by adapting slightly the uniqueness result for $s\frac 1 4+$. The main step is to obtain an $L^4_t(L^\infty_x)$ a
  priori bound on the solution which is then bootstrapped to get a
  bound in $X^{0, \frac 1 2}$; this turns out to be enough to handle
  the previous uniqueness argument (using that we have  $\frac 1 4 -$
  derivatives of scope).
\end{itemize}
In order to go further down, $s<\frac 1 4$, one has to split
$w$ in different parts and perform an iteration on the nonlinear
quadratic terms using this decomposition, which allows to recover the
seemingly hopeless $1/4$ loss from the gauge. We only sketch the
proof of the following result.
\begin{theoreme}
For any $s>0$, there exists a strong solution of the
Benjamin Ono equation~\eqref{eq:bo}, which is $C^0_{\text{loc}}(
\mathbb{R}_t; H^s( \mathbb{R}_x))$, and unique as a limit of smooth solutions.
\end{theoreme}
 Uniqueness for $0<s<\frac 1 4$ turns out to be much harder than for
 $s>\frac 1 4$ and we only provide an
outline of the proof in the case $s>3/20$; while obtaining existence
 for $s>\frac 1 8$ is a relatively straightforward modification of our
 main argument for $s>\frac 1 4$, the uniqueness part is quite
 involved and most likely not optimal; there remains a gap in the $0\leq s\leq \frac
 3 {20}$ range.

 Finally, we note that L. Molinet has (\cite{Mo05}) obtained (global) well-posedness for the
Benjamin-Ono equation on the torus, for $H^\frac 1 2(\mathbb{T})$
data.

\section*{Acknowledgements}
We would like to thank Jean-Marc Delort and Nikolay Tzvetkov for
numerous helpful and enlightening discussions.

\section{Statement of results}
\label{sec:resultats}
Before stating any results, we need to define several functional
spaces which will be of help. We start with (inhomogeneous) Besov
spaces (\cite{BL} for details). Let us recall that a Littlewood-Paley
decomposition is a collection of operator $(\Delta_j)_{j\in \Z}$
defined as follows: let $\phi \in \mathcal{S}(\mathbb{R}^{n})$ such that $\widehat\phi =
1$ for $|\xi|\leq 1$ and $\widehat\phi= 0$ for $|\xi|>2$,
$\phi_{j}(x)= 2^{nj}\phi(2^{j}x)$, $ S_{j} = \phi_{j}\ast\cdot$,
$\psi_j(x)=(\phi_{j+1}-\phi_j)(x)$, $\Delta_{j} = S_{j+1} -
S_{j}=\psi_j \ast \cdot$. For notational convenience, we may sometimes
refer to $S_0$ as $\Delta_{-1}$. We shall denote  by $u_j= \Delta_j u$
and $u_{\prec j} = S_{j-1}u$. Finally, we define the paraproduct
between two functions $f,g$ as
\begin{equation}
  \label{eq:bony}
  T_g f=\sum_j S_{j-1}(g) \Delta_j(f),
\end{equation}
which captures the low frequencies (of $g$) vs high frequencies (of $f$) interaction in
the product $gf$.
\begin{definition}
\label{d1}
 Let $f$ be in $\mathcal{S}'(\mathbb{R}^{n})$, $s\in \R$ and $1\leq
 p,q\leq +\infty$.  We say $f$ belongs to $B^{s,q}_{p}$ if and only if
\begin{itemize}
\item $S_0 f\in L^p$.
\item The sequence $(\varepsilon_{j})_{j\in \N}$ with $\varepsilon_j = 2^{js}\| \Delta_{j} (f)\|_{L^{p}}$
belongs to $l^{q}$.
\end{itemize}  
\end{definition}
Two modifications will be of interest, to handle the additional
time variable.
\begin{definition}
  \label{d2}
Let $u(x,t)\in \mathcal{S}'(\mathbb{R}^{n+1})$. We say
that $u\in \LB \rho s p q $ if and only if, for all $j\geq -1$,
\begin{equation}
  \label{eq:raaah12}
  2^{js}\|\Delta_j u\|_{L^\rho_t(L^{p}_x)} =\varepsilon_j \in l^q.
\end{equation}
\end{definition}
\begin{definition}
  \label{d3}
Let $u(x,t)\in \mathcal{S}'(\mathbb{R}^{n+1})$. We say
that $u\in \BL s p q \rho$ if and only if, for all $j\geq -1$,
\begin{equation}
  \label{eq:raaah12bis}
  2^{js}\|\Delta_j u\|_{L^p_x(L^{\rho}_t)} =\varepsilon_j \in l^q.
\end{equation}
\end{definition}
Finally, we define conormal spaces: set
\begin{equation}
    \label{eq:defbloc}
\Delta_{jk}^\pm v (x,t) =   \mathcal{F}_{\tau,\xi}^{-1}
  (\chi_{\pm \xi\geq 0} \psi^{\pm}_{jk}(\tau,\xi))\mathcal{F}_{t,x}( v)),\qquad \Delta_{jk} v (x,t)= \Delta_{jk}^+ v (x,t)+ \Delta_{jk}^- v (x,t)
\end{equation}
with
$$
\psi^{\pm}_{jk}(\tau,\xi)=\psi_j(\xi)\psi_k(\tau\mp
  \xi^2),
$$
which (dyadically) localizes $|\xi|$ at $2^j$ and $|\tau\mp \xi^2|$ at
$2^k$.
\begin{definition}
  Let  $u(x,t)\in \mathcal{S}'(\mathbb{R}^{n+1})$, $s,b\in \R$ and
  $1\leq q\leq +\infty$. We say
that
 $u\in X^{s,b,q}$ if and only if, for all $j\geq -1$,
\begin{equation}
  \label{eq:xsbbo}
   \|\Delta_{jk} u\|_{L^2_{t,x}} \lesssim
  2^{-js-kb} \e_{jk},\,\,\,(\e_{jk})_{jk}\in l^q.
\end{equation}
$X^{s,b,q}$ is endowed with its natural norm. An alternative definition is as follows:

\begin{equation}\label{eq.defalt}
u \in X^{s,b,q} \Leftrightarrow S_0(-t)u(t, \cdot)\in B^{b,q}_{2,t}; B^{s,q}_{2,x}
\end{equation}
where $S_0(t)= e^{-tH\partial_x^2}$ is the free evolution group. 
\end{definition}
\begin{definition}\label{def.5}
For $T>0$, we say that  $u\in X^{s,b,q}_T$ if $u\in X^{s,b,q}_{loc(t)}$ and 
\begin{equation}\label{eq.defX}
 \|u\|_{X^{s,b,q}_T}= \inf \{\|v \|_{X^{s,b,q}} \text{ with } v\in X^{s,b,q},\, (v-u) \mid_{t\in (-T,T)}=0, v\mid_{t\notin [-2T, 2T]} =0\}
\end{equation}
\end{definition}
We shall use the following result to estimate norms in $X^{s,b,q}_{T}$
\begin{lemme}\label{lem.substit}
Consider $u$ the solution of 
$$
 (\partial_t + H \partial^2_x)u=f, \qquad u \mid_{t=0}=u_0.
$$
Then for any $s$ and $0<b<1$
$$
 \|u\|_{X^{s,b,q}_T}\leq C_T (\|u_0\|_{B^{s,q}_2}+ \|f\|_{X^{s,b-1,q}_T}).
$$
\end{lemme}
\dem Take a sequence $f_n$ realizing the inf in~\eqref{eq.defX} and $u_n$ solution of 
$$
 (\partial_t + H \partial^2_x)u_n=f_n, \qquad u \mid_{t=0}=u_0.
$$
Clearly $(u_n-u)\mid_{|t|\leq T} =0$. Let us study first the
contribution to $u_n$ of the low conormal frequencies
($k=0$), i.e. let $u_n^0$ be the solution of
$$
(\partial_t + H \partial^2_x)u_n^0=\sum_j \Delta_{0,j}f_n= f_n^0,
\qquad u \mid_{t=0}=u_0.
$$
 Then
\begin{multline}
 \|\Delta_j u_n^1\|_{X^{s,0,2}_T}\leq C\|\Delta_j u_n\|_{L^2(-2T, 2T); H^s}\leq C(\|\Delta_ju_0\|_{H^s}+ \|\Delta_jf_n^0\|_{L^1(-2T, 2T); H^s})\\
\leq C(\|\Delta_ju_0\|_{H^s}+ \|\Delta_j f_n^0\|_{L^2(-2T, 2T);
H^s})\leq C(\|\Delta_ju_0\|_{H^s}+ \|\Delta_j f_n^0\|_{X^{s,-1,2}})
\end{multline}
As a consequence
$$
 \|u_n^1\|_{X^{s,b,q}}\leq C(\|u_0\|_{{B^s,q}_2}+\|
f_n^0\|_{X^{s,b-1,q}}).
$$
 Let us now study the contributions of the
high conormal frequencies. Let $u_n^1= Tf_n^1$ be the solution of
$$(\partial_t + H \partial^2_x)u_n^1=\sum_{k\neq 0,j} \Delta_{k,j}f_n=
f_n^0, \qquad u \mid_{t=0}=0.$$ Then the same argument as above shows
that if $g\in X^{s,0,2}$, $v=Tg$ is bounded in $X^{s,0,2}$ and as a
consequence, for $\phi\in C^\infty_0 (\mathbb{R})$ equal to $1$ on
$[-T,2T]$, $w= \phi(t) v$ being a solution of
$$ (\partial_t + H \partial^2_x)w=g+ \phi'(t) v$$ satisfies
$$ i(\tau - |\xi|\xi)\widehat{w}(\tau, \xi)=\widehat{g+ \phi'(t) v}(
\tau, \xi)$$ and
$$ \|w\|_{X^{s,1,2}} \leq C (\|g\|_{X^{s,0,2}}+ \|\phi'(t)
v\|_{X^{s,0,2}})\leq C (\|g\|_{X^{s,0,2}})$$ This shows that $T$ is
continuous from $X^{s,0,2}$ to $X^{s,1,2}$. By duality (since $T^*$
has essentially the same form as $T$) it is also continuous from
$X^{s,-1,2}$ to $X^{s,0,2}$ and by (real) interpolation in $k$ at
fixed $j$ and $l^q$-summing the $j$, from
$X^{s,b-1,q}$ to $X^{s,b,q}$ (for $0<b<1$).  Finally we obtained
$$
 \|u\|_{X^{s,b,q}_T}\leq \|u_n\|_{X^{s,b,q}}\lesssim \|u_0\|_{{B^s,q}_2}+\|
f_n\|_{X^{s,b-1,q}}\rightarrow_{n\rightarrow +
\infty} \|u_0\|_{{B^s,q}_2}+\| f\|_{X^{s,b-1,q}}.
 $$ \cqfd

Using~\eqref{eq.defalt} and writing 
$$
u(t) = S_0(t) \frac 1 {2\pi}\int_{\tau}
e^{it\tau}\mathcal{F}_{t\rightarrow \tau} ( S_0(-t)u(t, \cdot))d\tau.
$$
 we see easily that  $X^{s,\frac 1 2, 1}$ inherits the properties of
the {\em solutions of the linear} equation $S_0(t) u_0$. We now recall
several of them which will be of interest.

\begin{proposition}
\label{foliate}
The following properties hold true:
\begin{itemize}
\item{Strichartz} 
\begin{equation}\label{eq.strich}
 \|u \|_{\mathcal{L}^4_t(B^{s,1}_{\infty})} \lesssim \|u \|_ {X^{s,\frac 1 2, 1}}.
\end{equation}
\item{Maximal function} 
\begin{equation}\label{eq.max}
 \| u \|_{B^{0,1}_4(\mathcal{L}^\infty_t)} \lesssim \|u\|_ {X^{\frac 1
 4 ,\frac 1 2, 1}}.
\end{equation}
\item{Bilinear smoothing}
\begin{equation}\label{eq.bilin}
 \|S_{j-1} u \Delta_jv \|_{L^2_{t,x} }\lesssim 2^{-\frac j 2} \|u\|_
 {X^{0 ,\frac 1 2, 1}} \|v\|_ {X^{0,\frac 1 2, 1}}.
\end{equation}
\item{Smoothing}
\begin{equation}\label{eq.smooth}
\|  u \|_{B^{s+\frac 1 2,1}_\infty(\mathcal{L}^2_t)} \lesssim  \|
u\|_{X^{s, \frac 1 2,1}}.
\end{equation}
\end{itemize}
\end{proposition}
\dem Strichartz, smoothing and maximal function estimates are all
classical in the Schr\"odinger context and extend routinely to the
Benjamin-Ono situation. The bilinear estimate is the 1D case of a
generic $L^2_{t,x}$ bilinear estimate, namely $\|\Delta_k e^{it\Delta} u_0 \Delta_j
 e^{\pm it\Delta} v_0\|_{L^2_{t,x}} \lesssim 2^{\frac{n-1} 2 k-\frac j 2}\|\Delta_k u_0\|_2
\|\Delta_j v_0\|_2 $, where $n$ is the space dimension. However one can provide
a very simple proof of the 1D case through the following identity
(which works equally well without a complex conjugation on one factor):
\begin{eqnarray*}
  \left| \int e^{-it|\xi-\eta|^2} \bar g(\eta-\xi) e^{it|\eta|^2}
  f(\eta)\,d\eta\right|^2 & = & \int e^{2it\xi(\eta-\lambda)} f(\eta)\bar
  f(\lambda) \bar g(\eta-\xi) g(\lambda-\xi)\,d\eta\, d\lambda\\
 \int \left| \int e^{-it|\xi-\eta|^2} \bar g(\eta-\xi) e^{it|\eta|^2}
  f(\eta)\,d\eta\right|^2dt & = & \int \delta(2\xi(\eta-\lambda)) f(\eta)\bar
  f(\lambda) \bar g(\eta-\xi) g(\lambda-\xi)\,d\eta\, d\lambda\\
 \int \left| (e^{-it|\cdot|^2} \bar g\ast e^{it|\cdot|^2}
  f)(\xi)\,d\eta\right|^2dt \,d\xi & =  &\int  f(\eta)\bar
  f(\lambda) \bar g(\eta) g(\lambda)\,\frac{d\eta \,d \lambda}{2|\eta-\lambda|},
\end{eqnarray*}
from which the result follows whenever $f$ and $g$ have disjoint
compact support.\cqfd

After these preliminaries, we can proceed with our main result in the
case where $s\geq \frac 1 4$.
\begin{theoreme}
\label{t1}
Let $u_0\in \B {\frac 1 4} 2 1$. Then there exists a local in
time solution $u\in C[-T,T];\B {\frac 1 4} 2 1)$ to \eqref{eq:bo},
which is unique provided that
\begin{equation}
  \label{eq:classe}
  u\in \BL {\frac 3 4} \infty 1 2\cap \BL
  0 4 1 \infty \cap
  X^{0,\frac 1 2,1}.
\end{equation}
Moreover the map $S(t):\,u_0 \rightarrow u $ is continuous from the ball
$$
 B(0,R)= \{ u_0\in B^{\frac 1 4,1}_2; \|u_0\|_{B^{\frac 1 4,1}_2}\leq
  R\}
$$
to the space
$$
 \LB \infty {\frac 1 4} 2 1 \cap \LB 4 {\frac 1 4} \infty 1 \cap\BL {\frac 3 4} \infty 1 2\cap \BL
  0 4 1 \infty \cap
  X^{0,\frac 1 2,1}.
$$
\end{theoreme}
 Furthermore, we have propagation of regularity.
\begin{theoreme}
\label{t2}
if $u_0\in H^\sigma$, $\sigma>\frac 1 4$, then $S(t) u_0\in C((-T, T); H^\sigma)$ and
$$
 \sup _{t\in (-T, T) } \|S(t) u_0\|_{H^\sigma} \leq C(\|u_0\|_{\B
 {\frac 1 4} 2 1})
\|u_0\|_{H^\sigma}.
$$
In particular, when $\sigma=\frac 1 2$, the flow map is defined for
all $t\in \R$ and we have a global bound on the $H^\frac 1 2$ norm.
\end{theoreme}
The uniqueness part of Theorem \ref{t1} relies on estimates which are
of independent interest when studying the regularity of the flow map.
\begin{theoreme}\label{th.3}
  Let $u$ and $\tilde u$ be two solutions to \eqref{eq:bo} verifying
  \eqref{eq:classe}. Then There exist $T>0$ (depending on the size in $B_2^{1/4,1}$ of $u_0, \tilde u_0$) such that we have
  \begin{align}
    \label{eq:lip}
    \|u-\tilde u\|_{C_t(\B {-\frac 1 2} 2 1)} \leq & C(u_0,\tilde u_0)
    \|u_0-\tilde u_0\|_{\B {-\frac 1 2} 2 1}\\
    \|u-\tilde u\|_{X_T^{-\frac 3 4,\frac 1 2,1}} \leq & C(u_0,\tilde u_0)
    \|u_0-\tilde u_0\|_{\B {-\frac 1 4} 2 1}
  \end{align}
which implies uniqueness; moreover, the flow map is $C^\frac 1 4$ from
$L^2$ to
$X^{0,\frac 1 2,1}$ and $C^\frac 1 3$ from $L^2$ to $C_t(L^2)$.
\end{theoreme}
We proceed with statements for the $0<s<\frac 1 8$ range.
\begin{theoreme}
\label{t1bis}
Let $u_0\in \B {s} 2 1$, $s>0$. There exists a local in
time solution $u\in C[-T,T];\B {\s} 2 1)$ to \eqref{eq:bo}, and the
flow map is continuous as in Theorem \ref{t1}. Moreover, uniqueness
holds when
\begin{equation}
  \label{eq:classebis}
  u\in \BL {\frac {13}{20}} \infty 1 2\cap \BL
  {-\frac 1 {10}} 4 1 \infty \cap
  X^{-\frac 1 {10},\frac 1 2,1}\cap X^{0,\frac 1 4,1}.
\end{equation}
\end{theoreme}
We make a few remarks before proceeding with the proofs. The
Benjamin-Ono equation \eqref{eq:bo} is invariant under rescaling: if
$u(x,t)$ is a solution,
$$
u_\la=\la u(\la x,\la^2 t)
$$
is a solution for any $\la>0$. As such, the scale-invariant Sobolev
norm is $\dot H^{-\frac 1 2}$, and consequently, we are in a
subcritical situation w.r.t. scaling. For example, $\|u_\la\|_2=\la^{\frac 1
  2}\|u\|_2$; if $u$ has lifespan $T$, $u_\la$ will have lifespan
$\la^{-2} T$. Therefore, local in time existence can be reduced to
time $O(1)$ existence for small data, provided we work with
subcritical norms, by taking $\la=\e$. Note that
$$
\|u_\e\|_{\dot B^s_2}\sim \e^{\frac 1 2+s}.
$$
In the supercritical setting, i.e. $u_0\in H^s$ with $s<-\frac 1 2$,
one may prove \eqref{eq:bo} to be ill-posed in the Hadamard sense, see \cite{BuPl05}.

\section{Gauge and conormal spaces}\label{sec.gauge}
Suppose we have a smooth solution $u$ to \eqref{eq:bo}. 
One would like to obtain an a priori estimate on $u$ in $Y^\frac 1 4 \cap X^{0,\frac 1 2,1}$,
where we set
\begin{equation}
  \label{eq:space}
  Y^s =\LB \infty {s} 1 2\cap \LB 4 {s} \infty 1 \cap \BL
  {s+\frac 1 2} \infty 1 2 \cap \BL {s-\frac 1 4} 4 1 \infty,
\end{equation}
and emphasize the embedding $X^{s,\frac 1 2,1}\hookrightarrow Y^s$
 (Proposition \ref{foliate}). Observe that we will mainly use $s=\frac
 1 4$, and that while we seek $u\in Y^{\frac 1 4}$, in the conormal
 spaces, $u$ will be only $X^{0,\frac 1 2,1}$. This $\frac 1 4$ loss
 is a side effect of a gauge transform: due to the ``bad''
term $T_u \partial_x u$ in the nonlinearity, we need to renormalize, setting
$$
w^\pm=T_{\exp{\mp \frac i 2\int^x u}}u^\pm,
$$
where the antiderivative of $u$, $\int^x u = U(t,x)$ is defined as follows :
Consider $\Psi \in C^\infty_0( \mathbb{R})$ such that $\int \Psi(y) dy =1$ and define
$$ U(t,x)= \int_y \Psi(y) \int^x_y u(t, z) dz dy + G(t)$$ with $G$ to be fixed later. Obviously 
$$\partial_x U(t,x)= \int_y \Psi(y)  u(t, x)  dy = u(t,x)$$
and 
\begin{equation}\label{eq.integrationbis}
\begin{aligned}
\partial_t U(t,x) &= \int_y \Psi(y) \int^x _y \partial_t u(t,z) dz dy  + G'(t)\\
& = \int_y \Psi(y) \int^x _y -  H \partial_z^2 u(t,z) - \frac 1 2 \partial_z (u^2(t,z)) dz dy + G'(t)\\
& = - H \partial_x u(t,x) - \frac 1 2 u^2(t,x) + \int_y \left(H\Psi'(y)  - \frac 1 2 \Psi(y)\right) u(t,y) dy + G'(t).
\end{aligned}
\end{equation}
Now we choose 
$$ G(t) = \int_{s=0}^t \int_y \left(-H\Psi'(y)  + \frac 1 2 \Psi(y)\right) u(s,y) dy ds$$
so that 
\begin{equation}\label{eq.integration}
 \partial_t U(t,x) =  - H \partial_x u(t,x) - \frac 1 2 u^2(t,x).
\end{equation}
Remark that the construction of the anti derivative of $u$ makes sense
for $u \in L^2_{x,t \text{loc}}$.

One has to deal with the imaginary exponential and the $X^{s,\theta,1}$
spaces. There will be several terms which involve both, and a typical
one would be (where $\phi^+=P^+\phi$ is the projection on $\{\xi>0\}$)
$$
T_{\exp(i\int^x u)} \phi^+,
$$
where $\phi$ may be the nonlinear term in the new unknown, say
$\phi^+=P^{+}(\partial_x v^{-} v^{+})$, or simply $\phi \in X^{s,b,1}$. As such, one has to consider the
following situation:
$$
\phi \in X^{s,\pm\frac 1 2,1}_+ \text{ and } u\in X_T^{0,\frac 1
  2,1}.
$$
The key idea is that we may be able to perform this (para)-product at the cost of a $\frac 1 4$ spatial
 derivative. Hence, if one is looking for the gauged function $v\in X_T^{\frac 1 4,\frac 1 2,1}$,
 inverting the gauge transform yields $u\in X_T^{0,\frac 1 2,1}$, and closing
 a nonlinear estimate in $v$ requires the estimate which follows.
\begin{proposition}\label{prop3.4}
Let $u\in X^{0, \frac 1 2,1}_T$ be a solution of~\eqref{eq:bo}. Denote by $K$ the
 operator of para-multiplication by $e^{i\int^xu}$, and let $0<b<1$. Then $K$ maps
 $X^{s, \pm b,1}_{+,T}$ to $X^{s-\frac b 2, \pm b,1
 }_{+,T}$ with norm smaller than 
$$
 C (1+\|u\|_{X^{0, \frac 1 2,1}_T}).
$$   
\end{proposition}
\begin{rem}\label{rem.2} As will be clear from the proof, this
  result above extends to the case where $u$ is a finite sum of solutions
  of~\eqref{eq:bo}. The result extends as well to the composition
  of several such paraproducts (i.e. we lose $b/2$ of a derivatives
  once and for all) or more generally to operators of the following kind
$$
 w \mapsto \sum_j S_{j-2} ( e^{i\int^x u_1}) S_{j-2}( e^{i\int^x
 u_2}) \Delta_j w
$$
with a bound 
$$
 C (1+\|u_1\|_{X^{0, \frac 1 2,1}_T})(1+\|u_2\|_{X^{0, \frac 1
 2,1}_T})
$$
where $u_1$ $u_2$ are two solutions of~\eqref{eq:bo}.
 This fact will be of importance in the proof of uniqueness (see section~\ref{sec.unique}).
\end{rem}
We first remark that the para product preserves the $\xi$ localization. As a consequence, by interpolation the result reduces to proving the
following two cases:
\begin{itemize}
\item First, $K$ maps $X^{s,0,q}_{+,T}$ to $X^{ s,0,q}_{+,T}$ given the
  exponential factor is imaginary, hence bounded (notice we do not use
  any other information in this case: the continuity constant is
  $O(1)$).
\item Second,  $K$ maps $X^{s,1,2}_{+,T}$ to $X^{ s,1,2}_{+,T}$, with
  $q=1,\infty$. 
\end{itemize}
Indeed, if these two points are proven, then, to recover all cases $1 \leq q\leq + \infty$, it suffices to use the classical interpolation result (see~\cite[Theorem 5.6.1]{BL})
\begin{lemme} Let $A$ be a Banach spaces and $l^s_q(A)$ be the space of sequences $(a_n)_{n\in \mathbb{N}}$ such that 
$$ \|(a_n)\|_{l^s_q(A)}= \left(\sum_n (2^{sn}\|a_n\|_{A})q\right)^{1/q} <+ \infty
$$
Then for any $s_0 \neq s_1$, $0<q_{0,1,2}\leq + \infty$, $s= \theta s_0 + (1- \theta) s_1$ the real interpolation spaces satisfy 
$$ (l^{s_0}_{q_0}(A), l^{s_1}_{q_1}(A))_{\theta, q_2} = l^s_{q_2}(A)$$ 
\end{lemme}
In fact we frequently use the following (easier)
\begin{lemme}\label{lem.interpo} Let $l^s_q$ be the space of complex valued sequences $(a_n)_{n\in \mathbb{N}}$ such that 
$$ |(a_n)|_{l^s_q(A)}= \left(\sum_n (2^{sn}|a_n|)q\right)^{1/q} <+ \infty
$$
Then for any $s_0 <s= \theta s_0 + (1-\theta)s_1<s_1$, and any sequence $(a_n)$ in $l^{s_0}_{\infty}\cap  l^{s_1}_{\infty}$, 
$$\|(a_n)\|_{l^s_1}\leq \|(a_n)\|_{l^{s_0}_\infty}^\theta \|(a_n)\|_{l^{s_1}_{\infty}}^{(1-\theta)}$$ 
(we exchange $l^\infty$ bound for $l^1$ summability). 
\end{lemme}
\dem We have 
$$|a_n| \leq 2^{ns_j}\|(a_n)\|_{l^{s_j}_{\infty}}$$
and to estimate $\sum_n |a_n 2^{ns}$, we use the $s_0$ bound to bound
the sum for $n\leq N$ and the $s_1$ bound for the sum for $n>N$ and
optimize on $N$.\cqfd

The second point is the most difficult one and we proceed with it. In fact, we only deal with the case $q=2$ and recover the other cases
by (yet another) interpolation. Notice also that due to the para-product structure,
the spatial regularity is irrelevant (and we shall fix it to
$s=0$). Consequently, it suffices to  work at fixed $j$. Finally, we first prove a global (i.e. without the index $T$) version of the estimate (assuming that the functions are compactly supported in $(-2T, 2T)$), the local estimate follows by a limiting procedure. We consider $
S_{j-1}(e^{i\int^x u})\Delta_j(\phi)$ and apply the Schr\"odinger
operator,
$$
(i\partial_t-\partial_x^2)(S_{j-1}(e^{i\int^x u}) \Delta_j(\phi))=F_1+F_2,
$$
with
$$
F_1=(i\partial_t-\partial_x^2)(S_{j-1}(e^{i\int^x  u}))\Delta_j(\phi)
+S_{j-1}({e^{i\int^x  u}}) (i\partial_t-\partial_x^2) \Delta_j(\phi),
$$
and
$$
F_2=-2\partial_x S_{j-1}(e^{i\int^x  u})\partial_x \Delta_j(\phi)
=-2iS_{j-1}( e^{i\int^x  u}u )\Delta_j(\partial_x \phi) .
$$
We estimate all these terms in $L^2$, knowing $\phi\in X^{0,1,2}$. The
last term in $F_1$ is trivially ok: we control it by the norm of
$(i\partial_t -\partial^2_x)\Delta_j(\phi)$ in $L^2$. To
deal with the first
term, one simply recall the definition of $\int^x u$ and due
to~\eqref{eq.integration}, we can replace $\partial_t \int^x u$ by $-
H \partial_x u -\frac 1 2 u^2$. As a consequence, the first
term in $F_1$ is equal to
$$
 T_{e^{i\int^x u}(\frac {3 u^2} 2 + (H-i) \partial_x u)} \phi.
 $$
 This procedure yields two terms: the first one, namely $
 S_{j-1}(u^2 e^{i\int^x u})\Delta_j (\phi)$, is cubic, and
 by Strichartz inequality~\eqref{eq.strich}, we get $u^2$ in
 $L^4_t (L^2_x)$ and $\phi $ in $L^4_t (L^\infty_x)$. Hence, the contribution
 of this term is bounded by
$$
  \|u\|_{L^\infty_t(L^2_x)}\|u\|_{L^4_t(L^\infty_x)}
  \|\Delta_j(\phi)\|_{L^4_t; L^\infty_x}\lesssim \|u\|_{X^{0,\frac 1 2
  ,1}}^2\|\Delta_j \Phi\|_{X^{0, \frac 1 2, 1}}.
$$ 
The other term is equal to $S_{j-1}((H-i)\partial_x u {e^{i\int^x u}})
 \Delta_j(\phi)$. Thus, it is essentially the same thing as $F_2$,
 except for the distribution of derivatives; given the paraproduct
 structure, it will be easier to deal with than $F_2$, and
 consequently we shall focus only on $F_2$. We have 
\begin{equation}\label{eq.paraprod}
 S_{j-1} (e^{i\int^x  u}u) \Delta_j \partial_x\phi  = 
 e^{i\int^x  u} S_{j-1}u \Delta_j \phi  + ([S_{j-1},e^{i\int^x  u}] u) \Delta_j
 \partial_x\phi.
\end{equation}

The first term is the main one, we simply use~\eqref{eq.bilin} to deal with this term and obtain
\begin{equation}
\begin{aligned}
\|e^{i\int^x  u}  S_{j-1}u \Delta_j \partial_x\phi\|_{L^2_{t,x}}&= \|
S_{j-1}u \Delta_j \partial_x\phi\|_{L^2_{t,x}}\\
&\lesssim  \|u\|_{X^{0, \frac 1 2 , 1}} 2^{-\frac j 2} \| \Delta_j
\phi\|_{X^{1, \frac 1 2 , 1}}\\
&\lesssim  \|u\|_{X^{0, \frac 1 2 , 1}} 2^{\frac j 2} \| \Delta_j
\phi\|_{X^{0,\frac 1 2 , 1}}
\end{aligned}
\end{equation}
 which yields the $1/4$ loss in Proposition~\ref{prop3.4} . For
the remaining commutator, we use the following classical lemma:
\begin{lemme}
\label{l1}
  Let $f\in L^p$, $\nabla g \in L^\infty$, then
  \begin{equation}
  \label{eq:17}
\|   [S_j,g]f\|_p \lesssim 2^{-j} \|\nabla
g\|_\infty \|f\|_p.
\end{equation}
\end{lemme}
We provide the (trivial) proof for sake of completeness:
given $S_j f=2^{nj}\phi(2^j \cdot)\star f$, one writes
\begin{eqnarray*}
 [S_j,g]f(x) & = &\int 2^{nj}\phi(2^j(x-y))(g(y)-g(x))f(y)dy\\
  | [S_j,g]f(x)|&\leq &\int 2^{nj}|\phi|(2^j(x-y)) |x-y|\|\nabla
    g\|_\infty |f|(y)dy,\\
 & \leq & 2^{-j} \|\nabla g\|_\infty \int  2^{nj} \theta (2^j(x-y))
    |f|(y)dy,\\
 \| [\Delta_j,g]f\|_p & \leq & 2^{-j} \|\nabla g\|_\infty \|\theta\|_1 \|f\|_p,
\end{eqnarray*}
since $\theta(x)=|x||\phi|(x)\in L^1$. \cqfd

Hence, this term is treated using Strichartz estimate~\eqref{eq.strich}, exactly as
the cubic terms above, namely
\begin{eqnarray*}
   \|([S_{j-1},e^{i\int^x  u}] u) \Delta_j
   \partial_x\phi\|_{L^2_{t,x}} \lesssim & 2^{-j} \| u\, e^{i\int^x
   u}\|_{L^4_t (L^\infty_x)} \|u\|_{L^\infty_t (L^2_x)}\, 2^j \|\Delta_j
   \phi\|_{L^4_t (L^\infty_x)}\\
 \lesssim & \|u\|^2_{X^{0,\frac 1 2,1}} \|\Delta_j \phi\|_{X^{0,\frac 1 2,1}}.
\end{eqnarray*}
 Collecting all the estimates yields
 \begin{eqnarray*}
 \|(i\partial_t- \partial_x^2) \left(S_{j-1}(e^{i\int^xu})
\Delta_j(\phi)\right)\|_{L^2_t,x} & \lesssim (1+ \|u\|_{X^{0, \frac 1 2, 1}}^2)
 2^{\frac j 2} \|\Delta_j \phi\|_{X^{0, 1, 2}},
 \end{eqnarray*}
given the embedding $X^{0,1,2}\hookrightarrow X^{0,\frac 1 2,1}$ (at
fixed $j$). Therefore,
$$
\| S_{j-1}(e^{i\int^xu}) \Delta_j(\phi)\|_{X^{0,1,2}}\lesssim 2^{\frac
  j 2} (1+
\|u\|_{X^{0, \frac 1 2, 1}}^2) \|\Delta_j \phi\|_{X^{0, 1,2}}.
$$
 On the other hand, we clearly have
$$
 \| S_{j-1}(e^{i\int^xu}) \Delta_j(\phi)\|_{L^2_{t,x}}= \|
S_{j-1}(e^{i\int^xu}) \Delta_j(\phi)\|_{X^{0,0,2}}\lesssim \|\Delta_j
\phi\|_{X^{0, 0, 2}}.
$$ 
We can now decompose all terms according to the
conormal scale ($k$) and by real interpolation on these $k$ sequences
get
$$
 \| S_{j-1}(e^{i\int^xu}) \Delta_j(\phi)\|_{X^{0,\frac 1 2,q}}
\lesssim 2^{\frac j 4} (1+
\|u\|_{X^{0, \frac 1 2, 1}}) \|\Delta_j \phi\|_{X^{0, \frac 1 2 , q}},
$$
which, after summing in $j$, is the desired result: namely, for any $1\leq
q\leq +\infty$, $K$ maps $X^{0,\frac 1 2,q}$ to $X^{-\frac 1 4,\frac 1
  2,q}$. We now have to show that this result still hold with the
local in time Bourgain spaces $X_T^{0,\frac 1 2,q}$ to $X_T^{-\frac 1
  4,\frac 1 2,q}$. For this we come back to \eqref{eq.3.3} and replace
all occurrences in the right hand side of $w$ and $u$ respectively by
$u^n$ and $w^n$ where these sequences are minimizing sequences for
Definition~\ref{def.5}. Now we define $\tilde u_j^n$ to be equal to
the (new) left hand side. Obviously, $\tilde u^n = \sum_j \tilde
u^n_j$ is equal to $u$ on $[0,T]$ and applying the (global)  estimate
we just got yields
 \begin{multline}
   \label{eq:invertgauge}
    \|\tilde u^n\|_{X^{0,\frac 1 2, 1 }} \lesssim (1+ \|u^n\|_{X^{0,\frac 1
    2, 1 }})\|w^n\|_{X^{\frac 1 4, \frac 1 2, 1}} + \|u^n\|_{X^{0,\frac 1
      2,1}}^2+\|u_0\|_{L^2_x}\|u^n\|_{\LB 4 {\frac 1 4} \infty
    1}\\
+ \|u_0\|_{L^2_x}( \|u^n\|_{X^{0,\frac 1
    2,1}}+1)\|u^n\|_{X^{0,\frac 1 2,1}}.\text{\hskip 1cm {\cqfd}}
 \end{multline}
We shall also need to ``invert'' the estimates in
Proposition~\ref{prop3.4}:
\begin{proposition}\label{prop3.5}
Let $u\in X^{0, \frac 1 2,1}_T$ be a solution of~\eqref{eq:bo}. Let
$w$ be defined (as a real valued function) by
$$
 w= w^+ + w^-, w^\pm = P^\pm S_0 u+\sum_{j\geq 0} w_j^\pm, \qquad w_j^\pm= S_{j-1}(e^{\mp\frac i 2 \int^x u})
\Delta^+_j u .
$$
 Then we have
 \begin{eqnarray}
   \label{eq:invertgaugebis}
    \|u\|_{X^{0,\frac 1 2, 1 }} \lesssim &
 (1+ \|u\|_{X^{0,\frac 1
    2, 1 }})\|w\|_{X^{\frac 1 4, \frac 1 2, 1}} + \|u\|_{X^{0,\frac 1
      2,1}}^2+\|u_0\|_{L^2_x}\|u\|_{\LB 4 {\frac 1 4} \infty
    1}\nonumber\\
 & {}+ \|u_0\|_{L^2_x}( \|u\|_{X^{0,\frac 1
    2,1}}+1)\|u\|_{X^{0,\frac 1 2,1}}.
 \end{eqnarray}
\end{proposition}
The low frequency part is a trivial issue. As before, let us first prove the estimate in global spaces. Let us focus on the high
frequencies : from now on, denote by $F=\exp(i\int^x u)$. We write 
\begin{equation}
\begin{aligned}
 w_j^+ &= S_{j-1} F \Delta_j u^+\\
&= F \Delta_j u^+ +\sum_{k\geq j-1} \Delta_k ( F) \Delta_j u^+
\end{aligned}
\end{equation}
As a consequence
\begin{equation}\label{eq.3.3}
\begin{aligned}
 \Delta_j u^+  &= F^{ -1} w^+_j-  F^{-1} \sum_{k\geq j-1} \Delta_k (F) \Delta_j u^+\\
&= S_{j-3} (F^{-1}) w_j^+ + \sum_{k\geq j-3} \Delta_k( F^{-1}) 
w^+_j- F^{-1} \sum_{k\geq j-1} \Delta_k (F) \Delta_j u^+\\
&= S_{j-3} (F^{-1}) w_j^++ \tilde \Delta_j \left( \sum_{k\sim j} \Delta_k( F^{-1}) 
w^+_j\right)- \tilde \Delta_j \left(F^{-1} \sum_{k\geq j-1} \Delta_k (F) \Delta_j u^+\right),
\end{aligned}
\end{equation}
due to the frequency localization of $w_j^+$, with $\tilde \Delta_j$
an enlargement of $\Delta_j$; the first term is ok
according to Proposition~\ref{prop3.4} (where, obviously, the sign of
the phase term is irrelevant). We now proceed with the remaining terms,
and need to estimate
$$
 \|\Delta_j u^+-S_{j-3}( F^{-1}) \Delta_j w^+)\|_{X^{0,\frac 1 2
 ,1}}.
$$
For lack of a better alternative, we proceed as in the proof of
Proposition~\ref{prop3.4}. On one hand we have (recall $L^2_{t,x}=X^{0,0,2}$)
 using~\eqref{eq.strich}, $T\lesssim 1$, and $v$ being either $u$ or $w$,
\begin{multline}
\|\Delta_k F^{\pm 1} \Delta_j v^+\|_{X^{0,0,2}}\lesssim
2^{-k} \| \Delta_k \partial_x  ( e^{\mp i \int^x u})\|_{L^\infty_t (L^2_x)} \|\Delta_j v^+\|_{L^2_T(L^\infty_x)}\\
\lesssim 2^{-k} \|u\, e^{\mp i \int^x u}\|_{L^\infty_t(L^2_x)}\|\Delta_j
v^+\|_{L^4_t( L^\infty_x)}  \leq 2^{- k}\|u_0\|_{L^2}\|
\Delta_j v\|_{X^{0, \frac 1 2 ,1}}
\end{multline}
 and consequently (notice the spatial regularity gain !)
\begin{equation}
\label{reverseg0}
\|\Delta_j u^+-S_{j-3}( F^{-1}) \Delta_j w^+)\|_{X^{1,0,2}}\lesssim 
\|u_0\|_{L^2}(\| \Delta_j w\|_{X^{0, \frac 1 2 ,1}}+\| \Delta_j u\|_{X^{0, \frac 1 2 ,1}}).
\end{equation}
 On the other hand we compute
 \begin{equation}
   \label{eq:youyou}
   \|(\partial_t + H \partial_x^2)(\Delta_j u^+-S_{j-3}( F^{-1}) \Delta_j w^+)\|_{X^{0,0 ,2}},
 \end{equation}
for which we proceed differently: given we are on $\xi>0$ as well, we
may use the equation for the first term and a computation similar to
the proof of Proposition \ref{prop3.4} for the paraproduct. Obviously,
we have
$$
\|(\partial_t + H \partial_x^2)\Delta_j u^+\|_{L^2_{t,x}}\lesssim 2^j
\|\Delta_j^+(u^2)\|_{L^2_{t,x}},
$$
and using a paraproduct decomposition for $u^2$ and $T\lesssim 1$,
$$
\|(\partial_t + H \partial_x^2)\Delta_j u^+\|_{L^2_{t,x}}\lesssim 2^j
(2^{-\frac j 2}\|u\|^2_{X^{0,\frac 1 2,1}}+2^{-\frac j 4}\|u\|_{\LB 4
  {\frac 1 4} \infty 1}\|u_0\|_{L^2_x}).
$$
\begin{rem}
  One could deal with the remainder term differently and estimate it
  only with $X^{0,\frac 1 2,1}$ norms: due to support conditions, only
  opposite frequencies interactions occur, for which one may use the
  smoothing effect as on the paraproduct terms.
\end{rem}
The paraproduct term in \eqref{eq:youyou} requires distributing the $i\partial_t
-\partial^2_x$ operator:
\begin{itemize}
\item first, on the exponential factor, using the equation on $u$ and
  $F^{-1}\partial_x u=\partial_x(u F^{-1})-u\partial_x F^{-1}$,
\begin{multline*}
\|S_{j-3}\left( (i\partial_t - \partial_x^2) F^{-1}\right)  \Delta_j w^+\|_{L^2_{t,x}}\\
\begin{aligned}
\lesssim & \|S_{j-3}\left( F^{-1} u^2 \right)  \Delta_j
w^+\|_{L^2_{t,x}}+ \|S_{j-3}\partial_x\left( F^{-1} u \right)  \Delta_j w^+\|_{L^2_{t,x}}\\
 \lesssim & \|u_0\|_{L^2} \| u\|_{L^4_t( L_x^\infty)}\|\Delta_j
w\|_{L^4_t( L^\infty_x)}+ 2^j \| u\|_{L^\infty_t(L^2_x)}\|\Delta_j
 w^+\|_{L^4_t(L^\infty_x)}\\
\lesssim & \|u_0\|_{L^2}\|u\|_{X^{0, \frac 1 2 ,1}}\|w_j\|_{X^{0, \frac 1 2
    ,1}}+2^j \|u_0\|_{L^2}\|w_j\|_{X^{0, \frac 1 2 ,1}}.
\end{aligned}
\end{multline*}
\item For the next term, we use the (algebraic) computation \eqref{eq:wj} from the
  next section to remark that $w^+$ satisfies an equation which
  is no worse than $u$, hence, discarding the exponential factor,
\begin{equation*}
\begin{aligned}
\|S_{j-3} (F^{-1})  \Delta_j(\partial_t + H \partial_x^2)w^+\|_{L^2_{t,x}}
 \lesssim & \|\Delta_j\partial_x( u^2)\|_{L^2_{t,x}}
+\| S_{j-1}\partial_x (F u) u^+_j\|_{L^2_{t,x}}\\
& {}+\| S_{j-1} (F u^2) u^+_j\|_{L^2_{t,x}}+\| S_{j-1} (F u) \partial_x u^+_j\|_{L^2_{t,x}}\\
 \lesssim &  2^j
(2^{-\frac j 2}\|u\|^2_{X^{0,\frac 1 2,1}}+2^{-\frac j 4}\|u\|_{\LB 4
 {\frac 1 4} \infty 1}\|u_0\|_{L^2_x})\\
 & {}+2^j \|u\|_{L^\infty_t(L^2_x)}\|u^+_j\|_{L^4_T(L^\infty_x)}\\
 &{}+
 \|u\|_{L^\infty_t(L^2_x)}\|u\|_{L^4_t(L^\infty_x)}\|u^+_j\|_{L^4_t(L^\infty_x)}\\
 \lesssim &  2^j
(2^{-\frac j 2}\|u\|^2_{X^{0,\frac 1 2,1}}+2^{-\frac j
 4}\|u\|_{\LB 4 {\frac 1 4}\infty 1}\|u_0\|_{L^2_x})\\
 & {}+ \|u_0\|_{L^2}(2^j +\|u\|_{X^{0,\frac 1 2,1}})\|u_j\|_{X^{0,\frac 1 2,1}}.
\end{aligned}
\end{equation*}
\item Finally, the last term comes from distributing the laplacian,
\begin{equation}
\begin{aligned}
\|S_{j-3}\left(\partial_x F^{-1}\right)  \Delta_j\partial_x w^+_j\|_{L^2_{t,x}}
\lesssim & 2^j \|S_{j-3}\left( F^{-1} u \right)\|_{L^\infty_t(L^2_x)}
\|  \Delta_j w^+\|_{L^2_{T}(L^\infty_x)} \\
\lesssim &   2^j  \|u_0\|_{L^2}\|\Delta_j w\|_{X^{0, \frac 1 2 ,1}}.
\end{aligned}
\end{equation}
\end{itemize}
Collecting all terms 
\begin{eqnarray}
\label{reverseg1}
\|\Delta_j u^+- S_{j-3}\bigl ( F^{-1}\bigr )  w^+_j\|_{X^{0,1,2}}
\lesssim &  2^j \left ( 2^{-\frac j 4} \bigl ( \|u\|_{X^{0,\frac 1
      2,1}}^2+\|u_0\|_{L^2_x}\|u\|_{\LB 4 {\frac 1 4} \infty
    1}\bigr )\right .\\
 & {}+\left . \|u_0\|( \|u\|_{X^{0,\frac 1
    2,1}}+1)(\|w_j\|_{X^{0,\frac 1 2,1}}+\|u_j\|_{X^{0,\frac 1 2,1}})\right).\nonumber
\end{eqnarray}
By H\"older for $k$-sequences,
$$
 \|f\|_{X^{0, \frac 1 2, 1}}\leq \left( \|f\|_{X^{0, 0,
 2}}\|f\|_{X^{0, 1, 2}}\right)^{1/2}
$$
and we obtain from \eqref{reverseg0} and \eqref{reverseg1}
\begin{eqnarray*}
   \|\Delta_j u^+- S_{j-3}( F^{-1}) \Delta_j w^+\|_{X^{0,\frac 1 2
 ,1}}
\lesssim & 2^{-\frac j 4} \bigl ( \|u\|_{X^{0,\frac 1
      2,1}}^2+\|u_0\|_{L^2_x}\|u\|_{\LB 4 {\frac 1 4} \infty
    1}\bigr )\\
 & {}+ \|u_0\|( \|u\|_{X^{0,\frac 1
    2,1}}+1)(\|w_j\|_{X^{0,\frac 1 2,1}}+\|u_j\|_{X^{0,\frac 1 2,1}}).
\end{eqnarray*}
Summing over $j$ provides our estimate.\cqfd

Next we can obtain similar results for any $L^p,L^q$ norm, but without
any spatial regularity loss: 
\begin{proposition}\label{prop.4}
For any mixed $L^pL^q$ norm (independently of the order $(t,x)$ or
$(x,t)$), we have
\begin{equation}
  \label{eq:equivLpLq}
 (1-\|u_0\|_2) \|u_j\|_{L^pL^q} \leq \|w_j\|_{L^pL^q}\leq \|u_j\|_{L^pL^q},
\end{equation}
and consequently
\begin{equation}
  \label{eq:equivY}
 (1-\|u_0\|_2) \|u\|_{Y^s} \leq \|w\|_{Y^s}\leq \|u\|_{Y^s}.
\end{equation}
\end{proposition}
\begin{rem}\label{rem.4} As the proof below shows, the result still hold if the renormalization is performed using another solution of ~\eqref{eq:bo} (or a sum of such solutions)
\end{rem}
\dem Recall
$$
w^+_j=F_{\prec j} u^+_j,
$$
$$
w^+_j= Fu^+_j+\sum_{k>j} (\Delta_k F) u^+_j,
$$
$$
u^+_{j} = F^{-1} w_{j}^+-\sum_{k>j} F^{-1} (\Delta_k F) u^+_{j}.
$$
Write
$$
\|u_j\|_{L^pL^q}\leq \|w_j\|_{L^pL^q}+\sum_{k>j}\| \Delta_k(F)
u_j\|_{L^pL^q},
$$
together with
$$
\|\Delta_k F\|_{L^\infty_{t,x}}\lesssim 2^{-k} \|\Delta_k \partial_x
F\|_\infty\lesssim 2^{-k}\|\Delta_k( u F)\|_\infty\lesssim  2^{-\frac
  k 2}\|u \|_{L^\infty_t(L^2_x)}\lesssim 2^{-\frac k 2} \|u_0\|_{L^2_x},
$$
we obtain the desired
control.\cqfd

\begin{rem}
The $1/4$ loss in the gauge transformation, when $b=\frac 1 2$, is responsible for the
$s>1/4$ assumption in our main Theorem. Further computations suggest
that this loss in the gauge transform is unavoidable with the rather
crude method we developed here. The improvements over $s=\frac 1 4$
utilize the lesser $\frac 1 8$ loss which occurs when dealing with the
$b=\frac 1 4$ case.
\end{rem}

\section{Existence}
\label{sec.key}\label{sec.4}
We now come back to considering a smooth  solution $u$ to \eqref{eq:bo}. From the renormalization estimates, we know that, provided we can estimate $w\in
X^{\frac 1 4,\frac 1 2,1}$, we will have $u\in X^{0,\frac 1 2,1}$ by
using the gauge estimate proved earlier, and moreover, $u$ in some
$L^p L^q$ space is equivalent to the same properties on $w$.

In the sequel we shall adopt the following convention.
For $w$ a (smooth) function, we shall denote by $w_{\prec j}$ any term
obtained by applying to $w$ a spectral cut-off supported in the set
$|\xi|\leq 2^{j-K}$ for a sufficiently large (but fixed) $K$. We shall
also denote by $w_{\sim j}$ any term obtained by applying to $w$ a spectral
cut-off supported in the set $2^{j-N}\leq |\xi|\leq 2^{j+N}$ for a
sufficiently large (but fixed) $N$. If the cut of is supported in the
set where $\pm \xi \geq 0$ then we will denote the result by
$u_j^\pm$. The convention will be taken that if we write $v_{\prec j}
w_j$ we have chosen $N \leq K-3$ so that this expression is still
localized in the set $|\xi|\sim 2^j$.

The main result in this section is the following.
\begin{theoreme}\label{t4}
  Let $u$ be a solution to \eqref{eq:bo} verifying
  \eqref{eq:classe}, and $w$ defined by Proposition \ref{prop3.5},
  with small $\|u_0\|_{L^2_x}$. Then we have (if $\|u_0\|_{B^{\frac 1 4,1}_2}$ is small enough)
  \begin{multline}
    \label{eq:close}
   \|w\|_{X_{T=1}^{\frac 1 4,\frac 1 2,1}} \lesssim \|w_0\|_{B^{\frac 1 4,1}_2}
+\|u\|_{X_{T=1}^{0,\frac 1 2,1}}\|w\|_{X_{T=1}^{\frac 1 4,\frac 1 2,1}}
+(1+\|u\|_{X_{T=1}^{0,\frac 1 2,1}})\|w\|^2_{X_{T=1}^{\frac 1 4,\frac 1 2,1}}
\\
+ \|u\|_{Y^{\frac 1 4}}^3 + \|u\|_{Y^{\frac 1 4}} ^2 \|w\|_{Y^{\frac 1 4}}+\|w\|^3_{X_{T=1}^{\frac 1 4,\frac 1 2,1}}.
  \end{multline}
\end{theoreme}
\dem Once again, the low frequencies ($|\xi|\lesssim 1$) is a trivial
issue: the derivative in the nonlinearity vanishes and one may use,
say, Strichartz to estimate the quadratic term. We proceed with higher
frequencies, and accordingly will sum only over $j\geq 0$. We begin by a paralinearization of the equation, starting with
\begin{equation}
\begin{aligned}
\Delta_j^+ (u^2)&=  \Delta_j^+(2 \sum_{k\sim j} 2 S_{k-1}u  \Delta_k u +
\sum_{j\lesssim k\sim k'}  u_k u_{k'})\\
&= \Delta_j^+( 2S_{j-1} u \sum_{k\sim j} \Delta_k u + \sum_{j\lesssim
  k}  u_{\sim k} ^2)
\end{aligned}
\end{equation}
where the additional terms coming from freezing $k=j$ in the $S_{k-1}$
operator are transferred to the remainder term. Now, taking further advantage of support considerations,
\begin{equation}
  \label{eq:35}
  \partial_x \Delta_j^+ (u^2)=    2u_{\prec j} \partial_x u_j^++
2[ \Delta_j^+, u_{\prec j}] \partial_x u^+_{\sim j} +2\Delta_j^+((\partial_x u_{\prec j})
u^+_{\sim j}) +  \sum_{j\lesssim k} \Delta_j^+\partial_x (u_{\sim k}
u^+_{\sim k})
\end{equation}

As a consequence, we can localize equation~\eqref{eq:bo},
\begin{eqnarray*}
 \partial_t u_j^+-i\partial^2_x u_j^++ u_{\prec j}\partial_x u_j^+ = &
-\Delta_j^+((\partial_x u_{\prec j})u^+_{\sim j})- \frac 1 2 \sum_{j\lesssim k} \Delta_j^+\partial_x (u_{\sim k}
u^+_{\sim k})\\
 &{}- [ \Delta_j^+, u_{\prec j}] \partial_x u^+_{\sim j} \\
 = & f_{j,1} + f_{j,2} + f_{j,3}= f_j^+
\end{eqnarray*}
Recall that $w$ is defined by
$$
w^{+}_j=S_{j-1}(e^{-\frac i 2 \int^x u})\Delta_j^+(u),
$$
which yields the equation on $w_j^+$:
\begin{align}
\label{eq:wj}
\partial_t w^+_j-i\partial^2_x w_j^+= & S_{j-1}(e^{-\frac i 2 \int^x u}) f^+_j+u_j^+
 S_{j-1}\bigl(( \frac 1 2 H \partial_x u + \frac 1 4 u^2) e^{-\frac i 2 \int^x u}\bigr)\\
&+\left(S_{j-1}( u e^{-\frac i 2 \int^x
 u})-S_{j-1} u S_{j-1}(e^{-\frac i 2 \int^x u}) \right)\partial_x u_j^+.\nonumber
\end{align}
The origin of the second term is clear: it comes from the linear
operator hitting the exponential factor: we use~\eqref{eq.integration} to exchange time derivatives for space derivatives. Note that the $u^2$ term is really cubic, hence it
will be easier to deal with. The last term and the
$H\partial_x u$ term are, up to commutators, $w_j\partial_x
u_{\prec j}$. We will see when dealing with $f_j$ that we also get a
term like this and a whole set of commutation terms which are supposed
to be ``better'' in that they require no conormal spaces to deal with
them.

Heuristically, when summing over $j$, we have eliminated (up to commutators !) the
worse term, namely $T_u \partial_x u$, and are left with $T_{\partial_x
  u} w$: but from the product rules in the Appendix,
$$
\|T_{\partial_x  u} w\|_{X_{T=1}^{\frac 1 4,-\frac 1 2,1}} \lesssim
\|u\|_{X_{T=1}^{0,\frac 1 2,1}} \|w\|_{X_{T=1}^{\frac 1 4,\frac 1 2,1}},
$$
which will be ok to close the
estimate.

We now proceed with estimating all right-handside terms in \eqref{eq:wj} in the space $X_{T=1}^{\frac 1 4, - \frac 1 2, 1}$. This term has the form $F(u, w, \partial_x u, \partial_xw)$. To make the exposition more clear, we shall in a first step proceed as if $w$ were in $X^{\frac 1 4,  \frac 1 2, 1}$ (the global Bourgain space) and estimate the right hand side in $X^{\frac 1 4, - \frac 1 2, 1}$ (and $u$ in $X^{0,\frac 1 2,1}\cap Y^{\frac 1 4}$). We shall give at the end of the section the modifications required to handle the argument.
\begin{itemize}
\item We first deal with the ones
  coming from the exponentiation/paraproduct, starting with
 $P_1=u^+_j S_{j-1}(e^{-\frac i 2\int^x u}H\partial_x u)$. Denoting
  (again) by $F= e^{-\frac i 2 \int^x u}$, we have
  \begin{multline}
\label{gna}
     P_1 = u_j^+ F_{\prec j} H\partial_x u_{\prec j}+u_j^+
 [S_{j-1},F_{\prec j}]H\partial_x u_{\prec j} + u_j^+ S_{j-1} ( \sum
 _{j\lesssim k} F_{\sim k} H\partial_x u_{\sim k})\\
 + u_j^+ S_{j-1}(F_{\sim j} S_{j-1} H\partial_x u+F_{\prec j}
 H\partial_x u_{\sim j})
  \end{multline}
because the other terms vanish by support considerations.
 The
contribution to the RHS of~\eqref{eq:wj} of the first of these terms
 (which is nothing but $w_j^+ H\partial_x u_{\prec j}$) is  estimated by conormal
product laws~\eqref{eq:BH3} ($T_{X^{-1,\frac 1 2,1}}: X^{\frac 1
 4,\frac 1 2,1}\mapsto X^{\frac 1 4, -\frac 1 2 , 1} $), noticing that
 the Hilbert transform is harmless with conormal spaces (which depend
 only on the size of the Fourier transform).

Terms which have high frequencies of $F$ will be ``cubic'': heuristically,
$$
\Delta_k (e^{-\frac i 2 \int^x u})= F_k\approx 2^{-k}
\partial_x(F_k)\approx 2^{-k} \Delta_k( u F),
$$
hence we have
$$
\|F_k\|_{L^4_x(L^\infty_t)}\leq 2^{-k} \|u\|_{L^4_x(L^\infty_t)}.
$$
Pick the third term in \eqref{gna}, call it $G_3$, we use
$$
\|\partial_x u_k\|_{L^\infty_x L^2_t}\leq 2^{\frac k 4}
\eta_k \|u\|_{\BL {\frac 3 4} \infty 1 2},
$$
with $\sum_j \eta_j \lesssim 1$, $(\eta_j)_j$ a generic sequence which may change from line
to line, and
$$
\|u^+_j\|_{L^4_x L^\infty_t}\leq \|u\|_{L^4_x(L^\infty_t)},
$$
to get
$$
\|G_3\|_{L^2_{t,x}}\lesssim 2^{-\frac 3 4 j} \eta_j
\|u\|^2_{L^4_x(L^\infty_t)}  \|u\|_{\BL {\frac 3 4} \infty 1 2},
$$
which means the sum over $j$ is in $X^{\frac 3 4,0,\infty}\hookrightarrow X^{\frac 1
  4,-\frac 1 2,1}$ which is ok. This ``trick'' of estimating $F_k$ for
$k\geq j$ by $2^{-k} \partial_x F_k$ yielding another factor $u$ will
be used several times in the sequel and such terms will be referred as
``cubic terms''.

The next term we dispose of is $G_2=u_j^+ [S_{j-1},F_{\prec j}]
H\partial_x u_{\prec j}$: in a similar way,
\begin{eqnarray*}
  \|G_2\|_{L^2_{t,x}} \lesssim & \|u_j^+\|_{L^\infty_t(L^2_x)} 2^{-j}
  \|u\|_{L^4_t(L^\infty_x)} \| H\partial_x u_{\prec
  j}\|_{L^4_t(L^\infty_x)}\\
 \lesssim & \eta_j \|u\|_{\LB \infty {\frac 1 4} 2 1}
  \|u\|_{X^{0,\frac 1 2 ,1}} \|u\|_{\LB 4 {\frac 1 4}  \infty 1}
  2^{-\frac j 2},
\end{eqnarray*}
which again means the sum over $j$ to be in  $X^{\frac 1 2,0,\infty}\hookrightarrow X^{\frac 1
  4,-\frac 1 2,1}$.

The third term from \eqref{gna} splits into two terms: the first one, which
has $F_{\sim j}$, can be treated exactly like the high frequencies
interactions, to end up in $X^{\frac 3 4,0,\infty}$. The very last one
can be rewritten as
$$
 u_j^+ S_{j-1} (F_{\prec j}
 H\partial_x u_{\sim j})=u_j^+ F_{\prec j} S_{j-1}(H\partial_x u_{\sim
 j})+u_j^+ [S_{j-1},F_{\prec j}] H\partial_x u_{\sim j},
$$
which are two terms identical to the very first one and second in
\eqref{gna}, up to the replacement of $u_{\prec j}$ by respectively $S_{j-1}
u_{\sim j}$ and $u_{\sim j}$, which is harmless. To recap,
$$
\|\sum_j P_1\|_{X^{\frac 1 4,-\frac 1 2,1}}\lesssim \|w\|_{X^{\frac 1
    4,\frac 1 2,1}}\|u\|_{X^{0,\frac 1 2,1}}+\|u\|_{X^{0,\frac 1
    2,1}}\|u\|_{Y^\frac 1 4}^2+\|u\|^3_{Y^\frac 1 4}.
$$
\item Term $P_2=(S_{j-1}( u F)-S_{j-1} u S_{j-1} F) \partial_x
 u^+_j$: we have
  \begin{multline}
\label{gnabis}
     P_2 = \partial_x u_j^+
 [S_{j-1},F_{\prec j}] u_{\prec j} + \partial_x u_j^+ S_{j-1} ( \sum
 _{j\lesssim k} F_{\sim k} u_{\sim k})\\
 + \partial_x u_j^+ S_{j-1}(F_{\sim j} S_{j-1} u+F_{\prec j} u_{\sim j})
  \end{multline}
which is essentially the same term as before (summing a
$2^{-k}$ will kill the shifted derivative), and therefore
can be estimated in the same way. The only term which requires a
slightly different treatment is the very last one, namely
$$
P_{24}=\partial_x u_j^+ S_{j-1}(F_{\prec j} u_{\sim j}).
$$
Rewrite again
\begin{eqnarray*}
  P_{24} = & \partial_x u_j^+ [S_{j-1},F_{\prec j}] u_{\sim j}+\partial_x
u^+_j F_{\prec j} S_{j-1} u_{\sim j}\\
 = & \partial_x (w^+_j S_{j-1} u_{\sim j})-w_j^+\partial_x S_{j-1}
 u_{\sim j} - u^+_j S_{j-1} u_{\sim j} \partial_x F_{\prec j},
\end{eqnarray*}
and now the first two terms are ok by conormal product rules (in a
sense, they are both high-high frequencies interactions which are
already present in the $f_j^+$ term) and the last one is again cubic.
\item We now proceed with the $ F_{\prec j} f_j^+$ term:
 the very first term $F_{\prec j} f_{j,1}$ sums up to $T_{\partial_x u} w$. By using the conormal
 estimates~\eqref{eq:BH3}, we can estimate the para-product in
 $X^{\frac 1 4 , -\frac 1 2, 1}$. The second one is (where we only
 retain the diagonal term for notational convenience)
$$
P_3=F_{\prec j} \partial_x (\Delta_j \sum_{j\lesssim k} u^+_k u_k).
$$ Again,
  we would like to have $w_k$ rather than $u_k$: recall
$$
u^+_k=F^{-1} w^+_k-F^{-1} \sum_{k\lesssim l} F_l u^+_k.
$$
Obviously, the simplest case is when we have $u^+ u^-$, for the $F$
factors cancel (because $\bar F= F^{-1}$ as an imaginary exponential) and we obtain
$$
P_{31}=F_{\prec j} \partial_x (\Delta_j \sum_{j\lesssim k} w^+_k
w^-_k)
$$
as the main term which, using~\eqref{eq:BH1} and
Proposition~\ref{prop3.4}, leads to
$$
\|\sum_j P_{31}\|_{X^{\frac 1 4, - \frac  1 2 ,1}} \lesssim 
C (1 + \|u\|_{X^{0, \frac 1 2 , 1}} ) \|w \|^2_{X^{\frac 1 4, \frac 1
    2, 1}}.
$$ 
The other term is
$$
P_{32}=F_{\prec j} \partial_x \left(\Delta_j \sum_{j\lesssim k} w^+_k (
\sum_{l>k} F_l u^-_k)\right),
$$
which is obviously cubic and can be dealt with as before, to end up in  
$$
 X^{\frac 3 4,0,\infty}\hookrightarrow X^{\frac 1 4,-\frac 1 2,1}.
$$
Let us now study the case when we have $u^+ u^+$. By support
considerations, the only term appearing is 
$$
   F_{\prec j} \partial_x \Delta^+_j( u_{\sim j}^+ u_{\sim
 j}^+)=  \partial_x \left( F_{\prec
 j} \Delta^+_j
(u_{\sim j}^+ u_{\sim j}^+)\right)-(\partial_x F_{\prec
 j}) \Delta^+_j
(u_{\sim j}^+ u_{\sim j}^+)
$$
the second term is cubic and can be estimated as
before, we only have to estimate
\begin{eqnarray*}
  \partial_x \left( F_{\prec
 j} \Delta^+_j
(u_{\sim j}^+ u_{\sim j}^+)\right)= & \partial_x \Delta^+_j
(F_{\prec j} u_{\sim j}^+ u_{\sim j}^+)+\partial_x \left( [F_{\prec
 j}, \Delta^+_j]
(u_{\sim j}^+ u_{\sim j}^+)\right),\\
 = &  \partial_x \Delta^+_j
(w_{\sim j}^+ u_{\sim j}^+)+\partial_x \left(\Delta_j^+(F_{\sim j}
 u_{\sim j}^+ u_{\sim j}^+)+ [F_{\prec
 j}, \Delta^+_j]
(u_{\sim j}^+ u_{\sim j}^+)\right),
\end{eqnarray*}
for which the second and third terms are, once again, cubic.

Finally, using product law ~\eqref{eq:BH1},
 the remaining term is estimated in $X^{\frac 1 4 , -\frac 1 2, 1}$.
\item Let us finally deal with the commutator which appear
  in $f_j$:
$$
C_j=F_{\prec j} [\Delta_j^+,u_{\prec j}]\partial_x u^+_{\sim
  j}.
$$
We commute $F_{\prec j}$ to obtain $F_{\prec j} u^+_{\sim j}$ which we
know is $w_{\sim j}^++ F_{\sim j} u_{\sim j}^+$. The contribution of
$F_j u_j^+$ will then be cubic, and the additional commutators, namely $[F_{\prec j},\Delta_j]$ or
$[F_{\prec j},\partial_x]$ all gain regularity and will yield cubic terms. Thus, we
are finally left with estimating (at worse !)
$$
\tilde C_j= [\Delta_j^+,u_{\prec j}]  
\partial_x w_{\sim j}^+,
$$
which we intend to deal with in conormal spaces. Let us postpone the
issue and turn to the detail of the cubic terms:
\begin{eqnarray*}
  C_j = & F_{\prec j} \left ( \Delta_j^+(u_{\prec j} \partial_x
    u_{\sim j}^+)- u_{\prec j} \partial_x u^+_j\right)\\
 = &   \Delta_j^+( F_{\prec j} u_{\prec j} \partial_x
    u_{\sim j}^+)- u_{\prec j} F_{\prec j} \partial_x u^+_j+[F_{\prec
    j},\Delta_j^+] (u_{\prec j} \partial_x u^+_{\sim j})\\
 = & \Delta_j^+(u_{\prec j} \partial_x (F_{\prec j} u_{\sim
    j}^+))-u_{\prec j} \partial_x (F_{\prec j} u^+_j)\\
 & {}-\Delta_j^+(u_{\prec j} (\partial_x F_{\prec j}) u_{\sim
    j}^+)+u_{\prec j} (\partial_x F_{\prec j}) u^+_j +[F_{\prec
    j},\Delta_j^+] (u_{\prec j} \partial_x u^+_{\sim j})\\
 = & [\Delta_j^+,u_{\prec j}] \partial_x w^+_{\sim
    j}+\Delta_j^+(u_{\prec j}\partial_x (F_{\sim j} u^+_{\sim j}))\\
 & {}-\Delta_j^+(u_{\prec j} (\partial_x F_{\prec j}) u_{\sim
    j}^+)+u_{\prec j} (\partial_x F_{\prec j}) u^+_j +[F_{\prec
    j},\Delta_j^+] (u_{\prec j} \partial_x u^+_{\sim j})
\end{eqnarray*}
for which the last for terms are cubic. Hence we are left with $\tilde
 C_j$: we will rely on the
next lemma, which tells us we really have an harmless variant of
 $T_{\partial_x u} w$, and we are done.
\end{itemize}
\begin{lemme}
\label{com}
  Let us consider $G$ which is spectrally localized at $|\xi|\leq
  2^{j-1}$ and $F$ which is spectrally localized at $|\xi|\sim
  2^j$. Then one may estimate $[\Delta_j,G] \partial_x F$ in conormal
  spaces as if it were $(\partial_x G) F$.
\end{lemme}
\dem This will follow from a careful rewriting of the commutator:
\begin{align*}
  C = & -\int_y 2^j \phi(2^j(x-y))(G(x)-G(y))\partial_y F(y)\,dy,\\
 = &  -\int_0^1 \int_y 2^j(x-y) \phi(2^j(x-y))G'(y+\theta(x-y))\partial_y F(y)\,dy\,d\theta.
\end{align*}
Set $\psi(z)=z\phi(z)$ and $\psi_j(z)=2^j\psi(2^jz)$, $G'=g$ and $2^-j
\partial_x F=f$, 
$$
C= -\int_0^1 \int_y \psi_j(x-y)g(y+\theta(x-y))f(y)\,dy\,d\theta.
$$
Let us denote by $I_\theta(x)$ the integral over $y$, with a fixed $\theta$. By
Plancherel, and through changes of variable,
\begin{align*}
I_\theta(x)= & \int_{\xi,\eta} e^{ix\xi} \hat \psi_j(\xi) e^{i\eta
  \frac{\theta}{1-\theta} x} \frac{1}{1-\theta}\hat
g(\frac{\eta}{1-\theta}) \hat f(\xi-\eta)\,d\eta \,d\xi,\\
I_\theta(x)= & \int_{\lambda,\mu} e^{ix(\mu+(1-\theta)\lambda+\theta \lambda)} \hat \psi_j(\mu+(1-\theta)\lambda) \hat
g(\lambda) \hat f(\mu)\,d\mu \,d\lambda,\\
I_\theta(x)= & \int_{\xi,\eta} e^{ix\xi} \hat \psi_j(\xi-\eta+(1-\theta)\eta) \hat
g(\eta) \hat f(\xi-\eta)\,d\eta \,d\xi,\\
 = &  \mathcal{F}^{-1}_\xi \left(\int_\eta  \psi_j(\xi-\theta\eta) \hat
g(\eta) \hat f(\xi-\eta)\,d\eta\right).
  \end{align*}
Hence, $I_\theta(x)$ is the inverse Fourier transform of a restricted
convolution between $\hat g$ and $\hat f$: but all conormal spaces
estimates are proven using Plancherel and cutting the Fourier space
into carefully chosen blocks: here we only get part of them, as the
presence of $\hat \psi_j(\xi-\theta \eta)$ reduces the number of
situations where the convolution is non zero. Therefore, we can
estimate $I_\theta(x)$ as if it were $g f\sim G' F$, independently of $\theta$, which is the
desired result.\cqfd

Collecting all terms, we obtain that our source term (after summing
over $j$) is controlled in $X^{\frac 1 4,-\frac 1 2,1}$ by
$$
 \|w\|_{X^{\frac 1
    4,\frac 1 2,1}}\|u\|_{X^{0,\frac 1 2,1}}+\|u\|_{X^{0,\frac 1
    2,1}}\|u\|_{Y^\frac 1 4}^2+\|u\|^3_{Y^\frac 1
 4}+(1+\|u\|_{X^0,\frac 1 2,1})\|w\|^2_{X^{\frac 1 4,\frac 1
 2,1}}+\|u\|^2_{Y^\frac 1 4}\|w\|_{Y^\frac 1 4}.
$$
Using \eqref{eq:equivY}, and inverting the linear operator, we obtain
our a priori estimate.\cqfd

This estimate, when combined with Proposition \ref{prop3.5}, yields an
a priori bound (for small data) on the norm in $X^{\frac 1 4, -\frac 1 2,1}$ of the  r.h.s. of~\eqref{eq:wj}
by 
$$ \|u\|_{X^{0,\frac 1 2,1}}\|w\|_{X^{\frac 1 4,\frac 1 2,1}}
+(1+\|u\|_{X^{0,\frac 1 2,1}})\|w\|^2_{X^{\frac 1 4,\frac 1 2,1}}
+\|w\|^3_{X^{\frac 1 4,\frac 1 2,1}} 
$$
We want local Bourgain spaces instead of global ones (i.e. we want
norms in $X^{s, \frac 1 2, 1}_{T=1}$). For this we have to take
sequences $u_n$ and $w_n$ with supports in $[-2T,2T]$ and equal to $u$
and $w$ respectively minimizing~\eqref{eq.defX} and define $\widetilde
{w}_j^+$ to be the solution of the non linear Schr\"odinger equation
with initial data $w_j^+\mid_{t=0}$ and with a r.h.s obtained by
substituting in the r.h.s of~\eqref{eq:wj} every occurrence of $u$  by
$u_{n}$ and every occurrence of $w$  by $w_n$. In fact, in the analysis
above, there were also parts of $u$ and $w$ for which we used the
classical norms $Y^{\frac 1 4}$. In that case we keep $u$ and $w$ in
the r.h.s. This means that we make the substitution only on the dyadic
parts for which we used $X^{s,\frac 1 2,1}$ norms. We now remark that
$\widetilde{w}_j^+= w_j^+ $ for $|t| \leq T$. Passing to the limit $n
\rightarrow + \infty$ and using lemma~\ref{lem.substit} gives
~\eqref{eq:close}.

Remark that several different $X^{s,b,q}$ norms could have been used
for the same function $u$ (for which the minimizing sequences can differ). This does not matter, as long as we make the substitution
with the sequence corresponding to the norm which is used. This will
be used in the uniqueness Theorem.

 To prove the existence part of our Theorem for $s>\frac 1 4$, we only have to set up a bootstrap argument. Since we have
fixed $T=1$, we can not use any bootstrap on time, but rather will use
again the scale non invariance of the $L^2$ and $H^{1/4}$ spaces: we
fix $u_0\in B^{\frac 1 4 ,1} _2$ and consider $u_\lambda =
\lambda^{\frac 1 2} u(\lambda^2t, \lambda x)$. Then if $u$ is smooth
the norm of $u_\lambda$ and $w_\lambda$ in the spaces above  tend to
$0$ as $\lambda$ tends to $0$, which allow to apply the usual
bootstrap argument.  Existence is then obtained through a limiting
procedure from smooth solutions (notice that passing to the limit in
the equation is trivial, given our a priori control). Finally, continuity of the flow map is a simple
consequence of the classical Bona-Smith argument: for example, one can
implement it exactly as in \cite{KT1} and we therefore skip it.

It remains to prove Theorem \ref{t2}, which is nothing but persistence
of regularity. This requires to carefully check that all nonlinear estimates can
be rewritten with one factor in $(X_{T=1}^{s,b,2}+X_{T=1}^{s,\frac 1 2,1})\cap Y^s$, with $s>\frac 1
4$ and $b>\frac 1 2$. This is certainly obvious on all cubic terms,
and follows from the product rules in the Appendix for quadratic
terms. We leave the details to the reader.

\section{Uniqueness}
\label{sec.unique}
We now prove Theorem~\ref{th.3}.
Suppose we have two solutions $u$ and $v$ to \eqref{eq:bo}, such that
\begin{equation}
  u,v \in Y^{\frac 1 4}\cap X_{T=1}^{0,\frac 1 2,1}.
\end{equation}
Remark that in order to prove Theorem~\ref{th.3}, we can assume that
one of the solutions (say, $u$) is the one we have just constructed (and
consequently $T_{e ^{\mp \frac i 2 \int^x u}} u^{\pm}= w^\pm$ enjoys better
estimates, namely $w\in X_{T=1}^{\frac 1 4,\frac 1 2,1}$).
Later, we will use both informations: the $\frac 1 4$ regularity and the weaker conormal space. Before proceeding, we shall recall that $v$ (and $u$) satisfy
\begin{gather}\label{eq.apriori1}
\|v\|_{X^{0, \frac 1 2, 1}} \leq C\\
\|\Delta_j (v)\|_{L^\infty_x; L^2_t} \leq 2^{-\frac 3 4 j} c_j, \qquad c_j \in l^1\label{eq.apriori2}\\
\Delta_j(v)\|_{L^4_x; L^\infty_t}\leq c_j \in l^1\label{eq.apriori3}
\end{gather}
Furthermore, since $u$ is the solution we just constructed, its renormalized version 
\begin{equation}
 w_j^\pm= S_{j-1} ( e^{\mp \frac i 2 \int^x u}) \Delta_j^\pm (u).
\end{equation}
satisfy the additional estimate
\begin{equation}
\label{eq.apriori4} \|w \|_{X^{\frac 1 4, \frac 1 2, 1}} \leq C\end{equation}let us define
$\delta=u-v$. Then
\begin{equation}
  \label{eq:diffuni2}
  \partial_t \delta+H\partial^2_x \delta+\frac 1 2 \partial_x ((u+v)\delta)=0.
\end{equation}
One would like to obtain an a priori estimate on $\delta$ which would
imply uniqueness, or, even better, Lipschitz dependence in a possibly
weaker norm. However, one cannot process directly with this equation,
again for the same reasons that required a renormalization procedure:
indeed, a typical troublesome term is $2u \partial_x\delta$, or more precisely the
paraproduct $T_u \partial_x \delta$.
\begin{rem}
One may hope to get away with the problem by using a weaker norm,
  namely a norm with negative spatial regularity: a good candidate
  appears to be $X_{T=1}^{-\frac 1 2,\frac 1 2,1}$.
  One can check that we would later need a $T_{X_{T=1}^{0,\frac 1 2}}
  X_{T=1}^{-\frac 3 2,\frac 1 2}\ra X_{T=1}^{-\frac 1 2,-\frac 1 2}$ estimate which
  unfortunately fails. In fact, one particular
  term in the conormal decomposition ends up in $X_{T=1}^{-1,0}$ and no better.
\end{rem}
 We deal with the problem by
another renormalization. For existence, we renormalized the low
frequencies. Here, we lost the symmetry in the nonlinear term, and we
would like to take advantage of the additional properties of $u$, the
``good'' solution. Rewrite the equation, using $v=u-\delta$,
\begin{equation}
  \label{eq:diffuni3}
  \partial_t \delta+H\partial^2_x \delta+ (u-\delta)\partial_x
  \delta+\delta \partial_x u=0.
\end{equation}
which suggests a renormalization using $u-\delta=v$ as the exponential
factor. However, one would like to leave the derivative acting on
products of high frequencies, hence we rewrite once more (setting
$V=u+v$), with paraproduct notations,
\begin{equation}
  \label{eq:diffuni3bis}
  \partial_t \delta+H\partial^2_x \delta+ T_{\partial_x v}
  \delta+T_{v} \partial_x \delta+\partial_x\left(T_\delta
   u+\frac 1 2 R(V,\delta)\right)=0,
\end{equation}
Localizing in frequencies in~\eqref{eq:diffuni3}, we get
\begin{multline}
  \label{eq:diffuni3ter}
  \partial_t \Delta_j^+ (\delta)-i\partial^2_x \Delta_j^+ (\delta)
  +S_{j-1} v\partial_x \Delta_j^+ (\delta)\\
\begin{aligned}
=&-\Delta_j^+\partial_x(\delta_{\prec
  j} u^+_{\sim j}) - \Delta_j^+((\partial_x v_{\prec j})  \delta_{\sim
  j}^+) -\frac 1 2 \Delta_j^+
  \partial_x(\sum_{j\lesssim k} V_{\sim k} \delta_{\sim k})- [\Delta_j^+,v_{\prec j}] \partial_x \delta_j^+ \\
 = & f_{j,1}^++f_{j,2}^++f_{j,3}^++f_{j,4}^+.
\end{aligned}
\end{multline}
Then define $\omega^+$ on $\xi>0$ by $\omega^+_0=P^+ (S_0(\exp{-\frac i 2 \int^x v})S_0
\delta)$ and
$$
\omega^+=\omega^+_0 + T_{\exp{-\frac i 2 \int^x v}} \delta^+, \text{ i.e. }
\omega_j^+ = S_{j-1} (e^{-\frac i 2 \int^x v}) \Delta_j^+\delta \text{
  and }
\omega^+ = \omega^+_0+\sum_{j\geq 0} \omega^+_j.
$$
One may then define  $\omega$ by symmetry with $\omega^-=\bar \omega^+(-\xi)$ on the $\xi<0$
part, so that the low frequencies of $\delta$ and $\omega$ are the
same and $\omega$ is real valued.

From the renormalization estimates, we know that $\delta\in
X_{T=1}^{0,\frac 1 2,1}\hookrightarrow  X_{T=1}^{-\frac 1 4,\frac 1 2,1}$, hence $\omega\in X_{T=1}^{-\frac 1
  2,\frac 1 2,1}$ by using the gauge estimate in
Proposition~\ref{prop3.4}. We will estimate $\omega$ in $X_{T=1}^{-\frac
  1 2, \frac 1 2, 1}$ which, using Proposition~\ref{prop3.5}, smallness
of the data and Proposition~\ref{prop3.4} again, yields estimate~\eqref{eq:lip}.

The equation on $\omega_j$ is 
\begin{multline}
\label{eq:wjbis}
i\partial_t \omega_j^+ +\partial^2_x \omega_j^+=  S_{j-1}(e^{-\frac i 2 \int^x v}) (f_{j,1}^++f_{j,2}^++f_{j,3}^++f_{j,4}^+)\\ +\delta_j^+
 S_{j-1}((v^2+H \partial_x v) e^{-\frac i 2 \int^x v})
+\Bigl ( S_{j-1}(- \frac i 2  v e^{- \frac i 2 \int^x
 v})-iS_{j-1} (v)S_{j-1}(e^{-\frac i 2 \int^x v})\Bigr ) \partial_x \delta_j^+.\\
= P_1+ P_2+ P_3+ P_4+P_5+P_6
\end{multline}
This is essentially the same algebraic calculation as for the
existence part. Note that the $v^2$ gives a cubic
contribution, hence it
will be easier to be dealt with ($\delta$ and $\omega$ can
be estimated at the same regularity level, without conormal spaces). The last term and the
$\partial_x v$ term are, up to commutations, $\omega_j\partial_x
v_{\prec j}$. As before, the whole set of commutation terms are somewhat ``better'' in that they require no conormal spaces to deal with
them; however, we are at a lower regularity level and lost symmetry,
which lead to additional difficulties.

From our set of product estimates (see~\eqref{eq:BH3}), we know that all terms $T_{\partial_x
  v} \omega$ will be ok, meaning
$$
\|T_{\partial_x  v} \omega\|_{X_{T=1}^{-\frac 1 2,-\frac 1 2,1}} \lesssim
\|v\|_{X_{T=1}^{0,\frac 1 2,1}} \|\omega\|_{X_{T=1}^{-\frac 1 2,\frac 1 2,1}}.
$$
In this section, for conciseness, we shall denote by $F^\pm=e^{\pm \frac i 2 \int^x v}$.
We now proceed with estimating in $X_{T=1}^{-\frac 1 2, - \frac 1 2, 1}$ all right-handside terms in \ref{eq:wjbis}.
\begin{itemize}\item Term $P_1=F_{\prec j} \partial_x(\delta_{\prec j} u^+_j)$. As in the existence section, we shall in a first step work in global Bourgain spaces as if $u$ and $v$ were also in the global spaces.

We need to deal with the relationship between $\delta_j$ and $\omega_j$. 
$$
\omega^+_j=F^+_{\prec j} \delta^+_j = F^+\delta^+_j+\sum_{k>j} \Delta_k F^+ \delta^+_j,
$$
\begin{equation}\label{eq.inv}
\delta^+_{j'} = F^- \omega^+_{j'}-\sum_{k>j'} F^- \Delta_k F^+ \delta^+_{j'},
\end{equation}
where we used $F_+^{-1}= F_-$, and summing over $j'\leq j-1$,
$$
\delta^+_{\prec j} = F^- \omega^+_{\prec j}-\sum_{k>j',j'<j} F^- \Delta_k
F^+ \delta^+_{j'}.
$$
Now, the contribution of $\delta_{\prec j}^+$ to $G$, $G^+$, is 
\begin{equation}\begin{aligned}
  P_1^+ = & \partial_x(F^+_{\prec j}\delta_{\prec j} u^+_j)+\delta_{\prec j}
u^+_j S_{j-1}(-\frac i 2 v F^+)\\
= & \partial_x(F_+\delta^+_{\prec j} u^+_j)-\partial_x
\left((F^+-F^+_{\prec j})\delta^+_{\prec j} u^+_j \right)+\delta^+_{\prec j}
u^+_jS_{j-1}(-\frac i 2 v F^+)\\
= & \partial_x(\omega^+_{\prec j} u^+_j)- \partial_x\sum_{k>j', j>j'}
F^+_{k} \delta_{j'}^+ u_j^+ -\partial_x \left((F-F_{\prec
    j})\delta^+_{\prec j}u^+_j\right)+\delta_{\prec j} u^+_jS_{j-1}(-\frac i 2 v F^+)\\
 = & P_{1,1}^++P_{1,2}^++P_{1,3}^++P_{1,4}^+.
\end{aligned}
\end{equation}
(remark that the $F_-$ factor in~\eqref{eq.inv} cancels with the $F_+$ factor in front of $\delta_{\prec j} u^+_j$).
The first term is $\partial_x (T _{\omega^+} u^+)$ which is estimated
according to~\eqref{eq:BH3},
$$
\|\partial_x(T_{\omega^+} u^+)\|_{X^{-\frac 1 2, - \frac 1 2,
    1}}\lesssim \|u\|_{X^{0,\frac 1 2,1}} \|\omega\|_{X^{-\frac 1
    2,\frac 1 2,1}}.
$$  
\begin{rem}
  This term is the one which forces us to go down to $-\frac 1 2$ in
  regularity, apparently wasting $\frac 1 4$.
\end{rem}
Let us show how to estimate the remainder. The very last term
($P_{1,4}^+$) is cubic without derivatives and can be estimated
easily: 
$$
\|\delta_{\prec j} u^+_jS_{j-1}(i v F)\|_{L^2_{t,x}}\lesssim 2^{-\frac
    3 4 j} \|\delta_{\prec j}\|_{L^4_x(L^\infty_t)} 2^{+\frac 3 4 j}
  \|u^+_j\|_{L^\infty_x(L^2_t)} \|v\|_{L^4_x(L^\infty_t)}.
$$ 
But 
$$2^{-\frac
    3 4 j} \|\delta_{\prec j}\|_{L^4_x(L^\infty_t)} \leq 2^{-\frac
    3 4 j} \|\omega_{\prec j}\|_{L^4_x(L^\infty_t)}\leq C \|\omega\|_{X^{-\frac 1 2, \frac 1 2,1}}
$$
and the other terms are controlled by our a priori assumptions~\eqref{eq.apriori2} and~\eqref{eq.apriori3}. Consequently the sum over $j$ will be in $X^{0,0,\infty}\hookrightarrow
X^{-\frac 1 2,-\frac 1 2,1}$.

 The last but one term, $P_{1,3}^+$, is again ``cubic'' because
we can derive $F$ and kill
the $\partial_x$ with it, getting a $v$ instead, hence the same
estimate. The worst one is $P_{1,2}^+$. We have (throwing away $1/4$
regularity on
the $F_k$ factor because of the $j'$ sum)
$$
2^{\frac 3 4 j} u_j \in l^1_j (L^\infty_x L^2_t),\,\,
2^{-\frac 3 4 j'}\delta_{j'}\in L^4_x(L^\infty_t),\,\,
2^{\frac 3 4 k} \Delta_k F\in L^4_x(L^\infty_t),
$$
where we used the usual trick $\Delta_k F=\partial^{-1} \Delta_k
\partial_x F$ on the last term. Collecting everything yields
$$
 2^{-\frac j 4} P_{1,2}^+ \in l^\infty_j (L^2_t L^2_x),
$$
which gives that $P_{1,2}^+\in X^{-\frac 1 4,0,\infty}
\hookrightarrow X^{-\frac 1 2,-\frac 1 2,1}$.

We turn to the $\delta^-_{\prec j}$ contribution  to $P_1$ which when we substitute
doesn't kill the $F$ factor in $F_{\prec j} \partial_x (
\delta_{\prec j}^- u^+_j)$. Recall that $u$ is the solution we
constructed: its renormalized
version $w^+ = T_{e^{ - \frac i 2 \int^x u}} u^+$ belongs to $X^{\frac 1
  4,\frac 1 2,1}$. Furthermore, 
$$
 F^+= e^{ - \frac i 2 \int^x v}= e^{ \frac i 2 \int^x  \delta }e^{ - \frac i 2 \int^x u},
$$
 Thus, if one sees the $e^{- \frac i 2 \int^x  u}$ acting on $u$, this term
becomes after substitution and up to additional ``cubic'' terms as above (to
split the low frequencies of the exponential in the product of low
 frequencies of exponentials)
$$
\partial_x(\tilde F_{\prec j} \omega^-_{\prec j} w^+_j) \text { with }
\tilde F=  e^{\frac i 2 \int^x \delta}.
$$
 We can now apply exactly the same strategy as above to estimate
 $\partial_x( \omega^-_{\prec j} w^+_j)$ and end up in
$ X^{-\frac 1 4, - \frac 1 2,1}$ (again, notice the gain of a quarter of derivatives from substituting
 $w$ to $u$, together with control of $w$ in $X^{\frac 1 4, \frac 1 2, 1}$). Using
 Proposition~\ref{prop3.4} (which loses the quarter of derivatives we
 just gained) we  estimate 
$$
\|T_{e^{\frac i 2 \int^x \delta}}
 \left(\partial_x( T_{\omega^-} w^+)\right)\|_{X^{-\frac 1 2,
     -\frac 1 2, 1}} \lesssim C(\delta) \|\omega\|_{X^{-\frac 1 2,\frac 1
     2,1}}\|w\|_{X^{\frac 1 4,\frac 1 2,1}}.
$$
All other terms we discarded are cubic again and easily disposed of.
\item Let us study $ P_2$. This
  term
  is nothing (up to more cubic terms) but $T_{\partial_x v} \omega^+$ which is estimated
  using~\eqref{eq:BH3} in $X^{-\frac 1 2, -\frac 1 2,1}$.
\item Let us study the contribution of the third term
$P_3=F^+_{\prec j} \partial_x (\Delta_j^+ \sum_{j\lesssim k} V_k \delta_k)$. Again,
  we would like to have $\omega_k$ rather than $\delta_k$. We have 
\begin{equation}
P_3= \partial_x \Delta_j \sum_{j\lesssim k} (F^+_{\prec k} V_k \delta_k
+ (F^+_{\prec j} - F^+_{\prec k}) V_k \delta _k)  + [F^+_{j},
\partial_x \Delta_j] \sum_{j\lesssim k} V_k \delta_k= P_{3,1}+P_{3,2}+P_{3,3} .
\end{equation}
\begin{itemize}
\item term {$P_{3,1}$}: due to support conditions, in the product $V_k
  \delta_k$, the interactions $-,-$ cancel. As a consequence (since
  $V_k= 2u_k- \delta_k$), it is enough to estimate (distributing the
  $F^+_{\prec k}$ factor to $\delta_k^+$)
$$
  \partial_x \Delta_j \sum_{j\lesssim k} (F_{+,\prec k} u_k \delta_k^+=
\partial_x \Delta_jR( u, \omega^+)
$$
as well as (distributing the $F^+_{\prec k}$ factor to $u_k^+$)
$$
  \partial_x \Delta_j \sum_{j\lesssim k} (F_{+,\prec k} u^+_k \delta_k^-=
  \partial_x \Delta_jR( w, \delta).
$$
and
$$
  \partial_x \Delta_j \sum_{j\lesssim k} (F_{+,\prec k} \delta_k^+ \delta_k^-=
  \partial_x \Delta_jR(  \delta, \omega^+).
$$

But since $u$ is the solution we constructed in the previous section, it is bounded in $X^{0,\frac 1 2, 1}$ and its renormalized version, $w$ is bounded in $X^{\frac 1 4, \frac 1 2, 1}$. On the other hand, by assumption, $\omega$ is bounded in $X^{-\frac 1 2, \frac 1 2, 1}$ and consequently, according to Proposition~\ref{prop3.4}, $\delta$ is bounded in $X^{-\frac 3 4, \frac 1 2, 1}$. Finally, $\delta$ is a priori bounded in $X^{0, \frac 1 2, 2}$ and $\omega$ in $X^{-\frac 1 2,\frac 1 2,1}$ by assumption. Now the estimate on $K_1$ follow from product laws, namely ~\eqref{eq:BH1bis}.
\item term $P_{3,2}$: this will be a variation on the cubic term, as
$$ P_{3,2}= \partial_x \Delta_j \sum_{j\leq l \leq k} F^+_{l} V_k \delta _k.
$$
Recall that, using~\eqref{eq.smooth} and Proposition~\ref{prop.4}, 
$$
\|\delta_k\|_{L^\infty_x; L^2_t}\leq C\|\omega_k\|_{L^\infty_x; L^2_t}c_k \| \omega \|_{X^{-\frac 1 2, \frac 1 2, 1}}, \qquad (c_k)_k \in l^1
$$
 and, according to our {\em a priori} assumptions on $v$ and $u$~\eqref{eq.apriori3}, 
$$
\|V_k\|_{L^4_x(L^2_t)} \leq c_k \in l^1
$$

 Hence, 
$$
\|\delta_k V_k\|_{L^4_x(L^2_t)}\leq \eta_k \in l^1_k.
$$
This together with
$$
F_l \approx 2^{-l} ( v F)_l,\Rightarrow \|F_l\|_{L^4_x(L^\infty_t)}\sim 2^{-l} \|v\|_{L^4_x( L^\infty_t)}\sim 2^{-l}
$$
finally we get an estimate in
$$
 X^{0,0,\infty}\hookrightarrow X^{-\frac 1 2,-\frac 1 2,1}.
$$
\item term $P_{3,3}= [F^+_j, \partial_x \Delta_j] \sum_{j\lesssim k} V_k \delta_k$
The following Lemma (which appeared in a slightly different form in~\cite{BuPl04} shows that we can estimate this term as if we had $\partial_x F_j^+$ in place of the commutator (and consequently as we estimated $P_{3,2}$)
\begin{lemme}
\label{lem.commutation}
  Let $g(x,t)$ be such that $\|\partial_x
  g\|_{L^{p_1}_x(L^{q_\infty}_t)}<+\infty$, and $f(x,t)\in
  L^{p_\infty}_x(L^{q_2}_t)$, with $\frac 1 {p_1}+\frac 1
  {p_\infty}=\frac 1 2$ and $\frac 1 {q_\infty}+\frac 1 {q_2}=\frac 1 2$, then
  $h(x,t)=[\Delta_j,g] f$ is bounded in $L^1_x(L^2_t)$ by
$$ C 2^{-j}\|\partial_x
  g\|_{L^{p_1}_x(L^{q_\infty}_t)}\|f(x,t)\|_{  L^{p_\infty}_x(L^{q_2}_t)}
$$
\end{lemme}
\dem We first take $p_1=2$, $p_\infty=\infty$: set $h(x)=[\Delta_j,g]
f$, recall $\Delta_j$ is a convolution by $2^j\phi(2^j \cdot)$, and
denote $\psi(z)=z|\phi|(z)$:
\begin{equation*}\begin{aligned}
    h(x) & = \int_y 2^j\phi(2^j(x-y))(g(y)-g(x))f(y)dy \\ & = \int_{y
,\theta\in [0,1]} 2^j\phi(2^j(x-y))(x-y)g'(x+\theta(y-x))f(y) d\theta
dy \\ |h(x)| & \leq 2^{-j} \int_{y,\theta\in [0,1]}
2^j\psi(2^j(x-y))|g'(x+\theta(y-x))||f(y)| d\theta dy
\end{aligned}
\end{equation*}
and then take successively time norms and space norms,
\begin{equation*}
\begin{aligned}
 \|h(x,t)\|_{L^2_t} & \leq 2^{-j} \int_{y,\theta\in [0,1]}
   2^j\psi(2^j(x-y))\|g'(x+\theta(y-x,t))\|_{L^{q_\infty}_t}\|f(y,t)\|_{L^{q_2}_t}
   d\theta dy\\  \|h(x)\|^2_{L^2_t,x} & \leq 2^{-j}
   \|f\|_{L^\infty_x(L^{q_2}_t)}\int_{\smash{\theta\in [0,1]}} \|\int_y
   2^j\psi(2^j(x-y))\|g'(x+\theta(y-x))\|_{L^{q_\infty}_t}  dy\|_{L^2_x}
   d\theta\\ & \leq 2^{-j}\int_{\smash{\theta\in [0,1]}} \|\int_z
   \frac {2^j} \theta \psi(\frac{2^j} \theta(x-z))\|g'(z)\|_{L^{q_\infty}_t}  dz\|_{L^2_x}
   d\theta \\
&\leq  2^{-j} \int_{\smash{\theta\in [0,1]}}\|g'(x+\theta
   z)\|_{L^2_x; L^{q_\infty}_t} d\theta
\end{aligned}
\end{equation*}
The case $p_1=\infty$, $p_\infty=2$ is identical, exchanging $f$ and
$g'$ (in fact, this would be the usual commutator estimate !). The
general case then follows by bilinear complex interpolation. \cqfd
\end{itemize}
\item Let us now deal with the commutator which appears
  in $f_j$, $P_4$, namely
$$
P_4=F_{\prec j} [\Delta_j,v_{\prec j}] \partial_x \sum_{k\sim j}
\delta_k^+.
$$
First, we may replace the $F_{\prec j}$ factor by $F$: the difference
$F-F_{\prec j}$ will lead to a cubic term where we do not need to take
advantage of the commutator structure. Then, we commute $F$ with
everything else to obtain $F \delta_k$ which we know is $\omega_k^++$''cubic
terms''. The additional commutators, namely $[F,\Delta_j]$ or
$[F,\partial_x]$ all gain regularity and yield cubic terms. Thus, we
are finally left with
$$
\tilde C_j= [\Delta_j,v_{\prec j}] \partial_x \sum_{k\sim j}
\omega_k,
$$
which can be estimated in conormal spaces, exactly as we did for
existence, through Lemma \ref{com}.
\item Let us study 
  $P_5=\delta_j S_{j-1}(\partial_x v F)$. Up to commutator terms
  which are cubic again, we have
$$
P_5=\delta^+_j F_{\prec j} \partial_x v_{\prec j}+\delta_j^+
S_{j-1}(\sum_{k>j} \partial_x v_k  F_k) + \delta_j^+ S_{j-1} ( F_{+, \prec j} v_j + v_{\prec j} F_{+,j}),
$$
because other terms vanish by support considerations. The
first of these terms is $T_{\partial_x v} \omega$ which is ok in conormal
spaces. The second term is ok because
$$
F_k\approx 2^{-k} \partial_x(F_k)=2^{-k} \Delta_k( v F),
$$
which leads to a ``cubic'' term and the last one is basically the same.
\item Term $P_6=(S_{j-1}(vF)-S_{j-1} v S_{j-1} F) \partial_x
  \delta_j$: again, up to commutators, we cancel the second term to
  get
$$
P_6\approx S_{j-1}(\sum_{k>j} v_k F_k) \partial_x \delta^+_j,
$$
which is essentially the same ``cubic'' term as before, and therefore
ends up in $X^{0,0,\infty}$.

\end{itemize}
Finally, we can collect all terms: recall
$$
\partial_t \omega +H\partial_x^2 \omega= K,
$$
and we just proved
$$
\|K\|_{X^{-\frac 1 2,-\frac 1 2,1}}
\lesssim C(u,v) \|\omega\|_{X^{-\frac 1 2,\frac 1 2,1}},
$$
which immediately yields using the equation
$$
\|\omega\|_{X^{-\frac 1 2,\frac 1 2,1}}\lesssim \|\omega_0\|_{B^{-\frac 1
    2,1}_2},
$$
As in the previous section we can now pass from the global Bourgain space to the correct local Bourgain space procedure and obtain
$$
\|\omega\|_{X^{-\frac 1 2,\frac 1 2,1}_{T=1}}\lesssim \|\omega_0\|_{B^{-\frac 1
    2,1}_2},
$$
which provides uniqueness if $\omega_0=0$ and the desired continuity by
inverting the gauge transform and interpolation.\cqfd

\section{Uniqueness for weak finite energy solutions}
\label{sec:nic} In this section we prove uniqueness in  the natural finite energy space, $L^\infty_t(H^{1/2})$.
\begin{theoreme}\label{th.5}
Assume that $v$ is an $L^\infty_t( H^{1/2})$ solution to the Benjamin Ono equation
$$ \partial_t v + H \partial^2_x v = - \frac{\partial_x} 2 ( v^2).$$
 Then $v\in C^1( \mathbb{R}_t; \mathcal{D}'( \mathbb{R}_x))$ and
 consequently $v\mid_{t=0}$ makes sense in $\mathcal{D}'$. Then we
 have 
$$
\|v_{|t=0}\|_{H^{1/2}}\lesssim \|v \|_{L^\infty( H^{1/2})}
$$
 and if $u$ is another such solution, there exist $C$, $T>0$ depending
 both only of $\|v\|_{L^\infty(H^{1/2})}$, $\|u\|_{L^\infty(H^{1/2})}$
 such that 
\begin{equation}
\label{eq.unicite12}
\|u-v\|_{X^{- \frac 1 2, \frac 1 2, 1}_T} \leq  C \|u_{|t=0} -
v_{|t=0}\|_{B^{- \frac 1 2,1}_2}\leq C\|u_{|t=0} -
v_{|t=0}\|_{H^{1/2}}
\end{equation}
\end{theoreme} 
The first part of the Theorem is clear as 
$$
\partial_t v = H \partial_x^2 v + \frac {\partial_x} 2 ( v^2) \in
L^\infty_t; \mathcal{S}'.
$$
To prove control of $v_{|t=0}$ in $H^{1/2}$, we simply take a
sequence $t_n>0 \rightarrow 0$ such that $\|v(t_n)\|_{H^{1/2} }\leq
\|v\|_{L^\infty; H^{1/2}}$ and use $v(t_n) \rightarrow v
\mid_{t=0}$ in $\mathcal{D}'$. The main point is of course to
prove~\eqref{eq.unicite12}. We will prove this estimate
assuming first that $u$ is the strong solution of~\eqref{eq:bo} we
constructed in Section~\ref{sec.key}. Then applying this estimate to
$v\mid_{t=0}=u_0$, it implies uniqueness and consequently the fact
that $u$ is a strong solution is no longer an assumption.

Applying $e^{t H \partial_x^2}$ to $v$, remark that $v$ is a solution of the Duhamel equation
$$
 v(t, \cdot)= e^{-t H \partial_x^2} v\mid_{t=0} - \int_0^t ( e^{-(t-s)
 H \partial_x^2} (\frac {\partial_x} 2 ( v^2)).
 $$ 
We  proceed to prove a few a priori bounds on $v$ (and $u$).
By rescaling, we consider again the small datum case, with $T\sim 1$.
\begin{lemme}\label{lem.apriori}
  Assume that $v\in L^\infty_T(H^{1/2})$, then for any $\epsilon>0$
  there exist $C>0$ such that for any $j$
$$
\|\Delta_j( v)\|_{L^4_T(L^\infty)}\leq C (1+T)^{1/4} 2^{j(-\frac 1 4
  +\epsilon)}\|v\|^{\frac 3 2}_{L^\infty(-1, T+1;
  H^{1/2})}(1+\|v\|^{\frac 1 2}_{L^\infty(-1,T+1; H^{1/2})}).
$$
As a by product, weak solutions are $L^4_{t,\text{loc}}( L^\infty)$.
\end{lemme}
\dem Consider $\chi\in C^\infty_0( \mathbb{R}_t)$. We have
$$
 \|\chi(2^j t) \Delta_j (v)\|_{L^2_{t,x}}\leq C 2^{-j}
$$
and  $\chi(2^j t) \Delta_j (v)$ satisfies
\begin{equation}
\label{eq.trunc} (\partial_t+ H\partial_x^2) \chi(2^j t)\Delta_j(v)=
2^{ j} \chi'(2^j t) \Delta_j(v) +\frac 1 2\chi(2^j t)\partial_x
\Delta_j \left( \sum_{k\leq j} v_kv_j   + \sum_{k>j} v_k^2 \right)
\end{equation}
and (using Sobolev injection on one $v_k$ factor and the fact that the
integrals in time are taken on intervals of length of size $2^{-j}$)
we estimate the right hand side in $L^2_{t,x}$ by
\begin{equation*}
  \|v\|_{L^\infty_t(H^{1/2})}+ (2^j\sum_{k\leq j} 2^{-j}+ 2^j\sum_{k\geq j} 2^{-k})\|v\|^2_{L^\infty_t( H^s)}
\lesssim
(1 +j) \|v\|_{L^\infty_t(H^{1/2})}  
\end{equation*}
As a consequence we obtain the following bounds
$$
 \|\chi({2^j}t) \Delta_j(v)\|_{X^{0,0,2}}\lesssim
 2^{-j}\|v\|_{L^\infty_t; H^{1/2}}, \qquad \|\chi(2^j t)
 \Delta_j(v)\|_{X^{0,1,2}} \lesssim (1+j)\|v\|_{L^\infty_t(H^{1/2})}.
$$
Interpolating between these two bounds ($j$ is fixed), we obtain
(remark that even though it is not essential, we regain the $l^1$
summability with respect to the $k$ index by using
Lemma~\ref{lem.interpo})
$$
 \|\chi(2^j t) \Delta_j(v)\|_{X^{0,\frac 1 2,1}}\lesssim \sqrt
 j2^{- \frac j 2}\|v\|_{L^\infty_t(H^{1/2})}.
$$
But, according to~\eqref{eq.strich}, we have 
$$
\|\chi(2^j t) \Delta_j(v)\|_{L^4_t( L^\infty)} \lesssim
\|\chi(2^j t) \Delta_j(v)\|_{X^{0,\frac 1 2,1}},
$$
and applying this inequality between times $0,2^{-j}$, $2^{-j},2\times
2^{-j},...,2^j\times 2^{-j}=1$, we obtain
$$
 \|\Delta_j(v)\|_{L^4_t(L^\infty)}\lesssim \sqrt j 2^{- \frac j 4}
\|v\|_{L^\infty_t(H^{1/2)}}(1+ \|v\|^{\frac 1 2}_{L^\infty_t(H^{1/2})}),
$$
which is definitely summable over $j$: in fact, we obtained that $v\in
\LB 4 {\frac 1 4 -\epsilon} \infty 1 \hookrightarrow L^4_t
L^\infty_x$.\cqfd
\begin{cor}\label{cor.1}
For any $0\leq \theta <1$, we have
$$
\|v\|_{X^{\frac 1 2 -\theta, \theta , 2}}\lesssim \|v\|_{L^\infty_t(H^{1/2})}(1+ \|v\|^{\frac 1 2}_{L^\infty_t(H^{1/2})}).
$$
\end{cor}
We consider the equation satisfied by $\chi(t) v$
$$ (\partial_t + H \partial_x^2) \chi \Delta_j(v) =  \chi'( t ) \Delta_j(v) +\frac 1 2\chi(t)\partial_x \Delta_j \left( \sum_{k\leq j} v_kv_j   + \sum_{k>j} v_k^2 \right)  
$$
and we estimate the right hand side in $L^2_{t,x}$ (using the
$L^4(L^\infty)$ estimate we just proved) by
$$ 
2^{-\frac j 2} + 2^{j}\sum_{k\leq j} \sqrt k 2^{k(- \frac 1 4)}
2^{-\frac j 2} + 2^j \sum_{k\geq j} \sqrt k 2^{- \frac {3k} 4}
$$
As a consequence we obtain
$$ \|\chi(t)  \Delta_j(v)\|_{X^{0,0,2}}\leq C 2^{-\frac j 2}C\|v\|_{L^\infty_t(H^{1/2})}, \qquad \|\chi( t) \Delta_j(v)\|_{X^{0,1,2}} \leq C2^{\frac j 2} \|v\|^2_{L^\infty_t(H^{1/2})}  
$$
which reads
\begin{equation}
\begin{aligned} \|\Delta^\pm_{j,k}(\chi(t)v) \|_{L^2_{t,x}}&\leq C
  2^{-\frac j 2}C\|v\|_{L^\infty_t(H^{1/2})} c_k; \qquad
  \|c_k\|_{l^2_k}\leq 1\\
\|\Delta^\pm_{j,k}(\chi(t)v) \|_{L^2_{t,x}}&\leq C 2^{\frac j 2-k}\|v\|^2_{L^\infty_t(H^{1/2})}c_k; \qquad \|c_k\|_{l^2_k}\leq 1
\end{aligned}
\end{equation}
and the use of Lemma~\ref{lem.interpo} (with fixed $k$) to
regain the $l^2$ summability in $j$) gives the Corollary.\cqfd

We now return to the proof of~\eqref{eq.unicite12}. From the lemma, we know that $\delta\in
X_{T=1}^{- \epsilon,\frac 1 2+ \epsilon,2}\hookrightarrow  X_{T=1}^{-\frac 1 4,\frac 1 2,1}$, hence $\omega\in X_{T=1}^{-\frac 1
  2,\frac 1 2,1}$ as before. We will estimate again $\omega$ in $X_{T=1}^{-\frac
  1 2, \frac 1 2, 1}$.

Before proceeding, recall that $v$ satisfy 
\begin{gather} \label{eq.apriori5}
\|v\|_{X^{- \epsilon, \frac 1 2, 1}} \leq C\Rightarrow \|\Delta_j(v)\|_{L^4_t( L^\infty_x)}\leq C 2^{\epsilon j }\\
\|v\|_{L^\infty( H^{1/2})}\leq C\Rightarrow \|\Delta_j(v)\|_{L_t^\infty( L^2)}\leq C 2^{-\frac j 2}\label{eq.apriori8}
\end{gather}
according to Corollary~\ref{cor.1}. Furthermore, since $u$ is the
solution we just constructed, it satisfies~\eqref{eq.apriori1},
\eqref{eq.apriori2} and~\eqref{eq.apriori3} and its renormalized
version satisfy the additional estimate~\eqref{eq.apriori4}
(one could actually prove it to be even better, namely $X^{\frac 1 2, \frac 1 2, 2}$.

Recall all terms $T_{\partial_x
  v} \omega$ are under control according to~\eqref{eq:BH4bis},
$$
\|T_{\partial_x  v} \omega\|_{X_{T=1}^{-\frac 1 2,-\frac 1 2,1}} \lesssim
\|v\|_{X_{T=1}^{0,\frac 1 2,2}} \|\omega\|_{X_{T=1}^{-\frac 1 2,\frac 1 2,2}}.
$$

We now proceed with estimating in $X_{T=1}^{-\frac 1 2, - \frac 1 2, 1}$ all right-handside terms in \eqref{eq:wjbis}. We follow closely the estimates in the previous section and will only point out the differences. As before we forget in the first step the local Bourgain spaces and work with global ones.
\begin{itemize}\item Term $P_1=F_{\prec j} \partial_x(\delta_{\prec j} u^+_j)$.
\begin{itemize}
\item Term $P_{1,1}^+=\partial_x (T _{\omega^+} u^+)$ which is estimated
according to~\eqref{eq:BH3}: the solution $u$ that we constructed in Section~\ref{sec.4} is in $X^{0, \frac 1 2,1}$ (it is actually even better since the initial data is in $H^{1/2}$, it is in $X^{\frac 1 4-0, \frac 1 2,1}$)
\item Term $P_{1,4}^+$. We remark that interpolating between~\eqref{eq.apriori5} and~\eqref{eq.apriori8}, for any $\theta>0$, $v$ is bounded in
$$ L^{\frac 4 {(1- \theta)}}_t( L^{\frac 2 \theta}).$$
We choose $\theta$ arbitrary close to $0$ and write
$$
\|\delta_{\prec j} u^+_jS_{j-1}(i v F)\|_{L^2_{t,x}}\lesssim 2^{-\frac 1 2j} \|\delta_{\prec j}\|_{L^{\frac 4{1- \theta}}_t(L^{ \frac 2 \theta})} 2^{\frac 1 2 j}
  \|u^+_j\|_{L^{\frac 2 \theta}_t(L^{\frac 2 {(1- 2\theta)}})} \|v\|_{L^{\frac 4 {(1- \theta)}}_t(L^{\frac 2 \theta})}.
$$ 
and using~\eqref{eq.apriori2} and~\eqref{eq.apriori3} we can estimate 
$$
 \|u^+_j\|_{L^{\frac 2 \theta}_x(L^{\frac 2 {(1- 2\theta)}}_t)}\leq 2^{(\epsilon-\frac 1 2) j}$$
As a consequence, $P_{1,4}^+$ is in $X^{-\epsilon,0,\infty}\hookrightarrow
X^{-\frac 1 2,-\frac 1 2,1}$.
\item Term $P_{1,3}^+$. It is again ``cubic'' because
we can derive $F$ and kill
the $\partial_x$ with it, getting a $v$ instead, hence the same
estimate.
\item Term $P_{1,2}^+$. We write $F_k \sim 2^{-k} ( vF)_k$
We remark that by support condition, we can restrict the sum to the set $k\leq j$. We have to estimate
$$
 2^ j\sum_{j'\leq k\leq j} 2^{-k} \|\delta_{j'}\|_{L^{\frac 4{1-
 \theta}}_t(L^{ \frac 2 \theta})}\|(vF)_k\|_{L^{\frac 4{1-
 \theta}}_t(L^{ \frac 2 \theta})} \|u_j\|_{L^{\frac 2 \theta}_t(
 L^{\frac 2 {1 - 2 \theta}}_t)}\leq 2^ j\sum_{j'\leq k\leq j}2^{\frac 1
 2 j'+(\epsilon - \frac 1 2) j}\\
\leq 2^{\epsilon j}$$ 
which gives an estimate in $X^{-\epsilon,0,\infty}$ 
\end{itemize}
\item The main terms in the $\delta^-_{\prec j}$ contribution  to $P_1$ are estimated as in the previous section whereas the cubic remainders are estimated as above.
\item Let us study $ P_2$. This
  term
  is nothing (up to more cubic terms) but $T_{\partial_x v} \omega^+$ which is estimated
  using~\eqref{eq:BH4bis} in $X^{-\frac 1 2, -\frac 1 2,1}$.
\item Let us study the contribution of the third term
$P_3=F^+_{\prec j} \partial_x (\Delta_j^+ \sum_{j\lesssim k} V_k \delta_k)$. Again,
  we would like to have $\omega_k$ rather than $\delta_k$. We have 
\begin{equation}
P_3= \partial_x \Delta_j \sum_{j\lesssim k} (F^+_{\prec k} V_k \delta_k
+ (F^+_{\prec j} - F^+_{\prec k}) V_k \delta _k)  + [F^+_{j},
\partial_x \Delta_j] \sum_{j\lesssim k} V_k \delta_k= P_{3,1}+P_{3,2}+P_{3,3} .
\end{equation}
\begin{itemize}
\item Term {$P_{3,1}$} is estimated as in the previous section
\item term $P_{3,2}$: this will again be a variation on the cubic term, as
$$ P_{3,2}= \partial_x \Delta_j \sum_{j\leq l \leq k} F^+_{l} V_k \delta _k.
$$
Recall that, using~\eqref{eq.strich}, 
$$
\|\delta_k\|_{L^4_t(L^\infty_x)}\lesssim \|\omega_k\|_{L^4_t(L^\infty_x)}\leq c_k 2^{\frac k 2} \| \omega \|_{X^{-\frac 1 2, \frac 1 2, 1}}, \qquad (c_k)_k \in l^1.
$$
On the other hand, using~\eqref{eq.apriori8},
$$ \|V_k\|_{L^\infty_t( L^2)} \leq C 2^{-\frac k 2}$$
and according to Lemma~\ref{lem.apriori}
$$ \|v\|_{L^4_t(L^\infty)} \lesssim 1$$
As a consequence
\begin{equation}
\|P_{3,2}\|_{L^2_{t,x}} \lesssim 2^j \sum_{j\leq l\leq k}
2^{-l}c_k\lesssim 1.
\end{equation}
This yields an estimate in
$$
 X^{0,0,\infty}\hookrightarrow X^{-\frac 1 2,-\frac 1 2,1}.
$$
\item Term $P_{3,3}$: using the following Lemma (a simplified version of Lemma~\ref{lem.commutation}), we can estimate this term, exactly as we estimated $P_{3,2}$.
\begin{lemme}
\label{lem.commutationbis}
  Let $g(x)$ be such that $\|\partial_x
  g\|_{L^{\infty}_x}<+\infty$, then we have
$$ \|[\Delta_j,g] f\|_{L^2_x}\leq C 2^{-j}\|\partial_x
  g\|_{L^{\infty}_x}\|f(x,t)\|_{  L^{2}_x}
$$
\end{lemme} 
\end{itemize}
\item Let us now deal with the commutator which appears
  in $f_j$, $P_4$, namely
$$
P_4=F_{\prec j} [\Delta_j,v_{\prec j}] \partial_x \sum_{k\sim j}
\delta_k^+.
$$
First, we may replace the $F_{\prec j}$ factor by $F$: the difference
$F-F_{\prec j}$ will lead to a cubic term where we do not need to take
advantage of the commutator structure. Then, we commute $F$ with
everything else to obtain $F \delta_k$ which we know is $\omega_k^++$''cubic
terms''. The additional commutators, namely $[F,\Delta_j]$ or
$[F,\partial_x]$ all gain regularity and yield cubic terms. Thus, we
are finally left with
$$
\tilde C_j= [\Delta_j,v_{\prec j}] \partial_x \sum_{k\sim j}
\omega_k,
$$
which can be estimated in conormal spaces, exactly as we did for
existence through Lemma~\ref{com}.
\item Let us study 
  $P_5=\delta_j S_{j-1}(\partial_x v F)$. Up to commutator terms
  which are cubic again, we have
$$
P_5=\delta^+_j F_{\prec j} \partial_x v_{\prec j}+\delta_j^+
S_{j-1}(\sum_{k>j} \partial_x v_k  F_k) + \delta_j^+ S_{j-1} ( F_{+, \prec j} v_j + v_{\prec j} F_{+,j}),
$$
because other terms vanish by support considerations. The
first of these terms is $T_{\partial_x v} \omega$ which is ok in conormal
spaces (according to~\eqref{eq:BH4bis}). The second term is ok because
$$
F_k\approx 2^{-k} \partial_x(F_k)=2^{-k} \Delta_k( v F),
$$
which leads to a ``cubic'' term and the last one is basically the same.
\item Term $P_6=(S_{j-1}(vF)-S_{j-1} v S_{j-1} F) \partial_x
  \delta_j$: again, up to commutators, we cancel the second term to
  get
$$
P_6\approx S_{j-1}(\sum_{k>j} v_k F_k) \partial_x \delta^+_j,
$$
which is essentially the same ``cubic'' term as before, and therefore
ends up in $X^{0,0,\infty}$.

\end{itemize}
Finally, we can collect all terms: recall
$$
\partial_t \omega +H\partial_x^2 \omega= K,
$$
and we just proved 
$$
\|K\|_{X^{-\frac 1 2,-\frac 1 2,1}}
\lesssim C(u,v) \|\omega\|_{X^{-\frac 1 2,\frac 1 2,1}},
$$
with a constant $C$ small if $u$ and $v$ are small in $L^\infty; H^{1/2}$ (remark that the norms of the strong solution $u$ involved are controlled by $\|u+0\|_{H^{1/2}}$ and hence by $\|u\|_{L^\infty; H^{1/2}}$). This immediately yields using the equation
$$
\|\omega\|_{X^{-\frac 1 2,\frac 1 2,1}}\lesssim \|\omega_0\|_{B^{-\frac 1
    2,1}_2},
$$
We now use the substitution procedure which gives 
$$
\|\omega\|_{X^{-\frac 1 2,\frac 1 2,1}_{T=1}}\lesssim \|\omega_0\|_{B^{-\frac 1
    2,1}_2},
$$
\section{Existence and uniqueness in the $s<\frac 1 4$ range}
We now provide an outline of the proof of Theorem \ref{t1bis}. A
detailed proof involves a complete rewriting of the previous sections, which, given the result from \cite{IoKe}, is not
worth the effort at this stage where the uniqueness part is not optimal. 

We start with existence in the $0<s<\frac 1 4$ range. We seek the renormalized function $w\in X^{\e,\frac 1 2,1}$. Assume
for now that the mapping from Proposition \ref{prop3.4} holds true, despite the fact that we now expect $u\in X^{-\frac 1 4+\e,\frac 1
  2,1}$. Moreover, we also assume $u\in  X^{-\frac 1 8+\e,\frac 1
  4,1}$, by using  Proposition \ref{prop3.4} when $b=\frac 1 4$ rather
than $\frac 1 2$ together with the a priori knowledge on $w$. Ignore as well the cubic terms in the equation for $w$, we
decompose
$$
w=w_L+w_1+w_2+w_3,
$$
where $w_L$ is the linear part, hence $w_L\in X^{\e,\infty,1}$. The
three other terms come from inverting the linear operator on the
bilinear term (which involves interactions between $u$ and $w$).
\begin{itemize}
\item We
set $w_1$ to be the term associated to source terms $R_1(u,w) \in X^{\e+(\e-\frac 1
  4)+\frac 1 2,-\frac 1 2,1}$, which gains $\frac 1 4$ spatial
regularity. These terms come from all bilinear interactions where, in
the notation of the Appendix, the output is $X^{s+s'+\frac 3 2,-\frac
  1 2,1}$.
\item Next, $w_2$ comes from  (reminder) source terms $R_2(u,w) \in X^{\e+(\e-\frac 1 4),-\frac 1 4,1}\cap X^{\e+(\e-\frac 1
  8),-\frac 1 2,1}\cap X^{\e+\e,-\frac 3 4,1}$; moreover, it is worth
noting that $R_2$ is localized in the $k''<2j''$ region, with notations from the
Appendix (and only involves the $j''\sim j' \sim j<k/2$ zone on the $u,w$ factors), hence
$w_2\in X^{2\e-\frac 1 4,\frac 3 4,1}\cap X^{2\e-\frac 1 8,\frac 1
  2,1}\cap X^{2\e,\frac 1 4,1}$. 
\item Finally,
$w_3$ arises from the ``worst'' same spatial frequency interactions $R_3(u,w)$,
which are (picking $u\in X^{\e-\frac 1 8,\frac 1
  4,1}$) in $X^{2\e-\frac 1 8,-\frac 1 2,1}$ and
localized at $k''\sim 2j''$, hence $w_3\in X^{2\e-\frac 1 8,\frac 1 2,1}\cap X^{\e+(\e-\frac 1
  8)+\frac 1 2,\frac 1 4,1}$, and moreover, its Fourier localization is in the
region $j''\sim 2k''$.
\end{itemize}
Obviously, $w_L\in X^{\e+(\e-\frac 1 4),\frac 3 4,1}$, and we may
relabel $w_2$ to be $w_L+w_2$. Similarly, all cubic terms will be
controlled in $X^{2\e,0,2}$ (notice in the $s=\frac 1 4$ case, we had
$\frac 1 2$ spatial regularity to spare, which now becomes only an
$\e$), and inverting the linear operator yields yet another term which
we can safely incorporate into $w_2$. Strictly speaking, both $w_L$
and cubic terms have a part in the $2j''<k''$, but there one may trade
conormal regularity for spatial regularity and add them to $w_1$. 

Given the regularity loss, we cannot close an a priori estimate at
this level. However, one may go back to worst the bilinear interactions,
which occur in $\partial_x R(u,w)$, and perform
further substitution of $u$ by $w$; the price to pay is another bunch
of harmless cubic terms, and new bilinear interactions
$K\partial_x R(w,w))$, where $K$ stands for a
renormalization operator like in Proposition \ref{prop3.4}, arising
from the gauge. In the case
$s=\frac 1 4$, there was nothing to gain from such a substitution, as
$\partial_x R(w,w)\in X^{\frac 1 2,-\frac 1 2,1}$,
and the $\frac 1 4$ gain would be lost when applying $K$. Let us check
what estimates we may get, recalling that $R=R_1+R_2+R_3$:
\begin{itemize}
\item on $\partial_x R(u,w_1)$, we do not perform
  substitution, but the usual product rules give 
$$X^{2\e+\frac 1
  4+\e-\frac 1 4,-\frac 1 2,1}$$
 which will close;
\item on $\partial_x R(u,w_2)$, we do not perform
  substitution on $R_1$ and $R_2$, but the usual product rules give $X^{\e-\frac 1
  4+2\e-\frac 1 4+\frac 1 2,-\frac 1 2,1}$ which will close (we used
  $b=\frac 3 4$ on $w_2$); on the remaining (``worst'') term $R_3$, we need
  to substitute.
  \begin{enumerate}
  \item  $\partial_x R_3(w_1,w_2)$: this worst part is $X^{\e+\frac 1 4+\e,-\frac 1 2,1}$, where we
  took advantage of $b=\frac 1 4$ on $w_2$, and this term will close
  as well, losing $1/4$ with $K$ (no spare spatial regularity);
\item   $\partial_x R_2(w_2,w_2)$:  assume one of the two factors
  is such that $j''\sim j<k$ (or $j''\sim j'<k'$), then it will be $X^{2\e,\frac 1 2}$
  by trading conormal for spatial, and this term will be
  $X^{4\e,-\frac 1 2}$. The gauge does not yield any loss when applied
  to such a localized ($2j''\sim k''$) term : in fact, in the term
  which yields the $1/2$ loss when $K$ acts on $X^{0,1}$, one may use
  the $1/2$ remaining conormal regularity (which we threw away) to
  compensate for the spatial loss. Hence we close as well. Now, if
  both factors are such that $k<j''$ and $k'<j''$, then this term
  actually vanishes by support considerations (it requires $k>2j''$ !).
\item   $\partial_x R_3(w_3,w_2)$: using
  $b'=\frac 1 4$ on $w_3$, we get $X^{2\e-\frac 1 8+\frac 1
  2+2\e-\frac 1 8,-\frac 1 2,1}$ which is tight.
  \end{enumerate}
Finally, we substitute on the whole $\partial_x R(u,w_3)$:
\item  consider $\partial_x R(w_1,w_3
)$: 
 \begin{enumerate}
  \item the easiest part $R_1(w_1,w_3)$ is $X^{2\e+\frac 1 4+2\e-\frac 1 8+\frac 1
  2,-\frac 1 2,1}$ which will close;
\item the next one  $R_2(w_1,w_3)$ will be $X^{2\e+\frac 1 4+2\e-\frac 1 8,-\frac 1 4,1}$, and $K$
will only lose $\frac 1 8$ which is tight.
\item the third part  $R_3(w_1,w_3)$ is $X^{\e+\frac 1 4+2\e-\frac 1 8+\frac 1 2,-\frac 1 2,1}$, where we
  took advantage of $b=\frac 1 4$ on $w_3$, and this term will close;
  \end{enumerate}
\item  consider $\partial_x R(w_2,w_3
)$: 
 \begin{enumerate}
  \item the easiest part  $R_1(w_2,w_3)$ is $X^{2\e-\frac 1 8+2\e-\frac 1 8+\frac 1
  2,-\frac 1 2,1}$ which will close;
\item the second part  $R_2(w_2,w_3)$ vanishes by support considerations (as $w_3$ is
  localized at $k\sim 2j\sim 2j''$ is does not contribute: we are in a
  situation where $k''<2j''$ but actually one may take $k''<2j''-100$
  in the bilinear estimate which yields $w_2$).
\item the third part  $R_3(w_1,w_3)$ is $X^{2\e-\frac 1 8+2\e+\frac 3 8,-\frac 1 2,1}$, where we  took advantage of $b=\frac 1 4$ on $w_3$, and this term will close;
  \end{enumerate}
\item  consider $\partial_x R(w_3,w_3
)$: 
 \begin{enumerate}
  \item the easiest part  $R_1(w_3,w_3)$ is $X^{2\e-\frac 1 8+2\e-\frac 1 8+\frac 1
  2,-\frac 1 2,1}$ which will close;
\item the next one,  $R_2(w_3,w_3)$, vanishes for the same reasons as in the previous case.
\item the third part  $R_3(w_3,w_3)$ is $X^{2\e-\frac 1 8+2\e+\frac 3 8,-\frac 1 2,1}$, where we  took advantage of $b=\frac 1 4$ on $w_3$, and this term will close;
  \end{enumerate}
\end{itemize}
It remains to to deal with the gauge: the operator $K\equiv
T_{\exp(\int^x u)}$ acting on $\phi$ yields a term $T_u \partial_x
\phi$, in which we substitute $w$ to $u$ (actually, similar terms
arise as well in the other factors, unlike before). The worse terms occur in the
form of 
$$
\sum_{l<j-2} \Delta_l(\sum_{l\lesssim k} \Delta_k F^{-1} \Delta_k F
\Delta_l u) \partial_x \Delta_j \phi,
$$
where $F$ is the imaginary factor. But given $FF^{-1}=1=T_F F^{-1}+
T_{F^{-1}} F+R(F,F)$, we may substitute the remainder term $R(F,F)$ to
  reduce to estimating
$$
\| \Delta_l F \Delta_l u \partial_x \Delta_j \phi\|_2\lesssim 2^{-l}
2^{\frac l 4} 2^{\frac l 4} 2^{\frac j 2},
$$
using the smoothing estimate on $u$ and the maximal function estimate
on $\partial_x F$ and $u$. This allows to close the gauge as well.
\begin{rem}
  It should be clear for the previous outline that whenever $s>\frac 1
  8$, one does {\bf not} need to iterate the bilinear interactions
  twice, but that an appropriate use of the gauge estimate in the
  $b=\frac 1 4$ case is enough to close, combined with the fact that only an estimate on $u\in
  X^{\sigma,\frac 1 4,1}$ is required for all ``worst'' terms
  $R_1,R_2,R_3$. Note that all Duhamel terms in $w$ end up in $X^{2\epsilon,\frac 1
  2,1}$. When $\epsilon$ is actually $\frac 1 8$, this translates into
  a regularity gain.
\end{rem}
We now turn to uniqueness: recall the equation after renormalization
by the ``worst'' solution $v$ (denoting by $K_v$ the gauge operator) is
$$
\partial_t \omega+H\partial^2_x \omega+T_{\partial_x v} \omega+K_v
\partial_x(T_{\delta} u+R(2u,\delta)-R(\delta,\delta))=0.
$$
Note that as long as $s\geq 0$, assuming the solution to be in $X^{s-\frac 1
  4,\frac 1 2,1}$, the gauge transform maps now $X^{\sigma,\pm \frac 1
  2,1}$ to $X^{\sigma-\frac 1 4-\frac 1 2 (s-\frac 1 4),\pm\frac 1
  2,1}$ and we may consider $\omega\in X^{-\frac 1 2,\frac 1 2,1}$ as
in the range $s>\frac 1 4$.

The worst term appears to be $K_v\partial_x T_{\delta} u$, on which a
complete substitution produces a term $K_\delta \partial_x T_{\omega}
w$, where $w$ is the renormalized solution $K_u u$. Unlike on the
other paraproduct, the regularity of the low frequencies is $s'=-\frac
1 2 $ for which we barely recover the front derivative : however, the
output before the gauge is then $X^{s+s',-\frac 1 4,1}$ and the gauge
loss is slightly better; hence, from
$w\in X^{s,\frac 1 2,1}$, we obtain a condition $s-\frac 1 8+\frac 1
4 (s-\frac 1 4)>0$, which is exactly $s>3/20$. Next, $R(2u,\delta)$
requires exactly the same treatment (full substitution) with the same
result. The two remaining terms (up to cubic terms) $T_{\partial_x v} \omega$ and
$\partial_x R(\delta,\omega)$ are easier, and the same goes for all
cubic terms (though some become actually quartic through frequency
decomposition of the derivative of the exponential from the gauge).
\appendix

\section{Bilinear estimates}
\subsection{Conormal dyadic blocs products}
\label{sec:cdbp}
Recall that we defined a localization w.r.t. $\xi$ and
$\tau-\xi|\xi|$, see \eqref{eq:defbloc}. Now set
\begin{equation}
  \label{eq:defblocpm}
  \Delta^{\pm}_{jk} u(x,t) =  \mathcal{F}_{\tau,\xi}^{-1}
  (\psi^{\pm}_{jk}(\tau,\xi)\mathcal{F}_{t,x}( u)),
\end{equation}
where we recall that
$$
\psi^{\pm}_{jk}(\tau,\xi)=\chi_{\log |\xi|\sim j}\chi_{\log |\tau\mp
  \xi^2|\sim k},
$$
where $\sim$ means equivalent (except for $j=-1$ or
$k=-1$ for which it means $\lesssim$). We define two spaces which are
the restriction to positive and negative spectrum of $X^{s,b,q}$.
\begin{definition}
  Let  $u(x,t)\in \mathcal{S}'(\mathbb{R}^{n+1})$, $s,b\in \R$ and
  $1\leq q\leq +\infty$. We say
that
 $u\in X_\pm^{s,b,q}$ if and only if $\supp \mathcal{F}_{t,x} u
  \subset \{\pm \xi\geq 0\}$ and for all $j\geq -1$,
\begin{equation}
  \label{eq:xsb}
   \|\Delta^\pm_{jk} u\|_{L^2_{t,x}} \lesssim
  2^{-js-kb} \e_{jk},\,\,\,(\e_{jk})_{jk}\in l^q.
\end{equation}
\end{definition}
\begin{rem}
If it were not for the $\pm \xi\geq 0$ restriction, we
would have defined the usual Schr\"odinger Bourgain spaces. As such,
all the estimates we need in $X^{s,b,q}$ can be deduced from estimates
on products of functions in $X^{s,b,q}_{\pm,T}$, by reducing to dyadic
pieces and sorting out all possible sign combinations. Therefore, all
subsequent estimates could be retrieved, one way or another, from the
existing literature, see \cite{kpvnls1d,Taomulti,CDKS}. We elected to
give a self-contained proof in order to streamline the reading and
highlight as best as possible what the optimal estimate is, depending
on the frequencies constraints we set.
\end{rem}
In our setting we are dealing with the Benjamin-Ono equation where
$u_0$ (and consequently $u$) is real-valued: hence, we only need to estimate
its positive spectrum part, $\mathcal{F}^{-1}_{\xi}(\chi_{\xi\geq 0} \hat
u(\xi))$, to recover $u$. 

Later on, we will be interested in $P^+( P^\pm v P^+ u)$,
where $P^\pm$ are the
spectral projectors on positive/negative frequencies. From the
discussion above, we are reduced to estimating
$$
 \Delta^{+}_{j''k''}( \Delta^{\pm}_{j'k'} u \Delta^+_{jk} v),
$$
while knowing a priori that $j'<j$, due to the outer $P^+$.

We set $j^\flat=\min(j,j',j'')$ and $j^\sharp=\max(j,j',j'')$ and the
remaining middle one is $j^\natural$. Similarly with $k$, we have
$k^\flat\leq k^\natural\leq k^\sharp$. In the next two lemmata we set
$\|\cdot\|=\|\cdot \|_{L^2_{t,x}}$.

\begin{lemme}[Sobolev]
  We have
  \begin{equation}
    \label{eq:sobolev}
    \|  \Delta^{+}_{j''k''}( \Delta^{\pm}_{j'k'} v \Delta^+_{jk}
u)\| \lesssim 2^{\frac{j^\flat} 2+\frac {k^\flat} 2} \|\Delta^{\pm}_{j'k'} v\| \| \Delta^+_{jk}
u)\|.
  \end{equation}
\end{lemme}
\dem The product can be written as a convolution in $(\tau,\xi)$,
which is then localized according to $\Delta^+_{j'',k''}$. We may then
use Bernstein inequalities in the ``right'' directions using the
support sizes and the shape of the boxes. The proof of the next lemma
implicitly contains this one, so we do not give any details.\cqfd\\
In some situations, we can do better than Sobolev inequalities: the
usual (spatial) paraproduct splitting implies the following relations between the indices for which our function does not vanish by support
considerations:
\begin{itemize}
\item $j<<j'$ and $j'\sim j''$, hence $j^\flat=j$ and $j^\sharp\sim
  j'\sim j''$.
\item $j'<<j$ and $j\sim j''$, hence $j^\flat=j'$ and $j^\sharp\sim
  j\sim j''$.
\item $j'' << j$ and $j\sim j'$, hence  $j^\flat=j$ and $j^\sharp\sim
  j'\sim j$.
\item $j\sim j'\sim j''$, hence $j^\flat\sim j^\sharp$.
\end{itemize}
\begin{lemme}[Conormal regularity]
 \begin{itemize}
  \item   Let $k''=k^\sharp$. if $j''<<j\sim j'$, 
  \begin{equation}
    \label{eq:+++2}
    \|\Delta^+_{j''k''}(\Delta^+_{j'k'} v \Delta^+_{jk} u)\| \lesssim
    2^{\frac{k^\flat}2+\frac{k^\natural-j^\natural} 2}
    \|\Delta^+_{j'k'} v\| \| \Delta^+_{jk} u)\|,
  \end{equation}
and
  \begin{equation}
    \label{eq:+-+2}
    \|\Delta^+_{j''k''}(\Delta^-_{j'k'} v \Delta^+_{jk} u)\| \lesssim
    2^{\frac{k^\flat} 2+\frac{k^\natural-j^\flat} 2}
    \|\Delta^-_{j'k'} v\| \| \Delta^+_{jk} u)\|.
  \end{equation}
If $j\sim j'\sim j''$, 
  \begin{equation}
    \label{eq:+++21}
    \|\Delta^+_{j''k''}(\Delta^+_{j'k'} v \Delta^+_{jk} u)\| \lesssim
    2^{\frac{k^\flat}2+\frac{k^\natural} 4}
    \|\Delta^+_{j'k'} v\| \| \Delta^+_{jk} u)\|.
  \end{equation}
  \item   Let $k=k^\sharp$. If $j\sim j'\sim j''$,
  \begin{equation}
    \label{eq:+-+0}
    \|\Delta^+_{j''k''}(\Delta^-_{j'k'} v \Delta^+_{jk} u)\| \lesssim
    2^{\frac{k^\flat} 2+\frac{k^\natural} 4}
    \|\Delta^-_{j'k'} v\| \| \Delta^+_{jk} u)\|.
\end{equation}
  \item In all remaining cases, most notably including $j'<<j\sim
  j''$ irrespective of $k^\sharp$,
\begin{equation}
    \label{eq:+-+1}
    \|\Delta^+_{j''k''}(\Delta^\pm_{j'k'} v \Delta^+_{jk} u)\| \lesssim
    2^{\frac{k^\flat} 2+\frac{k^\natural-j^\natural} 2}
    \|\Delta^\pm_{j'k'} v\| \| \Delta^+_{jk} u)\|.
  \end{equation}
  \end{itemize}
\end{lemme}
\begin{rem}
Recall $\tau-\xi^2-(\tau-\sigma-(\xi-\eta)^2)-(\sigma\mp\eta^2)=(\xi-\eta)^2\pm\eta^2-\xi^2$,
and assume $k^\sharp=k''$. The left handside is clearly bounded by
$2^{k''}+2^{k'}+2^{k}=O(2^{k^\sharp})$. On the right handside, we consider
different cases:
\begin{itemize}
\item if $j''<<j\sim j'$, then $\xi-\eta\sim-\eta$. In the $+++$ case,
  we get $2 j^\sharp\lesssim k^\sharp$. In the $+-+$ case, expanding
  leads to $\xi \eta$ and therefore $j^\flat+j^\sharp\lesssim
  k^\sharp$;
\item if $j'<< j\sim j''$, similarly one gets $ j^\flat+j^\sharp\lesssim
  k^\sharp$ in both cases;
\item if $j<< j'\sim j''$, one has to switch $u$ and $v$ and we are
  back in the previous case (note that for BO, this never happens in
  the $+-+$ case);
\item if $j\sim j'\sim j''$, we have $\xi\sim \xi-\eta\sim \pm \eta$
  and again $2^{j^\sharp}\sim j^\flat+j^\sharp\lesssim k^\sharp$.
\end{itemize}
Finally, if $k^\flat,k^\natural<<k^\sharp$, then the left handside is actually
exactly $O(2^{k^\sharp})$, in which case we get either $2j^\sharp \sim
k^\sharp$ ($+++$ case) or $j^\flat+j^\sharp\sim k^\sharp$($+-+$ case).
\end{rem}
Now a simple duality argument reduces the study to the case $k''=k^\sharp$:
\begin{itemize}
\item if $k=k^\sharp$, then $
\langle  \Delta^{+}_{j''k''}( \Delta^{\pm}_{j'k'} u \Delta^+_{jk}
v),\varphi\rangle = \langle  \Delta^+_{jk}
v, \Delta^+_{jk}( \Delta^{\mp}_{j'k'} \bar u
\Delta^{+}_{j''k''}\varphi)\rangle 
$
 and we have to deal with $\Delta^+_{jk}( \Delta^{\mp}_{j'k'} \bar u
\Delta^{+}_{j''k''}\varphi)$;
\item if $k'=k^\sharp$, by the same reasoning, one has to deal with
$\Delta^\pm_{j'k'}( \Delta^{-}_{jk} \bar v
\Delta^{+}_{j''k''}\varphi)$.
\end{itemize}
\begin{rem}
  Notice now we got all possible sign combinations. However, from the
  symmetry with respect to $\xi=0$, the $--+$ case is no different
  from $+-+$. To
  summarize, we will deal with $+\pm+$ with $k''=k^\sharp$.
\end{rem}
 We should now set $k$ and $k'$: if we have the $+++$ case, by symmetry
we may choose $k'\leq k$ and we are left with the paraproduct
trichotomy; otherwise we get 2 separate cases. All in all, we get
$3+3\times 2$ possible cases. Set
$$
I(\tau,\xi)=\psi^+_{j''k''}(\tau,\xi)\int f(\tau-\sigma,\xi-\eta)
 g(\sigma,\eta)d\sigma d\eta,
$$
with $f,g$ being the space-time Fourier transforms of two dyadic
blocks of $u,v$. We have 2 cases:
\begin{itemize}
\item if $g=\mathcal{F}_{t,x}\Delta^+_{j'k'} v$ ($+++$ case), $
\tau-\xi^2=\tau-\sigma-(\xi-\eta)^2+\sigma-\eta^2+2\eta(-\xi+ \eta),
$
 and we define $F^+(\xi,\eta)=2\eta(-\xi+
\eta)=(\xi-\eta)^2+\eta^2-\xi^2$.
\item
 if $g=\mathcal{F}_{t,x}\Delta^-_{j'k'} v$ ($+-+$ case),
 $
\tau-\xi^2=\tau-\sigma-(\xi-\eta)^2+\sigma+\eta^2+((\xi- \eta)^2-\eta^2-\xi^2)
$
 and $F^-(\xi,\eta)=(\xi-\eta)^2-\eta^2-\xi^2$.
\end{itemize}
We now proceed as follows.
\subsubsection{Case $j''<< j^\sharp$: remainder}
This is the product minus paraproducts and $j^\flat\sim j^\sharp$
terms. We split the support
of $\hat g$ into non-overlapping intervals of size $2^{j''}$ (recall
$j''=j^\flat$), and the support of $f$ will be
forced into a interval of comparable length but opposite with respect
to the $\xi=0$ frequency. As such, one has to deal with
$$
I=\sum_Q \phi^+_{j''k''} f_{-Q} \star g_{Q}=\sum_Q I_Q.
$$
We will later write (using $L^2_\xi$ orthogonality of $(f_Q)_Q,(g_Q)_Q$)
\begin{align*}
  \|I\| \leq & \sum_Q \|I_Q\|  \leq  \sum_Q M \|f_{-Q}\|\| g_Q\| \leq  M \bigl(\sum_Q \|f_{-Q}\|^2\bigr)^\frac 1 2 \bigl(
   \sum_Q \|g_Q\|^2\bigr)^\frac 1 2\\
    \lesssim & M  \|\sum_Q f_{-Q}\| \| \sum_Q g_Q\| \lesssim   M \|f\|\|  g\|,
\end{align*}
where $M$ will be obtained when evaluating $I_Q$ in the next step.
\begin{lemme}
 Introduce a parabolic level set
decomposition, where $l\in \Z$,\begin{equation}
  \label{eq:parabolic}
  \chi_l (\xi,\eta)=\chi_{l 2^{k^\natural} < F^\pm(\xi,\eta)< (l+1)
  2^{k^\natural}} \text{ and } I^l_Q= \psi^+_{j''k''}\int \chi_l
  \hat f_{-Q}(\xi-\eta,\tau-\sigma) \hat g_{Q}(\eta,\sigma) d\sigma
  d\eta.
\end{equation}
  Then the family $(I^l_Q)_l$ is almost orthogonal in $L^2_{\tau,\xi}$.
\end{lemme}
\dem Assume $(\xi,\tau)$ is in the support of $I^l_Q$ and $I^m_Q$. We
shall prove that $l\sim m$. Certainly there exist $(\sigma_1,\eta_1)$
such that
$$
l2^{k^\natural}<F(\xi,\eta_1)<(l+1)2^{k^\natural} \text{ and } |\sigma_1 \mp
\eta_1^2|<2^{k'} \text{ and }  |(\tau-\sigma_1)-
(\xi-\eta_1)^2|<2^{k}.
$$
Similarly, there exist $(\sigma_2,\eta_2)$
such that
$$
m2^{k^\natural}<F(\xi,\eta_2)<(m+1)2^{k^\natural}
\text{ and } |\sigma_2 \mp
\eta_2^2|<2^{k'} \text{ and }  |(\tau-\sigma_2)-
(\xi-\eta_2)^2|<2^{k}.
$$
Recall that, for any $(\sigma,\eta)$, we have
$\tau-\xi^2=\tau-\sigma-(\xi-\eta)^2+\sigma\mp\eta^2+F^\pm(\xi,\eta)$,
hence
$$
F^\pm(\xi,\eta_1)-F^\pm(\xi,\eta_2)=\tau-\sigma_2-(\xi-\eta_2)^2+\sigma_2\mp\eta_2^2-(\tau-\sigma_1-(\xi-\eta_1)^2+\sigma_1\mp\eta_1^2).
$$
Given the self-imposed bounds, the right handside is bounded by $2(2^{k}+2^{k'})=O(2^{k^\natural})$. This
bounds $|F^\pm(\xi,\eta_1)-F^\pm(\xi,\eta_2)|$, but by virtue of the
 lower bounds on $F(\xi,\eta_1)$ and $F(\xi,\eta_2)$, this difference is bounded from below by $(l-m-1) 2^k$ if $m<l$,
hence $l\sim m$.\cqfd

We can now perform Cauchy-Schwarz,
\begin{equation}
  \label{eq:CSxitau}
  I^l_Q(\xi,\tau)\leq \left(\int_{A(Q,j',k',\xi,\tau)}\chi_l(\xi,\eta) d\sigma
    d\eta\right)^{\frac 1 2} \left(\int_{\sigma,\eta} \chi_l(\xi,\eta)
    |f_{-Q}|^2 |g_Q|^2\right)^\frac 1 2,
\end{equation}
where 
$$
A(Q,j',k',\xi,\tau)=\{ (\sigma,\eta),\,\,\eta\in Q,\,\,|\eta|\sim
2^{j'},\,\,|\sigma\mp \eta^2|\sim
2^{k'},\,\,|\tau-\sigma-(\xi-\eta)^2|\sim 2^k\}.
$$
We need to bound $M^2=\int_A \chi_l$; we start with integration w.r.t. $\sigma$: one cannot improve upon what the support size
yields; namely, at fixed $(\eta,\xi,\tau$), the interval is at most
$\min(2^{k'},2^k)=2^{k^\flat}$. Now for the integration w.r.t. $\eta$, one has
an integral over a domain $B$, with fixed $(\xi,\tau)$,
$$
B=\{ \eta\in Q,,\,\,|\eta|\sim
2^{j'},\,\,l2^{k^\natural}<F^\pm(\xi,\eta)<(l+1)2^{k^\natural}\}.
$$
We have to deal separately with the two cases, knowing that at any
rate,  $\eta$ is in an interval of size $2^{j''}$ (which yields
Sobolev !):
\begin{itemize}
\item either $F^-$, namely $l2^{k^\natural}< -2\xi\eta< (l+1) 2^{k^\natural}$. As such,
  $\eta$ is in an interval of
  size $2^{k^\natural}/|\xi|\sim 2^{k^\natural-j''}$ and we
  obtain 
$
M\lesssim 2^{\frac{k^\flat} 2+\frac{k^\natural-j^\flat} 2}.
$

\item either $F^+$, namely $l2^{k^\natural} < 2\eta^2-2\xi\eta< (l+1)
  2^{k^\natural}$. Note that (as $|\eta|\sim 2^{j'}>> 2^{j''}\sim
  |\xi|$), $2^{2j'}\lesssim F^+$, and therefore $|\xi|^2<<l2^{k^\natural}$. We have
  \begin{equation*}
    |\eta_{\max}-\eta_{\min}|\leq  \sqrt{(l+1)
     2^{{k^\natural-1}}+\frac{\xi^2} 4}-\sqrt{l
     2^{{k^\natural}-1}+\frac{\xi^2} 4}
\leq  \frac{2^{k^\natural}}{\sqrt{(l+1)
     2^{{k^\natural}}+\frac{\xi^2} 4}}\lesssim \frac{2^{{
     k^\natural} }}{\sqrt{2^{2j'}}}.
  \end{equation*}
Therefore $\eta$ is  in an interval of
  size at most $2^{\frac {k^\natural-j} 2}$, and we
  obtain 
$
M\lesssim 2^{\frac{k^\flat} 2+\frac{k^\natural-j^\natural} 2}.
$
\end{itemize}
>From \eqref{eq:CSxitau}, we have
\begin{equation*}
\|I_Q\|^2\leq \sum_l \int_{\xi,\tau} I^2_l
     \leq  M^2 \int_{\xi,\tau,\sigma,\eta}
     |f_{-Q}|^2(\xi-\eta,\tau-\sigma) |g_Q|^2(\eta,\sigma )
    \leq  M^2 \|f_{-Q}\|^2 \|g_Q\|^2.
\end{equation*}
\subsubsection{Case $j^\flat\sim j^\sharp$}
We write (with $\eta=\xi/2+\lambda$)
$$
I(\tau,\xi)=\phi^+_{j''k''} \int f(\tau-\sigma,\frac \xi
2-\lambda) g(\sigma,\frac \xi 2+\lambda) d\sigma d\lambda.
$$
Now $F^+=2\lambda^2-\xi^2/2$ and $F^-=-2\lambda\xi-\xi^2$. Again,
define level sets  with $l2^{k^\natural}<F^\pm<(l+1)
2^{k^\natural}$ and with the characteristic function $\chi_l$. We are
led to $M^2=\int_{\sigma,\lambda} \chi_l$, and a situation which is
very similar to the reminder situation.
\begin{itemize}
\item If we have $F^-$, the $\lambda$ interval will be of size at most
  $2^{k^\natural}/|\xi|$, which yields $M\lesssim 2^{\frac{k^\flat}
  2+\frac{k^\natural-j^\natural} 2}$.
\item If we have $F^+$, the situation is worse: unlike in previous
  cases, $F^+$ may very well be close to zero when
  $\lambda$ varies. We have $l2^{k^\natural-1}+\frac{\xi^2}
  2<\lambda^2<(l+1)2^{k^\natural-1}+\frac{\xi^2}
  2$. Either $l=0$ and $\lambda$ varies in an interval of size
  $2^{\frac{k^\natural} 2}$, or $l\neq 0$ and either $\frac{\xi^2}4+l 2^{k^\natural-1}<2^{k^\natural}$,
and we conclude 
    \begin{equation*}
    |\eta_{\max}-\eta_{\min}|\leq  \left|\sqrt{(1+l)
     2^{{k^\natural-1}}+\frac{\xi^2} 4}-\sqrt{l
     2^{{k^\natural-1}}+\frac{\xi^2} 4}\right|
 \lesssim  2^{\frac{
     k^\natural} 2 }
  \end{equation*}
or $\frac{\xi^2}4+l 2^{k^\natural-1}>2^{k^\natural}$ and we conclude by
    \begin{equation*}
    |\eta_{\max}-\eta_{\min}|\leq \frac{2^{k^\natural}}{\sqrt{l
     2^{{k^\natural-1}}+\frac{\xi^2} 4}}
 \lesssim 2^{\frac{
     k^\natural} 2 }.
  \end{equation*}
Finally, we get $
M\lesssim 2^{\frac{k^\flat} 2+\frac{k^\natural} 4}.
$
\end{itemize}
\subsubsection{Case $j'<<j''\sim j\sim j^\sharp$: paraproduct}

We still proceed with the same computation as before (without
introducing the intervals $Q$). Here, we have $|\eta|\sim 2^{j'}$.

\begin{lemme}
We  consider again the parabolic level set
decomposition, 
$$
  \chi_l (\xi,\eta)=\chi_{l 2^{k^\natural} < F^\pm(\xi,\eta)< (l+1)
  2^{k^\natural}} \text{ and }
I^l= \psi^+_{j''k''}\int \chi_l 
f(\xi-\eta,\tau-\sigma) g (\eta,\sigma) d\sigma d\eta.
$$
  The family $(I^l)_l$ is almost orthogonal in $L^2_{\tau,xi}$.
\end{lemme}
the proof is word for word identical to the previous case: we never
used any support condition in $\xi$ or $\eta$.\cqfd

Perform Cauchy-Schwarz,
$$
I^l(\xi,\tau)\leq \left(\sup_{\xi,\tau} \int_{B(Q,j',k',\xi,\tau)}\chi_l(\xi,\eta) d\sigma
    d\eta\right)^{\frac 1 2} \left(\int_{\eta,\sigma} \chi_l(\xi,\eta)
    |f|^2\,
    |g|^2\right)^\frac 1 2,
$$
where 
$$
B(Q,j',k',\xi,\tau)=\{ (\sigma,\eta),\,\,|\eta|\sim 2^{j'}\,\,|\sigma\mp \eta^2|\sim
2^{k'},\,\,|\tau-\sigma-(\xi-\eta)^2|\sim 2^k\}.
$$
Again, at fixed $(\eta,\xi,\tau$), the interval in $\sigma$ is at most
$\min(2^{k'},2^k)=2^{k^\flat}$. For $\eta$, we have 
$$
B=\{ \eta\in Q,,\,\,|\eta|\sim
2^{j'},\,\,l2^{k^\natural}<F^\pm(\xi,\eta)<(l+1)2^{k^\natural}\}.
$$
Split the two cases, knowing $\eta$ is at most in an interval of size
$2^{j'}$ (which is Sobolev):
\begin{itemize}
\item either $F^-$, namely $l2^{k^\natural}< -2\xi\eta< (l+1) 2^{k^\natural}$. As such,
  $\eta$ is in an interval of
  size $2^{k^\natural}/|\xi|\sim 2^{k^\natural-j}$. This yields again
$
M\lesssim 2^{\frac{k^\flat} 2+\frac{k^\natural-j} 2}.
$
\item either $F^+$, namely $l2^{k^\natural} < 2\eta^2-2\xi\eta< (l+1)
  2^{k^\natural}$. Now, as $|\eta|\sim 2^{j'}<< 2^{j}\sim
  |\xi|$), $|F^+|\sim 2^{j+j'+1}$. 

Assume $\eta\xi<0$ (which never happens for BO)
  then $l>0$ and 
    \begin{equation*}
    |\eta_{\max}-\eta_{\min}|\leq  \sqrt{(l+1)
     2^{{k^\natural-1}}+\frac{\xi^2} 4}-\sqrt{l
     2^{{k^\natural}-1}+\frac{\xi^2} 4}
   \leq  \frac{2^{k^\natural}}{4\sqrt{l
     2^{{k^\natural}-1}+\frac{\xi^2} 4}}\lesssim {2^{\frac{
     k^\natural-j} 2  }}.
  \end{equation*}
 and $\eta$ is again in an interval of
  size at most $2^{\frac{k^\natural-j^\natural} 2}$ and
$
M\lesssim 2^{\frac{k^\flat} 2+\frac{k^\natural-j^\natural} 2}$.

Assume now $\eta\xi>0$, this forces $k^\natural\lesssim j+j'$ or the $\eta$-set
is empty, and $-l\sim 2^{j+j'+1-k^\natural}$. Call $m=-l>0$:
$$
-m2^{k^\natural} <2\eta^2-2\xi\eta< (-m+1)  2^{k^\natural},
$$
and remark that $(-m+1)2^{k^\natural-1}+\frac{\xi^2}4 \geq
(\eta-\frac{\xi}2)^2\sim 2^2j$. Hence
    \begin{equation*}
    |\eta_{\max}-\eta_{\min}|\leq  
   \leq  \frac{2^{k^\natural}}{4\sqrt{(-m+1)
   2^{{k^\natural}-1}+\frac{\xi^2} 4}}\lesssim {2^{\frac{
     k^\natural-j} 2  }}.
  \end{equation*}
and
$
M\lesssim 2^{\frac{k^\flat} 2+\frac{k^\natural-j^\natural} 2}$. Summing the $l$ pieces is the same as before (except we have no $Q$).
\end{itemize}
\subsubsection{Case $j''\sim j'\sim j^\sharp$: paraproduct}
One may just reverse the order of $f$ and $g$ in
the convolution to get the very same result as in the previous case:
this will be obvious for the $+++$ case, while for $+-+$, we have a shifted $F^-=2\xi\eta-2\xi^2$ to which the same computation applies.
\subsection{Product estimates}\label{sec:prod}
Recall the paraproduct
decomposition: $
uv=T_v u+T_u v+R(u,v),
$
 where we define
$$
T_v u=\sum_j S_{j-1} v \Delta_j u,\,\,\text{ and }\,\, R(u,v)=\sum_{|j-j'|\leq 1}
\Delta_j v \Delta_{j'} u.
$$
In this section we prove all important bilinear estimates. We first state the estimates required to obtain the existence and uniqueness in $H^s, s>1/4$.
\begin{proposition} Let $ u \in X^{ s , b,q }$ and $v \in X^{s' ,
    \frac 1 2 ,1 }$. We consider the mapping $w=P^+(P^{\pm}v\, P^+u)$.
  \begin{itemize}
  \item Assume moreover that $s+s'\geq -\frac 1 2$, we have
    \begin{itemize}
    \item if $b=\frac 1 2$, $q=1$:
    \begin{equation}
      \label{eq:BH1}
\partial_x P^+(R(P^{\pm}v\,P^+u))\in X^{s+s',-\frac 1 2,1}.
    \end{equation}
\item  if $s'=0$, $s>\frac 1 4$ $b>\frac 1 2$, $q=2$:
    \begin{equation}
      \label{eq:BH2}
\partial_x P^+(R(P^{\pm}v\,P^+u))\in X^{s,b-1,2}+X^{s,\frac 1 2,1}.
    \end{equation}
    \end{itemize}
\end{itemize}
\item Assume moreover that $s'+\frac 1 2\leq 0$, then
    \begin{itemize}
    \item if $b=\frac 1 2$, $q=1$:
    \begin{equation}
      \label{eq:BH3}
      P^+(T_{P^\pm v} P^+u)\in X^{s+s'+1,-\frac 1 2,1}.
    \end{equation}
\item  if $b>\frac 1 2$, $q=2$:
    \begin{equation}
      \label{eq:BH4}
      P^+(T_{P^\pm v} P^+u)\in X^{s+s'+1,b-1,2}+X^{s+s'+1,-\frac 1 2,1}.
    \end{equation}
    \end{itemize}
\end{proposition}
We shall also need the following refinement when dealing with the unconditional uniqueness in $L^\infty_t( H^{1/2})$.
\begin{proposition}
Let $ u \in X^{ s , \frac 1 2,2}$ and $v \in X^{-1 ,
    \frac 1 2 ,2 }$. Then we have
 \begin{equation}
      \label{eq:BH4bis}
      P^+(T_{P^\pm v} P^+u)\in X^{s,- \frac 1 2,1}.
    \end{equation}
Let  $ u \in X^{ - \frac 1 2 , \frac 1 2,1}$ and $v \in X^{0 ,
    \frac 1 2 ,2 }\cap X^{0,\frac 1 4,1}$. Then we have  
\begin{equation}
      \label{eq:BH1bis}
\partial_x P^+(R(P^{\pm}v\,P^+u))\in X^{- \frac 1 2,-\frac 1 2,1}.
    \end{equation}
\end{proposition} 
\dem We expand all functions dyadically, and are therefore left with estimating
$$\sum_{j'', k''} \Delta^+_{j''. k''} (\sum_{j',j, k',k}\Delta^{\pm}_{j',
  k'}v\times \Delta^+_{j,k}u).
$$ 
We set
$$
\|\Delta^+_{j,k}u\|\leq \beta_{jk}
2^{-js-kb},\,\,\,\|\Delta^{\pm}_{j',k'}v\|\leq \alpha_{j'k'}
2^{-j's'-k'b'} \text{ with } \alpha,\beta \in l^1(j,k).
$$
Obviously, in the $+-$ case, the first constraint is $j'<j$ otherwise by support
condition (the $P^+$ in front of the product) it vanishes. In the $++$
case, both functions $u,v$ play identical role and therefore we
consider only two terms in the usual
spatial paraproduct decomposition (discarding $T_u v$). Given that for $b>\frac 1 2$, we obviously have
$X_{\pm}^{s,b,1}\hookrightarrow X_{\pm}^{s,\frac 1 2,1}$, we will
perform interpolation to recover (most of the) $q>1$ cases (which are
useful for propagation of Sobolev regularity and unconditional uniqueness). In most instances,
$b=b'=\frac 1 2$, but we will occasionally set $b$ or $b'$ to be
$0,\frac 1 4,$ or $\frac 3 4$.
\subsubsection{The remainder of the spatial paraproduct: $j'\sim j\sim
  j^\sharp$ and $j''=j^\flat$} Here it should be obvious that only the
  sum of the two regularities will be of importance, given the $j\sim
  j'$ condition. We accordingly set $\sigma=s+s'+\frac 1 2$, the case
  $\sigma=0$ being the borderline (worst) case, and the situation
  improving with $\sigma>0$ which provides better summability in
  several instances.\\
Start with the $-+$ case: we have $j''+j\lesssim
  k^\sharp$ from the dispersion relation; either $j''+j\sim k^\sharp$
  or $k^\sharp\sim k^\natural$. We consider 3 cases, depending on
  $k^\sharp$.
  \begin{enumerate}
  \item $k''=k^\sharp$: then we have  $2j''\lesssim j+j''\lesssim
  k''$. We deal with both $j''<<j$ and $j''\sim j$. Pick the conormal factor from \eqref{eq:+-+2},
we get 
$$
c_{j'',k''}\lesssim \sum_{k,k'\lesssim k''\,\,\,}\sum_{j''\lesssim j
  \lesssim k''-j''} 2^{\frac{k+k'} 2-\frac{j''} 2-j(s+s')-k'b'-kb}\alpha_{j,k'} \beta_{j,k}.
$$
Given $b,b'\geq \frac 1 2$, we simply use $\alpha\in l^\infty(k')$ and $\beta \in
l^\infty(k)$ to get $k''^2$ for the $k,k'$ sums, and $s+s'+\frac 1
2\geq 0$ to get, if $s+s'>0$,
$$
c_{j''k''}\lesssim 2^{-\sigma j''} (k''- 2j'') k''^2 \lesssim 2^{( \frac {1} 2-\epsilon)k''}
2^{-(1-4\epsilon+\sigma)j''} (k''- 2j'') 
 k''^2 2^{-\frac \epsilon 2 k''-\epsilon j''},
$$
where, as $2j''\lesssim k''$, we trade conormal regularity
for spatial regularity, and regain $k'',j''$ summability ($\epsilon>0$
is small). If $s+s'<0$, 
\begin{multline}
c_{j''k''}\lesssim 2^{-\frac{j''}2} (k''- 2j'') k''^2 2^{-(s+s')(k''-j'')} \\
\lesssim 2^{( \frac {1} 2-\epsilon)k''}
2^{-(1-4\epsilon+\sigma)j''} 2^{\sigma(2j''-k'')} (k''- 2j'') 
 k''^2 2^{-\frac \epsilon 2 k''-\epsilon j''}.
\end{multline}
 Therefore, the output is $X^{\frac 3 2+s+s'-4\epsilon,-\frac 1
  2+\epsilon,1}\chi_{2j''\lesssim k''}$.
\item $k=k^\sharp$ then $k>> 2j''$, and we have to sum over
  $k>k''$ , $k'<k$ and $j''<j<k-j''$. First, we deal with $j''<<j$: we
  have
$$
c_{j'',k''}\lesssim\sum_{k''\lesssim k,k'\lesssim k} \sum_{\,\,\,\,\,j''\lesssim
  j\lesssim k-j''} 2^{\frac{k^\flat} 2+\inf(\frac{j''}2,\frac{k^\natural-j^\natural} 2) -j(s+s')-kb-k'b'}\alpha_{j,k'}\beta_{j,k}.
$$
\begin{enumerate}
\item We restrict the sum to $k'<k''$: 
  \begin{itemize}
\item if $k''<< k$, $j\sim k-j''$, and conormal is better ($k''-j \lesssim j''$),
$$
c_{j'',k''}\lesssim 2^{\frac{k''} 2+j''(\frac 1 2+s+s')}
\sum_{k>\sup(k'',2j'')} 2^{-k(s+s'+b+\frac 1 2)} \beta_{k-j'', k}\sum_{k'<k''} 2^{k'(\frac 1 2-b')}\alpha_{k-j'',k'}.
$$
Now, if $\alpha,\beta \in l^\infty$ and $b'=\frac 1 2$, the last sum is less than
$k''$, and the first sum is finite because $s+s'+b\geq 0$,
$$
c_{j''k''}\lesssim  2^{\frac{k''} 2+\sigma j''}2^ {
  -(\sigma+b)\sup(k'',2j'')} k''\lesssim 2^{k''(\frac 1 2-\epsilon)}
2^{-j''(1+\sigma-4\epsilon)} 2^{-k''\frac \epsilon 2-j''\epsilon}k'', 
$$
where again we may trade conormal regularity for spatial regularity,
depending on $\sup(k'',2j'')$. Therefore, the output is $X^{\frac 3 2+s+s'-4\epsilon,-\frac 1
  2+\epsilon,1}$.
\item if $k\sim k''$, $j<k''-j''$, which means Sobolev is better:
$$
c_{j'',k''}\sim \sum_{k'\lesssim k''} \sum_{\,\,\,\,\,j''\lesssim
  j\lesssim k''-j''} 2^{\frac{k'} 2+\frac{j''}2 -j(s+s')-k''b-k'b'}\alpha_{j,k'}\beta_{j,k''}.
$$
With $b'=\frac 1 2$ and $\alpha\in l^\infty$, the first sum yields
$k''$. Either $s+s'\geq 0$ and we get
$$
c_{j'',k''}\lesssim 2^{\frac{j''} 2-bk''-j''(s+s')} k''(k''-2j''),
$$
or $s+s'<0$ and we get
$$
c_{j'',k''}\lesssim 2^{\frac{j''} 2-bk''+(k''-j'')(-s-s')} k'',
$$
where we only used $\beta \in l^\infty$. Assume $j''<k''/2$, trading conormal
regularity for spatial regularity yields
$$
c_{j'',k''}\lesssim 2^{-(s+s'+2b+\frac 1 2-\epsilon)j''}2^{(\frac 1
  2-\epsilon)k''}(k''-2j'') k'' 2^{-\frac \epsilon 4 k''-\frac
  \epsilon 2 j''},
$$
hence  the output is $X^{\frac 3 2+s+s'-2\epsilon,-\frac 1
  2+\epsilon,1}$.

If $j''>k'' /2$, from $k\sim k''$ and $j<<k-j''$ we get $k''\sim
  2j''$, and $j\sim j''$, hence this case will be treated later
  (recall we are in a situation where $j''<<j$).
\end{itemize}
\item Now we deal with $k''<k'<k$:
 \begin{itemize}
  \item if $j''<k'/2$, pick Sobolev, and
$$
c_{j'',k''}=2^{\frac {j''}2+\frac {k''} 2} \sum_{\sup(k'',2j'')<k}
\sum_{\,\,\,\,\,\,j''<j<k-j''} 2^{-j(s+s')-kb}
\beta_{j,k}\sum_{\sup(k'', 2j'')<k'<k} 2^{-k'b'}
  \alpha_{j, k'}.
$$
When $s+s'>0$, one gets (discarding any summability in $\alpha,\beta$)
$$
 c_{j'',k''}\lesssim 2^{\frac{k''} 2} 2^{-j''(s+s'-\frac 1 2)}
 2^{-(b+b')\sup(k'', 2j'')},$$
while when $s+s'<0$, $s+s'+b\geq 0$ and $\beta\in l^1$, we get
$$
c_{j''k''}\leq 2^{k''/2+\sigma j''-\sup(k'',2j'')(s+s'+b+b')}
$$
therefore, the output is  $X^{s+s'+2b+2b'-\frac 1 2-4\epsilon,-\frac 1
 2+\epsilon,1}$.
\item if $j''>k'/2$, 
\begin{multline}
c_{j'',k''}\lesssim 2^{\frac {k''} 2} \sum_{\sup(k'',2j'')<k} 2^{-kb}
\sum_{j''<j<k-j''} 2^{-j(s+s')}\beta_{j,k}\\
 \sum_{k''<k'<\inf(k,2j'')}
2^{-b'k'} \alpha_{j,k'} 2^{\inf(\frac{j''}2,\frac{k'-j} 2)}.
\end{multline}
which is non-zero only if $k''<2j''$. 
\begin{itemize}
\item either $k'<<k$ and $j\sim k-j''$, conormal is better,
$$
c_{j''k''}\lesssim 2^{\frac {k''} 2+j''(\frac 1 2+s+s')} \sum_{k>2j''}
2^{-k(s+s'+b+\frac 1 2)}\beta_{k-j'', k} \sum_{k''<k'<k} 2^{k'(\frac 1 2-b')}\alpha_{k-j'',k'}.
$$
Assume $\alpha,\beta\in l^\infty$, then, as $s+s'+b+\frac 1 2>0$,
$$
c_{j''k''}\lesssim 2^{\frac {k''} 2+\sigma j''}2^{-2j''(\sigma+b)}k'',
$$
and the output is  $X^{s+s'+2b+\frac 1 2-4\epsilon,-\frac 1
 2+\epsilon,1}$.

\item Either $k'\sim k$ and $j<< k-j''$. Then Sobolev is better,
$$
c_{j''k''}\lesssim 2^{(k''+j'')/2} \sum_{\sup(2j'',k'')<k} 2^{-k(b+b')} \sum_{j''<j<k-j''} 2^{-j(s+s')}\beta_{j,k}\alpha_{j,k}.
$$
Discarding summability over $\alpha,\beta$ and setting $b+b'=1$, either $s+s'>0$ and we get
$$
c_{j''k''}\lesssim  2^{(k''+j'')/2} 2^{-\sup(2j'',k'')}2^{-(s+s')j''},
$$
which is $X^{s+s'+\frac 3 2-4\epsilon,-\frac 1 2+\epsilon,1}$, 

or $s+s'<0$ and we get (recall $s+s'+1>0$)
$$
c_{j''k''}\lesssim  2^{\frac {k''}2}
2^{-\sup(2j'',k'')(1+s+s')}2^{\sigma j''},
$$
which yields again  $X^{s+s'+\frac 3 2-4\epsilon,-\frac 1
  2+\epsilon,1}$.
\end{itemize}
  \end{itemize}
We are left with $j\sim j'\sim j''$: we have
$$
c_{j''k''}\lesssim \sum_{k'',k'\leq k} 2^{\frac{k^\flat}
  2+\inf(\frac{j''}2,\frac{k^\natural} 4)-j''(s+s')-kb-k'b'}
\alpha_{j''k'}\beta_{j''k}.
$$
\begin{itemize}
\item Assume $k'<k''$:
  \begin{itemize}
  \item if $j''<\frac{k''}2$, Sobolev is better,
$$
c_{j''k''}\lesssim 2^{-(s+s'-\frac 1 2)j''} \sum_{k'<k''\leq k}
2^{-kb+(\frac 1 2-b')k'}
\alpha_{j''k'}\beta_{j''k}
$$
which, discarding summability, yields
$$
c_{j''k''}\lesssim  2^{-(s+s'-\frac 1 2)j''} 2^{-bk''}k'',
$$
and trading regularity, the output is $X^{s+s'+\frac 3 2-4\e,-\frac 1
  2+\e,1}$;
\item if $j''>\frac{k''}2$, conormal is better, (notice we may use
  $k''<2j''-C$ with a large $C$ in this term and the previous one, should it be helpful)
$$
c_{j''k''}\lesssim 2^{\frac{k''} 4 -(s+s')j''} \sum_{k'<k''<2j''<k}
2^{-kb+(\frac 1 2-b')k'}
\alpha_{j''k'}\beta_{j''k}
$$
which, assuming $\alpha,\beta\in l^1$, $b'=\frac 1 2$ is
$X^{s+s'+2b,-\frac 1 4,1}$. Notice that if $b'=\frac 1 4$, then the
output is $X^{s+s'+2b,-\frac 1 2,1}$.
  \end{itemize}
\item Assume $k''<k'$:
  \begin{itemize}
  \item if $j''<\frac{k'}2$, Sobolev is better,
$$
c_{j''k''}\lesssim 2^{-(s+s'-\frac 1 2)j''+\frac{k''}2} \sum_{\sup(k'',2j'')<k'\leq k}
2^{-kb-b'k'}
\alpha_{j''k'}\beta_{j''k}
$$
which, discarding summability, yields
$$
c_{j''k''}\lesssim  2^{-(s+s'-\frac 1 2)j''} 2^{\frac {k''} 2-(b+b')\sup(k'',2j'')}k'',
$$
and the output is $X^{s+s'+\frac 3 2-4\e,-\frac 1
  2+\e,1}$;
\item if $j''>\frac{k'}2$, conormal is better,
$$
c_{j''k''}\lesssim 2^{\frac{k''} 2 -(s+s')j''} \sum_{k''<k'<2j''< k}
2^{-kb+(\frac 1 4-b')k'}
\alpha_{j''k'}\beta_{j''k}
$$
which, assuming $\alpha,\beta\in l^1$, $b'=\frac 1 2$ is
$X^{s+s'+2b,-\frac 1 4,1}$. Notice again that if $b'=\frac 1 4$, then the
output is $X^{s+s'+2b,-\frac 1 2,1}$.
\begin{rem}
  Both last situations where $j''>k^\flat/2$ define the $R_2$ term in
  the section on existence when $s<\frac 1 4$. The gain depends on $b$
  (regularity $s+s'+2b$) and one may set $b'=0,\frac 1 4,\frac 1 2$,
  with respective outputs having $b''=-1,-\frac 3 4,-\frac 1 2$
\end{rem}
  \end{itemize}
\end{itemize}
\end{enumerate}
\item Finally, $k'=k^\sharp$: notice the conormal factors are
  identical whether $j''<<j$ or $j''\sim j$, and  we obtain exactly the same
  result as $k=k^\sharp$ when $j''<<j$.
\end{enumerate}

Proceed with the $++$ case: we have, by support considerations,
  $j^\flat\sim j^\sharp$, and $2j''\lesssim
  k^\sharp$ from the dispersion relation and we consider 3 cases.
  \begin{enumerate}
  \item $k''=k^\sharp$: assume $k,k'<< k''$, then $2j''\sim k''$,
$$
c_{j''k''}\lesssim\sum_{k,k'<k''} \inf( 2^{\frac {k^\flat} 2+\frac{j''} 2}, 2^{\frac
  {k^\flat} 2+\frac{k^\natural} 4}) \alpha_{j''k'}
2^{-j''s'-k'b'}  \beta_{j''k} 2^{-j''s-kb}.
$$
Assume $k<k'$,
\begin{itemize}
\item if $k'<2j''$, the conormal factor is better,
$$
c_{j'',k''}\lesssim 2^{-j''(s+s')} \sum_{k<k'<k''} 2^{\frac{k'} 4+\frac{k} 2-k'b'-kb}\alpha_{j'',k'} \beta_{j'',k}.
$$
Given $b=\frac 1 2$ and $b' =\frac 1 4$ (or the other way around), use
$\alpha,\beta\in l^1$, to get (trading regularity)
$$
c_{j''k''}\lesssim  2^{-j''(s+s'+1)+\frac{k''}2} \lambda_{j''k''},
$$
and the output is $X^{s+s'+1,-\frac 1
  2,1}\chi_{2j''\sim k''}$;
\begin{rem}
  This defines the $R_3$ term in the section on existence below
  $s=\frac 1 4$.
\end{rem}
\item if $2j''<k'$, Sobolev is better,
$$
c_{j'',k''}\lesssim 2^{-j''(s+s'-\frac 1 2)} \sum_{k<k',\,2j''<k'} 2^{\frac{k} 2-k'b'-kb}\alpha_{j'',k'} \beta_{j'',k}.
$$
Given $b=b'=\frac 1 2$ and discarding summability we to get (trading regularity)
$$
c_{j''k''}\lesssim  2^{-j''(s+s'+2b'+\frac 1 2-\e)+(\frac{1}2-\e)k''} \lambda_{j''k''},
$$
with $\lambda\in l^1$, and the output is $X^{s+s'+\frac 3 2-\e,-\frac
  1 2+\e,1}$.
\end{itemize}
On the other hand, assume $k<k'\sim k''$, then $2j''<< k''$ and
Sobolev is better, hence
$$
c_{j''k''}\lesssim 2^{-j''(s+s'-\frac 1 2)} 2^{-b'k''} \sum_{k<k''}
2^{(\frac 1 2-b)k} \alpha_{j'',k''}\beta_{j'',k},
$$
and with $\alpha,\beta\in l^\infty$ and $b=1/2$, we get
$X^{s+s'+2b'+\frac 1 2-4\e,-\frac 1 2+\e,1}\chi_{2j''<k''}$.
\item $k=k^\sharp$, we have $ 2j''\lesssim k$. We have
$$
c_{j'',k''}\lesssim \sum_{k''<k,k'<k} 2^{\frac{k^\flat} 2+\inf(\frac{j''}2,\frac{k^\natural-j''} 2) -j''(s+s')-kb-k'b'}\alpha_{j',k'}\beta_{j,k}.
$$
\begin{enumerate}
\item We restrict the sum to $k'<k''$:
  \begin{itemize}
  \item 
if $j''<k''/2$, Sobolev is better than conormal,
$$
c_{j'',k''}\lesssim 2^{-j''(s+s'-\frac 1 2)} \sum_{k>k''} 2^{-kb} \beta_{j'', k}\sum_{k'<k''} 2^{k'(\frac 1 2-b')}\alpha_{j'',k'}.
$$
Now the last sum yields $k''$, and we get with $b\geq \frac 1 2$,

$$
c_{j''k''}\lesssim  2^{-(s+s'-\frac 1 2) j''}2^{-bk''}k'',
$$
and again we may trade conormal regularity for spatial regularity, as $k''\geq 2j''$: 
$$
c_{j''k''}\lesssim 2^{(\frac 1 2-\epsilon) {k''}}
2^{-j''(s+s'+1)} k'' 2^{-k''(\frac 1 2+b-\epsilon)+\frac
  3 2 j''}.
$$
Therefore, the output is $X^{s+s'+2b+\frac 1 2-4\epsilon,-\frac 1
  2+\epsilon,1}\chi_{2j''\lesssim k''}$.
\item
if $j''>k'' /2$, conormal is better (and recall $k>> 2j''$),
$$
c_{j'',k''}\lesssim 2^{\frac{k''} 2} 2^{-j''(s+s'+\frac 1 2)} \sum_{k>2j''}
2^{-b k} \beta_{j'', k}\sum_{k'<k''} 2^{k'(\frac 1 2-b')}\alpha_{j'',k'}.
$$
The last sum is again $k''$, and
$$
c_{j''k''}\lesssim  2^{\frac{k''}{2}} 2^{-j''(s+s'+\frac 1 2+2b)} k'',
$$
therefore, the output is $X^{s+s'+2b+\frac 1 2-\e,-\frac 1
  2+\e,1}\chi_{k''<2j''}$.
\end{itemize}
\item Now we deal with $k''<k'<k$:
 \begin{itemize}
  \item if $j''<k'/2$, Sobolev is better,
$$
 c_{j'',k''}=2^{\frac{k''} 2} 2^{-j''(s+s'-\frac 1 2)}\sum_{k>\sup(k'',2j'')}
 2^{-kb} \beta_{j'', k} 
\sum_{\sup(k'', 2j'')<k'<k} 2^{-b'k} \alpha_{j'', k'}.
$$
and 
$$
c_{j''k''}\lesssim 2^{k''/2-(s+s'-\frac 1 2) j''-(b+b')\sup(k'',2j'')},
$$
which is in $X^{s+s'+\frac 3 2-4\epsilon,-\frac 1 2+\e,1}$ as
$b,b'\geq \frac 1 2$.
\item if $j''>k'/2$, conormal is better,
$$
 c_{j'',k''}=2^{\frac{k''}2} 2^{-j''(\frac 1 2+s+s')}\sum_{k>\sup(k'',2j'')}
 2^{-kb} \beta_{j'', k} 
\sum_{k''<k'<k} 2^{(\frac {1} 2-b')k} \alpha_{j'', k'}.
$$
Therefore
$$
c_{j''k''}\lesssim 2^{(\frac 1 2-\epsilon)k'''-(s+s'+\frac 1 2+2(b-2\epsilon))j''},
$$
and the output is $X^{s+s'+\frac 1 2+2(b- 2\epsilon)-\e,-\frac 1 2+\epsilon,1}$.
  \end{itemize}

\end{enumerate}

\item Finally, $k'=k^\sharp$: we have $2j''\lesssim k'$, and use the conormal factor.
$$
c_{j''k''}= 2^{\frac{k''}2-j''(\frac 1 2+s +s')} \sum_{\sup(k'',2j'')<k'}\alpha_{j'', k'} 
2^{-b'k'} \sum_{k<k''} 2^{k(\frac 1 2-b)}\beta_{ j'', k}.
$$
The last sum is (again) $k''$, and
$$
c_{j''k''}\lesssim 2^{k''(\frac{1}2-\epsilon)-j''(\frac 1 2+s +s'-2\epsilon+2b')},
$$
therefore, the output is $X^{\frac 1
  2+2b'+s+s'-4\epsilon,-\frac 1 2+\epsilon,1}\hookrightarrow
X^{\frac 3 2+s+s'-4\epsilon,-\frac 1 2+\epsilon,1}$.
\end{enumerate}

We can now collect all cases to obtain that for any $\epsilon<\frac 1
{10}$,
$$
R(u,v)\in X^{1+s+s',-\frac 1 2+\epsilon,1}.
$$
Therefore, we have obtained the two cases of interest:
\begin{itemize}
\item If $s=\frac 1 4$, $s'=0$ and $b=b'=\frac 1 2$, $
\partial_x R(u,v) \in X^{\frac 1 4,-\frac 1 2,1}.
$
\item If  $s>\frac 1 4$ and $b>\frac 1 2$, $s'=0$ and
  $b'=\frac 1 2$, one has $
\partial_x R(u,v) \in X^{s,-\frac 1 2+\epsilon,1}$, except for very
  specific cases where $j\sim j'\sim j''$. Let us postpone them: by interpolation,
  $X^{s_\theta,b,2}\subset[X^{s_1,b,1},X^{s_2,b,1}]_\theta$, we get
  the desired estimate, $
\partial_x R(X^{s,b,2},X^{0,\frac 1 2,1}) \mapsto X^{s,-\frac 1 2+\epsilon,2}.
$
On the remaining terms, they really require only $\beta\in l^2$
  provided $b'>\frac 1 4$ (to recover $b-\frac 1 2$ at the end).
\end{itemize}
\subsubsection{The spatial paraproduct: $j\sim j''\sim j^\sharp$}\label{sec.A.2.2}
Here, the value of $s$ should be irrelevant: depending on $s'$, we
might gain regularity ($s'=-3/4$, gain $\frac 1 4$), just be even
($s'=-1$), or have a loss ($s'<-1$). Hence the final regularity should be $s+1+s'=\sigma$.\\
>From the dispersion relation we have $j'+j''\lesssim k^\sharp$. We
split between values of $k^\sharp$ as
before, treating both $\pm+$ cases.
  \begin{enumerate}
  \item Case $k''=k^\sharp$. Note that $j'$ and
  $k^\natural-j^\natural=k^\natural-j''$ are not comparable a priori,
  thus Sobolev may be better.
  \begin{itemize}
  \item If $k<k'<<k''$, we have $j'\sim k''-j''$ and the conormal
  factor is always better, as $k''-j''>k'-j''$, hence
$$
c_{j''k''}\lesssim \sum_{k<k'<k''} 2^{\frac k 2+\frac{k'-j''}2}
2^{-j''s-(k''-j'')s'} 2^{-bk-b'k'} \alpha_{k''-j'', k'} \beta_{j'',k},
$$
for which summing over $k,k'$ yields
$$
c_{j''k''}\lesssim \sum_{k<k'<k''} 2^{\frac k 2+\frac{k'-j''}2}
2^{-j''s-(k''-j'')s'} 2^{-bk-b'k'}\alpha_{k''-j'', k'} \beta_{j'',k},
$$
discarding any summability in $k,k'$ to get
$$
c_{j''k''}\lesssim 2^{j''(-\frac{1} 2-s+s')- {k''} s'}k''^2 \lesssim 2^{j''(-\frac{1} 2-s+s')- {k''}( s'+\frac 1
  2)+\frac{k''} 2}k''^2,
$$
and given that $k''<2j''$, with $s'+\frac 1 2\leq 0$ we obtain 
$$
c_{j''k''}\lesssim 2^{j''(-\frac{3} 2-s-s'-3\epsilon)+k''(\frac{1}
  2-\epsilon)}\lambda_{j'',k''}, \text{ with } \lambda\in l^1.
$$
Therefore the output is $X^{\frac 3 2+s+s'-3\epsilon,-\frac 1
  2+\epsilon,1}\chi_{k''<2j''}$.
\item If $k<k'\sim k''$, then Sobolev is better as $j'<k''-j''$, and
$$
c_{j''k''}\lesssim \sum_{k<k''\,\,\,\,} \sum_{j'<\sup(j'',k''-j'')} 2^{\frac k
  2 +\frac {j'} 2-j''s-j's'-bk-b'k''} \alpha_{j'k''}\beta_{j''k},
$$
therefore, with $-s'+\frac 1 2> 0$, and discarding the $k$ sum ($b\geq 1/2$)
$$
c_{j''k''}\lesssim 2^{-j''s-bk''} \alpha_{\inf(j'',k''-j''),k''}
2^{(\frac 1 2-s')\inf(j'',k''-j'')}k''.
$$
\begin{itemize}
\item If $j''<k''-j''$, then
$$
c_{j''k''}\lesssim 2^{-j''(s+s'+\frac 3 2-4\epsilon)+k''(\frac 1
  2-\epsilon)} \lambda_{k'',j''},\text{ with } \lambda \in l^1,
$$
and the output is $X^{s+s'+\frac 3 2-4\epsilon,-\frac 1
  2+\epsilon,1}\chi_{2j''<k''}$.
\item If $k''<2j''$, then
$$
c_{j''k''}\lesssim 2^{-j''(s-s'+\frac 1 2)+\frac{k''}
  2-(b'+s')k''}\alpha_{k''-j'',k''}.
$$
If $s'+\frac 1 2=0$, one gets $X^{s+1,-\frac 1 2,1}\chi_{k''<2j''}$.

If $s'<-\frac 1 2 $ and discarding summability in $\alpha,\beta$, then
$$
c_{j''k''}\lesssim 2^{-j''(s+s'+\frac 3 2)+\frac{k''}
  2+k''((-s')-b')-j''((-2s')-1)}k'',
$$
and with $b'=\frac 1 2$ and $k''<2j''$, we get $X^{\frac 3 2+ s + s'-3\epsilon,-\frac 1 2+\epsilon,1}\chi_{k''<2j''}$.
\end{itemize}
\item If $k'<k<<k''$, then $j'\sim k''-j''$, pick the conormal factor,
$$
c_{j''k''}\lesssim\sum_{k'<k<k''} 2^{\frac {k'} 2+\frac{k-j''}2}
2^{-j''s-(k''-j'')s'} 2^{-bk-b'k'}\alpha_{k''-j'', k'} \beta_{j'',k},
$$
and, exactly as before, discarding any summability
$$
c_{j''k''}\lesssim 2^{-j''(s-s'+\frac 1
  2)}2^{-s'{k''}}{k''}^2,
$$
taking advantage of $k''<2j''$, we get $X^{s+s'+\frac 3
  2-3\epsilon,-\frac 1 2+\e,1}\chi_{k''<2j''}$.
\item If $k'<k\sim k''$, then $j'<k''-j''$, Sobolev is better, and we
  have
$$
c_{j''k''}\lesssim \sum_{k'<k''} \sum_{j'<\inf(j'',k''-j'')}
  2^{\frac{k'}2+\frac{j'}2-j''s
    -j's'-bk-bk''}\alpha_{j',k'}\beta_{j''k''},
$$
which is even slightly better than the $k<k'\sim k''$.
  \end{itemize}
\item Case $k=k^\sharp$: we have
$$
c_{j''k''}=\sum_{k',k''<k,}\sum_{\,j'<k-j''} 2^{\frac{k^\flat} 2+\inf(\frac{j'}2,\frac{k^\natural-j''} 2)}
2^{-j''s-j's'-kb-k'b'}\alpha_{j', k'} \beta_{j'',k}.
$$
\begin{itemize}
\item If $k'<k''<<k$, then $j'\sim k-j''$, pick the conormal factor,
$$
c_{j''k''}\lesssim 2^{\frac{k''} 2-j''(s-s'+\frac 1
  2)}\sum_{k'<k''<k<2j''} 2^{-ks'-kb} 2^{\frac {k'} 2-b'k'}\alpha_{k-j'',k'}\beta_{j'',k}.
$$
The sum over $k'$ is irrelevant,  assuming again $s'+b\leq 0$,
we get, with no summability over $\alpha,\beta$
$$
c_{j''k''}\lesssim 2^{\frac{k''} 2-j''(s+s'+\frac 1 2+2b)}{2j''}k'',
$$
which means the output is $X^{s+s'+2b+\frac 1 2-3\epsilon,-\frac 1
  2+\epsilon,1}\chi_{k''<2j''}$.
\item If $k'<k''\sim k$, we are back in a case where $k^\sharp\sim k''$.
\item If $k''<k'<< k$, then $j'\sim k-j''$, conormal again,
$$
c_{j''k''}\lesssim 2^{\frac{k''} 2-j''(s-s'+\frac 1 2)}\sum_{k''<k'<k<2j''}
2^{-ks'-kb-k'(b'-\frac 1 2)}\alpha_{k-j'',k'}\beta_{j'',k}.
$$
The sum over $k'$ is irrelevant, and recalling $s'+b\leq0$, we are back to the previous case,
$$
c_{j''k''}\lesssim 2^{(\frac{1} 2-\epsilon)-j''(s+s'+\frac 1 2+2b-3\epsilon)}\lambda_{j'',k''},
$$
which is  $X^{s+s'+2b+\frac 1 2-3\epsilon,-\frac 1
  2+\epsilon,1}\chi_{k''<2j''}$.
\item If $k''<k'\sim k$,and $j'\lesssim k-j''$, so that $j''<k$ and
  Sobolev is better, hence
$$
c_{j''k''}\lesssim \sum_{\sup(k'',j'')<k,}\sum_{\,j'<\inf(\frac k 2,k-j'')}
2^{\frac{k''}2+\frac{j'}2-j''s-j's'-k(b+b')}
\alpha_{j',k}\beta_{j,k},
$$
discarding summability, setting $b=b'=1/2$,
$$
c_{j'',k''}\lesssim 2^{\frac{k''}2-j''s} 2^{-\sup(k'',j'')(\frac 3
  4-\frac{s'}2)}\lesssim 2^{\frac{k''}/2-j''(s+s'+\frac 3
  2)+j''(\frac{3s'}2+\frac 3 4)},
$$
which, for $s'\leq - \frac 1 2$ is $X^{s+s'+\frac 3 2-4\epsilon,-\frac 1
  2+\epsilon,1}$.
\end{itemize}

\item Finally, $k'=k^\sharp$.
  \begin{itemize}
\item If $k<k''\sim k'$ we are back to a $k^\sharp\sim k''$ case.
\item If $k''<k\sim k'$ we are back to
  the very last case of $k^\sharp=k'$.
  \item If $k,k''<<k'$, then $j'+j''\sim k'$, hence $j''\lesssim
  k'<2j''$ (and therefore $k,k''<2j''$ as well). We have
$$
c_{j''k''}\lesssim \sum_{k'',k<k'}
2^{\frac{k^\flat}2+\min(\frac{j'}2,\frac{k^\natural-j'} 2)-k's'-bk-b'k'} 2^{-j''(s-s')}\alpha_{k''-j'', k'} \beta_{j'',k}.
$$
 Discard summability and assume $b=b'=\frac 1 2$, $s'\leq -\frac 1 2$,
\begin{itemize}
\item Let $k''<k<<k'$, then $k'\sim j'+j''$,
  \begin{itemize}
  \item if $j'=k'-j''<k-j'$
$$
c_{j''k''}\lesssim 2^{\frac{k''}2-j''(s-s'+\frac 1
  2)}\sum_{k''<k\,\,}\sum_{k<k'<j''+\frac k 2}
2^{-s'k'-\frac k 2} \lesssim  2^{\frac{k''}2-j''(s+s'+\frac 3
  2)}j'',
$$
which yields $X^{s+s'+\frac 3 2-4\e,-\frac 1 2+\e,1}$;
\item if $k-j'<k'-j''$ we get, if $s'\leq -1$,
\begin{multline}
c_{j''k''}\lesssim
2^{\frac{k''}2-j''(s-s')}\sum_{k''<k<2j''\,\,}\sum_{j''+\frac k 2<k'<2j''}
2^{\frac{j''}2-k'(1+s')}\\
\lesssim 2^{\frac{k''}2-j''(s-s'-\frac 1 2)-2j''(1+s')}j'',
\end{multline}
and again the output is  $X^{s+s'+\frac 3 2-4\e,-\frac 1 2+\e,1}$,
while for $-1<s'\leq -\frac 1 2$, retaining the summability  $\beta\in l^1$,
\begin{eqnarray*}
c_{j''k''}  \lesssim &
2^{\frac{k''}2-j''(s-s')}\sum_{k''<k<2j''\,\,}\sum_{j''+\frac k 2<k'<2j''}
2^{\frac{j''}2-k'(1+s')} \beta_{j'',k}\\
  \lesssim &
2^{\frac{k''}2-j''(s-s'-\frac 1 2)-(j''+\frac
  {k''}2)(1+s')}\lambda_j,\,\,\text{ with }\,\lambda\in l^1,\\
 \lesssim & 2^{\frac{k''}4-j''(s+s'+1)-(j-\frac
  {k''}2)((-s')-\frac 1 2)}\lambda_j,
\end{eqnarray*}
which is $X^{s+s'+1,-\frac 1 4,1}$.
  \end{itemize}
\item Let $k<k''<<k'$, then $k'\sim j'+j''$,
  \begin{itemize}
  \item if $k'-j''<k''-j'$,
\begin{multline}
c_{j''k''}\lesssim  2^{-j''(s-s'+\frac 1 2)}
\sum_{k<k''\,\,}\sum_{k''<k'<j''+\frac{k''}2} 2^{k'(-s')}\\
\lesssim
2^{-j''(s+s'\frac 3 2)+\frac{k''} 2+(s'+1)j''-(s'+1)\frac{k''}2} k'',
\end{multline}
which is again  $X^{s+s'+\frac 3 2-4\e,-\frac 1 2+\e,1}$;
\item if $k''-j'<k'-j''$, either $s'\leq -1$,
$$
c_{j''k''}\lesssim  2^{-j''(s-s'-\frac 1 2)+\frac{k''}2}
\sum_{k<k''\,\,}\sum_{j''+\frac{k''}2<k'<2j''} 2^{k'(-s'-1)}\lesssim
2^{-j''(s+s'\frac 3 2)+\frac{k''} 2} k'',
$$
and we obtain   $X^{s+s'+\frac 3 2-4\e,-\frac 1 2+\e,1}$; or
$-1<s'\leq -\frac 1 2$ and  retaining the summability  $\beta\in l^1$,
\begin{eqnarray*}
c_{j''k''}  \lesssim &
2^{\frac{k''}2-j''(s-s')}\sum_{k<k''\,\,}\sum_{j''+\frac {k''} 2<k'<2j''}
2^{\frac{j''}2-k'(1+s')} \beta_{j'',k}\\
  \lesssim &
2^{\frac{k''}2-j''(s-s'-\frac 1 2)-(j''+\frac
  {k''}2)(1+s')}\lambda_j,\,\,\text{ with }\,\lambda\in l^1,\\
 \lesssim &  2^{\frac{k''}4-j''(s+s'+1)-(j-\frac
  {k''}2)((-s')-\frac 1 2)}\lambda_j,
\end{eqnarray*}
which is again $X^{s+s'+1,-\frac 1 4,1}$.
\begin{rem}
On  both this term and the previously similar one, one may notice that
the output is such that $b''=-\frac 1 4$. As such, applying a gauge
transform will result in a lesser loss ($-\frac 1 8$ in general,
slightly more in the uniqueness part below $1/4$).
\end{rem}
  \end{itemize}
\end{itemize}
\end{itemize}
\end{enumerate}
Finally, collecting all terms, we got all cases we are interested in.


\end{document}